\documentclass[hidelinks,onefignum,onetabnum]{siamart250211}

\usepackage{lipsum}
\usepackage{amsfonts}
\usepackage{graphicx,epstopdf} 
\usepackage[caption=false]{subfig} 
\usepackage{algorithmic}

\ifpdf
  \DeclareGraphicsExtensions{.eps,.pdf,.png,.jpg}
\else
  \DeclareGraphicsExtensions{.eps}
\fi

\numberwithin{equation}{section}
\newtheorem{remark}{Remark}[section]
\newtheorem{ex}{Example}[section]
\newtheorem{prob}{Problem}[section]

\def\dis{\noindent{\bf Discussion: }}

\def\cO{{\cal O}}
\def\sm{\\smallskip}
\def\sm{\smallskip}
\def\ms{\medskip}
\def\bs{\bigskip}
\def\noin{\noindent}
\def\ts{\thinspace}
\def\nb{{n_b}}
\def\ni{{n_i}}
\def\Hw{{\tilde H}}
\def\Gw{{\tilde G}}
\def\rw{{\tilde r}}
\def\ax{{a_x}} \def\ay{{a_y}} \def\az{{a_z}}
\def\bx{{b_x}} \def\by{{b_y}} \def\bz{{b_z}}
\def\kx{{k_x}} \def\ky{{k_y}} \def\kz{{k_z}}
\def\dx{{d_x}} \def\dy{{d_y}} \def\dz{{d_z}}
  
\def\nx{{n_x}} \def\ny{{n_y}} \def\nz{{n_z}}
\def\nt{{n_t}}

\def\Nbar{{\bar N}}
\def\Nhat{{\widehat N}}
 
\def\RR{\mathop{{\rm I}\kern-.2em{\rm R}}\nolimits}
\def\pO{{\partial \Omega}}
\def\tri{{\triangle}}
\def\cS{{\cal S}}
\def\tp{{tensor-product}}
\def\tpart{{tetrahedral partition}}

\def\frl#1{Fig.\,\ref{#1}(left)}%
\def\frr#1{Fig.\,\ref{#1}(right)}%
\def\Szd{{{\cal S}^0_d(\tri)}}
\def\Pd{{\cal P}_d}
\def\st{{\; : \;\,}}
\def\sb#1{\setbox0=\hbox {\hskip #1in}
\vbox\bgroup \def\xem{.7em}\hrule
  \hbox\bgroup\vrule\kern\xem \vbox\bgroup \kern\xem}
\def\eb{\kern\xem\egroup \kern\xem\vrule\egroup \hrule \egroup}
\def\ob{\hsize\wd0 \obeylines \parskip0pt \parindent0pt \tt}
\def\eword#1{{\tt #1}}

\headers{Solving
PDEs on Curved Domains in 3D with splines}{A. Greene, L. L. Schumaker}

\title{Solving PDEs on Curved Domains in 3D using the Immersed Penalized
Boundary Method with Trivariate Spline Spaces}

\author{Aussie Greene\thanks{Department of Mathematics, Vanderbilt University,
Nashville, TN \\
  (\email{aussie.m.greene@vanderbilt.edu}, \email{larry.schumaker@vanderbilt.edu}).}
\and Larry L. Schumaker\footnotemark[1]}

\usepackage{amsopn}

\nocite{soft}
\nocite{Holcombe12}

\begin{document}

\maketitle


\begin{abstract}
Second-order elliptic boundary-value problems (BVPs) defined on curved domains in 2D and 3D
arise frequently in practice. Considerable work has gone into creating numerical
methods to solve such problems.
One of the newest and most promising methods is the
immersed penalized boundary method (IPBM) introduced in [Schumaker, L. L., Solving
elliptic PDE’s on domains with curved boundaries with an immersed penalized boundary
method, J. Sci. Comp. 80(3) (2019), 1369–1394]. In this paper we present a software
package that extends the use of this method to complex curved domains in 3D by using
trivariate spline spaces. The package works with the well-known spaces of trivariate
tensor-product splines, as well as with spaces of polynomial splines on tetrahedral
partitions. In addition to describing functions and scripts in the package, we run
the software on a detailed set of examples  for a variety of BVPs defined on 
different shaped curved domains in $\RR^3$. We believe these examples show
that due to its high accuracy, low
computational cost, and ease of implementation,
IPBM methods  are an attractive alternative to existing methods.
\end{abstract}



\section{Introduction}
Suppose $\Omega$ is a domain in 3D which has a curved boundary.
In this paper we focus on the numerical solution of 
the following boundary-value problem for 
second order elliptic PDEs defined on such domains. Let
\begin{equation}   
Lu = a_1 u_{xx} + a_2 u_{yy} + a_3 u_{zz}   + 2a_4 u_{xy}
  + 2a_5u_{xz} + 2a_6u_{yz}, \label{L1} \end{equation}
where $a_1,\ldots,a_6$ are functions defined on $\Omega$.
The best known example of this type of operator is the
Laplace operator $Lu = \tri u: = u_{xx} + u_{yy} + u_{zz}$.

\begin{prob}  \label{Prob1}
Suppose $\Omega$ is a domain in $\RR^3$ and that $L$ is a linear second order partial
differential operator.    Given a function $f$ defined on $\Omega$ and
a function $g$ defined on the boundary $\partial\Omega$
of $\Omega$, find a function $u$ defined on $\Omega$ such that
\begin{align}
   Lu &= f, \hbox{\qquad on $\Omega$,} \label{Lu} \\
   u &= g, \hbox{\qquad on $\pO$}.  \label{Dir}
\end{align}
\end{prob}

Constructing numerical methods to efficiently and accurately solve this problem is
of great importance to engineers and scientists, and a wide variety of 
computational methods
have been suggested.  While much of the early work dealt with  domains with
a polygonal boundary, a lot of the recent work has focused on the considerably more 
difficult case where the boundary 
is curved or the domain has cavities and/or holes.
We do not attempt to review this literature here,
but see Remarks~\ref{immerse}-- \ref{RBF}  where we list
some papers and books that are related to our work. These resources 
provide extensive lists of additional references to the literature.

In this paper we discuss the implementation of two methods
introduced in \cite{IPBM} which solve the problem by immersing
the domain $\Omega$ in a larger computational domain $D$ and setting
up an appropriate least-squares problem with a quadratic functional which includes
a penalty term  involving
collocation of the boundary function $g$  at points on the boundary of $\Omega$.
To make use of this approach, we have to assume that the coefficients of $L$
and the function $f$ 
in (\ref{Lu}) are defined not only on $\Omega$ but on all of $D$,
see Remark~\ref{extension}. 
We believe this approach has
significant advantages compared with other methods in the literature,
see the discussion in Section~9.
Although immersing domains have been used in some earlier papers,
see Remark~\ref{immerse}, the novelty
of the methods studied here is the way in which the boundary conditions 
are dealt with.

It was shown in \cite{Scompb} that
these  methods are highly effective in solving
second order elliptic boundary values problems on curved domains in $\RR^2$.
Our aim in this paper is to  provide a  self-contained Matlab package
for using IPBM on curved 3D domains together with two types of
well-known spaces of multivariate splines. We will use it on
a variety of examples to show that the methods are also
highly effective in the 3D case.
The package provides
everything needed to run all of the examples in the paper, and
includes scripts  that the reader can use
to compare with their favorite methods or to 
solve problems with PDEs and domains of their own choosing.

The paper is organized as follows.
In Section~2  we describe our methods mathematically, while
in Section~3 implementation issues are addressed.
In Sections~4--5 we
illustrate the use of our methods with trivariate tensor-product splines, while
in Sections~6--7 we do the same for polynomial splines defined on
tetrahedral partitions of an immersing domain.
Section~8 is devoted to an approach to reducing computational time when working
with splines on tetrahedral partitions which is based on
trimming the partitions. Finally,
in Sections~9 and 10  we present conclusions and remarks.

\section{Two immersed penalized boundary methods in 3D}
Given a 3D curved domain $\Omega$ contained in an immersing domain $D$,
let $\cS$ be a finite dimensional linear space of approximating 
functions defined on $D$.
Suppose $\{\phi_j\}_{j=1}^n$ is a basis for $\cS$. Then every $s \in \cS$
can be written 
in the form 
\begin{equation}
 s = \sum_{j=1}^n c_j \phi_j. \label{s} 
\end{equation}
To find a coefficient vector $c$
so that $s$ is a good approximation to the solution of Problem~\ref{Prob1},
we can apply either of the following two methods introduced in \cite{IPBM},
see also \cite{Scompb}. We call the first one IPBF since it is inspired by
the finite-element method, and the second one IPBC since it is a
least-squares collocation method.

\subsection{The IPBF method}
Given a set of points
B:=$\{\xi_i,\eta_i,\zeta_i\}_{i=1}^\nb$ on $\pO$, and a parameter $\lambda >0$,
find $s\in \cS$ that minimizes
\begin{equation}
  \Phi_{F}(s) :=  \sum_{i=1}^n \big[\int_D  (L s - f) \phi_i\big]^2 
   + \lambda^2\sum_{i=1}^\nb  [s(\xi_i,\eta_i,\zeta_i) - g(\xi_i,\eta_i,\zeta_i)]^2.
     \label{IPBF}  
\end{equation}
Finding a minimum of this functional is equivalent to solving the following
system of equations in the least-squares sense:

\begin{align}
 \int_D L s\ts \phi_i &=\int_D f\ts \phi_i, 
 \qquad{i=1,\ldots,n},  \label{PBFa} \\
    \lambda\ts s(\xi_i,\eta_i,\zeta_i)  &= \lambda \ts g(\xi_i,\eta_i,\zeta_i), 
    \qquad{i = 1,\ldots,\nb}. \label{PBFb}
\end{align}
\goodbreak

Substituting (\ref{s}) in these equations 
leads to an overdetermined  linear system of equations 
of the form $Hc = r$  for the vector $c$ of coefficients of $s$.
This system can be solved in the least-squares sense by simply using the
Matlab backslash operator. When $H$ is of full rank,
it can also be solved by multiplying both sides by $H^T$ to get the equivalent
system $Gc = H^Tr$, where $G = H^TH$ is a 
symmetric and positive--definite matrix. Throughout this paper we will use
the backslash operator, see Remark~\ref{backslash} for some reasons why.
Note that the integrals appearing here are defined on the
domain $D$, not on $\Omega$.  The parameter $\lambda$ 
controls the relative weight given to the two sets of equations.

\subsection{The IPBC method}
Given sets of points
$\Gamma := \{\alpha_i,\beta_i,\gamma_i\}_{i=1}^\ni$ 
in $D$ and $B:=\{\xi_i,\eta_i,\zeta_i\}_{i=1}^\nb$ 
on $\pO$ and a real number $\lambda > 0$, find $s\in \cS$ that minimizes
\begin{equation}
  \Phi_{C}(s) :=  \sum_{i=1}^\ni 
[L s(\alpha_i,\beta_i,\gamma_i) - f(\alpha_i,\beta_i,\gamma_i)]^2 
   + \lambda^2\sum_{i=1}^\nb  [s(\xi_i,\eta_i,\zeta_i) - g(\xi_i,\eta_i,\zeta_i)]^2.
     \label{IPBC}  
\end{equation}
Finding a minimum of this functional is equivalent to solving the following
set of equations in the least-squares sense
\begin{align}
L s(\alpha_i,\beta_i,\gamma_i) &= f(\alpha_i,\beta_i,\gamma_i),  \qquad i =
1,\ldots,\ni, \label{IPBCa}  \\
  \lambda\ts s(\xi_i,\eta_i,\zeta_i)  &= \lambda \ts g(\xi_i,\eta_i,\zeta_i),
    \qquad{i = 1,\ldots,\nb}. \label{IPBCb}
\end{align}

Substituting (\ref{s}) in these equations
leads to a linear system  $\Hw c = \rw$  for the 
vector $c$ of coefficients of $s$. This system
can be solved in the least-squares sense using the Matlab backslash operator.
Assuming $n_i + n_b \ge n$
and that  $\Hw$ is of full rank, we can also convert it to the 
equivalent  system $\Gw c = \Hw^T\rw$, where $\Gw = \Hw^T\Hw$ is symmetric and positive--definite.
Throughout the paper we will work with $\Hw$ directly, see Remark~\ref{backslash}.

\section{Implementing our methods}

\subsection{Choosing an approximation space}
In this paper we are going to
take the immersing domain $D$ to be the bounding box associated with the domain
$\Omega$, i.e., the smallest rectangular box in $\RR^3$ that contains $\Omega$.
Although the above methods can be used with any finite dimensional linear space 
$\cS$ of functions on $D$, here we will focus on spaces of splines
(piecewise polynomials) defined on a partition of $D$. First we will consider
trivariate \tp{} splines where the partition of $\Omega$
consists of a collection of sub-boxes, see \cite{Scompb} and Section~4 below.
We call this a box partition.
As an alternative, we also discuss the use of polynomial splines defined 
on so-called type-5 tetrahedral partitions which are obtained by splitting
each sub-box in a box partition into five tetrahedra, see
Section~6.  Both types of splines have excellent approximation properties,
see Remarks~\ref{appxTP} and \ref{Appxtet}. Moreover,
a comprehensive set of Matlab software is readily available 
for computing with them, see \cite{Scompb}.

\subsection{Setting up the domain}

For our purposes a domain in 3D is 
a solid object consisting of a set of points that are enclosed in a watertight
piecewise smooth boundary.  We allow domains with cavities and with
holes passing through them.
The domains we are interested in may be
actual physical objects, or they may be defined mathematically, for
example by using some kind of CAD/CAM software.
It will be enough to have a finite point cloud of reasonably well-spaced points
lying on the boundaries of the domain. For physical objects we can create
such clouds using a scanner.  For mathematically defined surfaces,  
we can make use of existing software to create the clouds. 

Point clouds can be stored in various formats.
A commonly used choice is to work with 
a file giving the number of points and a list of the 
Cartesian coordinates of each point.  We use the extension .dat for such files.
For an example, see the file \eword{sphere50.dat} which is used in 
Example~\ref{Exsphere}. It  contains 50,000 points on the surface of the unit sphere.

Point clouds are also sometimes stored in what is called
the PLY format. 
In addition to a list of coordinates, such files typically contain
additional information that we don't need.  An example is the
file {\tt bunny.ply} which contains 35,947 points on the surface of
the classical Stanford bunny, see Example~\ref{Exbunny}.
Here is a function that prompts for and reads such files.

\ms 
\sb{1.7}\vbox{\ob
[xp,yp,zp] = setdomain3
}\eb
\ms\noin 
This function can read files in both the  .dat or .ply format, but also
files with the extension .stl.  Such files are called 
STL (stereo lithography) files. They describe the boundary of the
object as  a collection of (flat) triangles, and  
are particularly useful for 3D printing. We should note, however, that
STL files do not provide exact representations of objects with truly curved
boundaries. Matlab and its various tool boxes contain numerous
examples of 3D objects defined by STL files (often in binary form).
The ones used here
were selected primarily because we believe they define domains that might typically
arise in applications.

Matlab also includes a variety of tools for processing STL files.
Given a domain stored in STL format, to
form the set $B$ needed for our IPBM methods, we use
the Matlab function {\tt mesh2pc} (based on Poisson disk sampling)
to create an associated point cloud.

Another CAD/CAM approach to modelling 3D surfaces is to use
a collection of NURBS patches.  This creates objects with 
truly curved  boundaries.
This can be challenging for complicated domains, 
particularly if the domain has holes or cavities.
For example, the well-known NURBS model of the Utah teapot 
shown in \frl{Fig3} makes use of 32 individual NURBS patches.  
Several CAD/CAM systems include software for turning NURBS surfaces
into point clouds.  The above function reads these point clouds, and is not
designed to read NURBS surfaces directly.

\subsection{Choosing the set B}
To use our IPBM methods, we need to select a set 
B = (xb,yb,zb) of points on the boundary of
the domain to be used for fitting the Dirichlet boundary conditions.
We like to begin with point sets with at least 50,000 points in order
to make sure they do a good job of capturing the geometry of our domain.
We could choose B to be this original point cloud, but we can save time
in later computations by downsampling it first.  This can also help with
getting well-spaced points.  

The problem of downsampling a point cloud is a major research topic in
its own right. For our implementation
we employ the farthest point strategy, see \cite{Eldar}.
To get $\nb$ points, this algorithm picks a starting
point on the cloud, and then repeatedly adds a point on the
cloud that maximizes the minimum distance to all previously accepted points.
The iteration stops  when we have $\nb$ points.

Here is the function that we use to downsample the point cloud describing
the boundary of a domain to produce smaller well-spaced sets of
points on the boundary.

\ms
\sb{2.4}\vbox{\ob
[xb,yb,zb] = pcdown(xp,yp,zp,nb)
}\eb
\ms\noin
Here \eword{(xp,yp,zp)} are the vectors giving the point cloud representation of the
boundary, and \eword{nb} is the desired number of points
to be created. These points are selected from the points in the cloud
using the farthest points strategy to ensure that they are well-spaced.
In the supplement we include several examples to explore how the choice
of {\tt nb} affects the accuracy of our methods.
In the examples  presented here we typically choose 
{\tt nb} to be between 2,000 and 5,000.

\subsection{Computing dense sets of test points for error computations}
All of the examples in this paper solve BVPs where we know the true solution.
In that case it is generally of interest to compute a measure of how well
the spline solution $s$ fits the true solution. To do this we
create a large set of well-spaced points (xr,yr,zr) lying inside the domain,
and then compute the error vector u(xr,yr,zr) - s(xr,yr,zr). We then use
it to compute the maximum error over these points, and the corresponding RMS error.
Here is our function for generating such a set of points.

\ms
\sb{2.3}\vbox{\ob
[xr,yr,zr] = inpts3(xp,yp,zp,n)
}\eb
\ms\noin
The main inputs to this function are the vectors (xp,yp,zp) giving the point-cloud
representation of the boundary.  To create the desired points, this function
first constructs an  $n \times n \times n$ grid covering the bounding box surrounding
the domain, and then selects those points on this grid which cannot be seen
from a set of six sensor points located outside the boundary box, each
corresponding to a different  side of the box.
This code, which we include in our package,  makes use of the so-called
{\sl Hidden-Point-Removal operator}, see \cite{Katz}.  It is
based on the idea of identifying and removing all points in a domain that can be
seen from a given collection of view-points lying outside of $\Omega$.
This approach does not always
remove all of the grid points that lie outside of the domain,
particularly if the domain has holes.  To get a set of points
all of which lie inside $\Omega$ we may have to add and adjust sensor locations.
However, having a few points outside the domain does
not have much of an effect on the errors reported in our tables, see
Sect.~SM8 of the supplement.
\frr{Fig4} shows the interior points chosen by our algorithm
that lie in a slice through the vessel in \frl{Fig4}.
The interior points chosen for error calculations are shown in red.

\subsection{Setting up the PDE}
Throughout this paper we will be giving examples 
of the performance of IPBM methods
for various test BVPs. To simplify our scripts, 
we have written  the following function to input 
an exact solution $u$ and the corresponding functions $f$ and $g$ assuming
that $L$ is the Laplace operator.

\ms
\sb{1.8}\vbox{\ob
 [u,f,g] = setpde3L(pde)
}\eb
\ms\noin
The input to this function is a positive integer 
\eword{pde} less than or equal to ten.
The output includes \eword{u}, the corresponding
right-hand side \eword{f} of the PDE, and the function \eword{g}
giving the Dirichlet boundary conditions.
The table below lists the solutions corresponding to
\eword{pde} from one to ten.

For problems where $L$ is written in the more general
form (\ref{L1}), we have included the function

\ms
\sb{3.0}\vbox{\ob
 [u,f,g,a1,a2,a3,a4,a5,a6] = setpde3(pde)
}\eb
\ms\noin
For values of \eword{pde} up to 10, this function gives the same true solutions as
\eword{setpde3L}, and sets the coefficients of the differential operator $L$
to be (1,1,1,0,0,0) which
corresponds to the Laplace operator.  For {\tt pde} between 11 and 12, it
sets up BVPs which correspond to some PDEs with variable coefficients.
The true solutions are the same as for {\tt pde} = 4 and 6, respectively.
The second table below shows the values of the coefficients $a_1,\ldots,a_6$
for these values of {\tt pde}.

\bs
\sb{1.9}\vbox{\ob 
1) $x+2y + 3z$
2) $x^3 + 2y^3 + 3z^3$
3) $x^5 + 2y^5 + 3z^5$  
4) $\sin(5x) + \sin(5y) + \sin(5z)$
5) $\sin(5x)\sin(5y)\sin(5z)$
6) $\sin(8x)+\sin(12y)+\sin(14z)$
7) $|x+y-z|^3$
8) $40(x^2-x)(y^2-y)(z^2-z)$
9) $xy^2z^3$
10) $\sin(8x)\sin(12y)\sin(14z)$
}\eb
\ms\noin

\bs
\sb{3.5}\vbox{\ob
11) $\#4 :(10+\cos(x),e^y,10+\sin(x),\cos(x),\cos(y),\cos(z))$ 
12) $\#6 :(10+\cos(x),e^y,10+\sin(x),0,0,0)$
}\eb
\ms\noin

\section{Using trivariate \tp{} splines}
To use our immersed penalized boundary  methods with trivariate \tp{} splines,
we choose the immersing domain $D$ to be the bounding box surrounding
the domain $\Omega$. 
For details on computing with trivariate \tp{} splines, along with
associated Matlab software, see
Chapter~5 of \cite{Scompb}. Here we follow the notation
of Section~15.2 of that book which we briefly review here.

Let
$D:=  [a_x,b_x] \otimes [a_y,b_y] \otimes [a_z,b_z]$
be a rectangular box in $\RR^3$.  
Given positive integers $\kx,\ky,\kz$, let
\begin{align}
 \tri_x &:= \{\ax = x_0 < x_1 < \cdots < x_{\kx+1} = \bx\}, \label{kx} \\
\tri_y &:= \{\ay = y_0 < y_1 < \cdots < y_{\ky+1}  = \by\}, \label{ky} \\
\tri_z &:= \{\az = z_0 < z_1 < \cdots < z_{\kz+1}  = \bz\}\label{kz} . 
\end{align}
Given a positive integer $\dx$ let $\cS^{\dx-1}_\dx(\tri_x)$
be the associated space of univariate polynomial splines of degree $\dx$.
It is of dimension $\nx = \kx + \dx + 1$, and is spanned by the
normalized B-splines $\{N^{\dx+1}_i\}_{i=1}^\nx$.
Similarly, for given $\dy$ and $\dz$, let
$\{\Nbar^{\dy+1}_i\}_{i=1}^\ny$ be the B-spline basis for
$\cS^{\dy-1}_\dy(\tri_y)$ and
$\{\Nhat^{\dz+1}_i\}_{i=1}^\nz$ be the B-spline basis for
$\cS^{\dz-1}_\dz(\tri_z)$. Now to get a space $\cS$ of trivariate tensor-product
splines, we take the span of the products of these univariate B-splines in each of the
three variables. It is
shown in  Lemma~5.2 of \cite{Scompb} that these trivariate basis functions are
linearly independent, and thus
every spline $s \in \cS$ can be written
uniquely in the form
 $$  s(x,y,z) = \sum_{i=1}^\nx \sum_{j=1}^\ny \sum_{k=1}^\nz c_{ijk} 
     N^{\dx+1}_i(x) \Nbar^{\dy+1}_j(y) \Nhat^{\dz+1}_k(z) \label{bform}$$
for some matrix
$C := [c_{ijk}]_{i=1,j=1,k=1}^{\nx,\ny,\nz}$  of coefficients.

\subsection{Implementing IPBF with trivariate \tp{} splines}
To set up the equations (\ref{PBFa}) when using \tp{} splines, we need to calculate
integrals of \tp{} basis functions and their derivatives
over the box $D$.  This can be done by accumulating
the integrals over the sub-boxes of $D$ defined by the knots (\ref{kx}) -- (\ref{kz}).
For this we can use tensor-products of standard univariate Gaussian quadrature rules,
see Section~1.12.1 of \cite{Scomp}.
Both the IPBF and IPBC methods are coded using \tp{} splines with
arbitrarily spaced knots, but in all of the examples here
we have chosen equally-spaced
knots.  Here is our function for computing the coefficient
matrix \eword{C} of such a \tp{} spline.

\ms 
\sb{4.0}\vbox{\ob
[C,H,t1,t2] = ipbftp3(dx,dy,dz,xe,ye,ze,f,xb,yb,zb,...
~~~    a1,a2,a3,a4,a5,a6,gb,lambda)
}\eb
\ms\noin Here \eword{dx,dy,dz} are the degrees of the \tp{} splines,
and \eword{xe,ye,ze} are the extended knot sequences
in each of the three variables. The vectors \eword{xb,yb,zb} contain
Cartesian coordinates of the points in the set $B$, and \eword{gb} contains
the values of the function $g$ at these points. The input \eword{lambda} is
a nonnegative real number, and the other inputs describe the PDE.
For diagnostic purposes, we also output 
the matrix {\tt H} of the system being solved.

We have also written a simplified version of this function that can be used when the
PDE is the Laplace operator.

\ms 
\sb{4.6}\vbox{\ob
[C,H,t1,t2] = ipbftp3L(dx,dy,dz,xe,ye,ze,f,xb,yb,zb,gb,lambda)
}\eb
\ms\noin 
The inputs here are the same as in the previous function, minus the coefficients
{\tt a1,..,a6}. The output {\tt H} is  the matrix of the
linear system used to compute the coefficients of the {\tp} spline.
We have explored how to choose the parameter $\lambda$ when using these two
functions, and have concluded that $\lambda = .01$ is a good default value,
see the supplement.
To test these functions we have written the scripts \eword{TP3ipbf} and
\eword{TP3ipbfL}. To explore rates of convergence of these methods we
can use  \eword{TP3ipbfr} and \eword{TP3ipbfLr}.

\subsection{Implementing IPBC with trivariate \tp{} splines}
The implementation of the IPBC method based on the spaces of \tp{} splines introduced
in the previous section is nearly the same as for IPBF. In particular,
we can choose the set $B$ on the boundary
of the domain and the parameter $\lambda$ in the same way as above.
But now we also have  to choose a set $\Gamma$ of collocation points.
Throughout this paper we take $\Gamma$ to be the union of sets
of  $m_c^3$ points lying on equally-spaced grids 
covering each sub-box of the \tp{} partition. 
Here $m_c$ should be chosen large enough to make the matrix $H$ have full rank.
To get a rule of thumb to achieve this, we can compare the
number of collocation points to the dimension of the \tp{} spline space being used.
In the examples below for simplicity we will take {dx = dy = dz = d} and use $m$ equally-spaced
knots in each direction. Under these assumptions, the dimension of the \tp{} spline
space is $(d+m-1)^3$ while the number of points in $\Gamma$ is approximately
$[(m-1)(mc-1)]^3$ which leads to the condition $m_c \ge 1+(d+m-1)/(m-1)$. 
Experiments in the supplement show that as we increase $m_c$, we get larger
computational times without too much change in condition numbers or accuracy.
Since the examples below use only $d=4$ or $d = 5$ and $m \ge 5$, this suggests
that $m_c = 3$ is a good default value. We will use it unless otherwise stated.

Given functions $f$ and $g$ as in
Problem~\ref{Prob1} and a set of points \eword{(xb,yb,zb)} on the boundary
of $\Omega$, let gb = g(xb,yb,zb). 
Suppose the differential operator $L$ is as in (\ref{L1})
with coefficients \eword{a1,a2,a3,a4,a5,a6}. Here is  a
function for computing the coefficients of the IPBC \tp{} spline.
\ms

\sb{4.1}\vbox{\ob
 [C,H,t1,t2] = ipbctp3(dx,dy,dz,xe,ye,ze,f,xb,yb,zb,...
~~~   a1,a2,a3,a4,a5,a6,gb,lambda,mc)
}\eb
\ms\noin 
The inputs are the same as for \eword{ipbftp3} except for the parameter
\eword{mc} which controls the selection of the collocation points. 
This function does not require the use of  Gauss quadrature 
since no integrals are being computed.
The main output of the function is the array \eword{C}
containing the coefficients of the IPBC \tp{} spline.  For diagnostic purposes,
we also output the matrix
\eword{H} of the system being solved by least-squares, and the times to assemble and
solve it.

We have also written a simplified version of this function that can be used when the
PDE is the Laplace operator.

\ms
\sb{4.75}\vbox{\ob
[C,H,t1,t2] = ipbctp3L(dx,dy,dz,xe,ye,ze,f,xb,yb,zb,gb,lambda,mc)
}\eb
\ms\noin
Based on examples in the supplement, we will choose 
$\lambda = .01$ as a  default value when calling these functions.
To test these functions we have written the scripts \eword{TP3ipbc} and
\eword{TP3ipbcL}. To explore rates of convergence of these methods we
can employ the scripts  \eword{TP3ipbcr} and
\eword{TP3ipbcLr}.

\section{Performance of the methods using \tp{} splines}
In this section we present several examples to
illustrate the use of both the IPBF and IPBC methods.
For each example we choose a domain and a PDE with a known true solution $u$.
We then set $f = Lu$ and $g = u|_\pO$.  Since we know the true solution,
we can compute and tabulate errors by comparing it with our spline
solution  on a large set of points covering $\Omega$. 
The  examples  all include tables giving the computational times,
the numbers $m$ of grid lines in each direction forming the box partition,
the times to assemble the systems and solve them, the numbers nc of coefficients, 
and the
max and RMS errors measured on a large set of points covering $\Omega$.
We run each example for  a sequence of  \tp{} spline spaces with decreasing mesh 
sizes. This allows us to compute estimated rates of convergence, see Remark~\ref{Rates}.
They are displayed in the last two columns of the tables, where 
the first refers to max errors, and the second to RMS errors.
The examples all show that 
the methods produce highly accurate results
in relatively short computational times.  

To compute the coefficients of an approximating spline with either
IPBF or IPBC, we have to solve
an overdetermined linear system of equations with a sparse matrix $H$.
To gauge how well this system is being solved, we could compute 
its 2-condition number, see Remark~\ref{CN}. However, as explained in the
remark this can be problematic for large systems, and so 
in all of our examples we report an estimate for it, namely
the square root of the 1-norm condition number of $H^TH$.

\subsection{Comparison of IPBF with IPBC}
In this section we present two examples where we run both IPBF and IPBC 
on the same problem.

\begin{figure} 
\centering
\subfloat{\includegraphics[width = .5\linewidth, trim=0 80 0 80, clip]{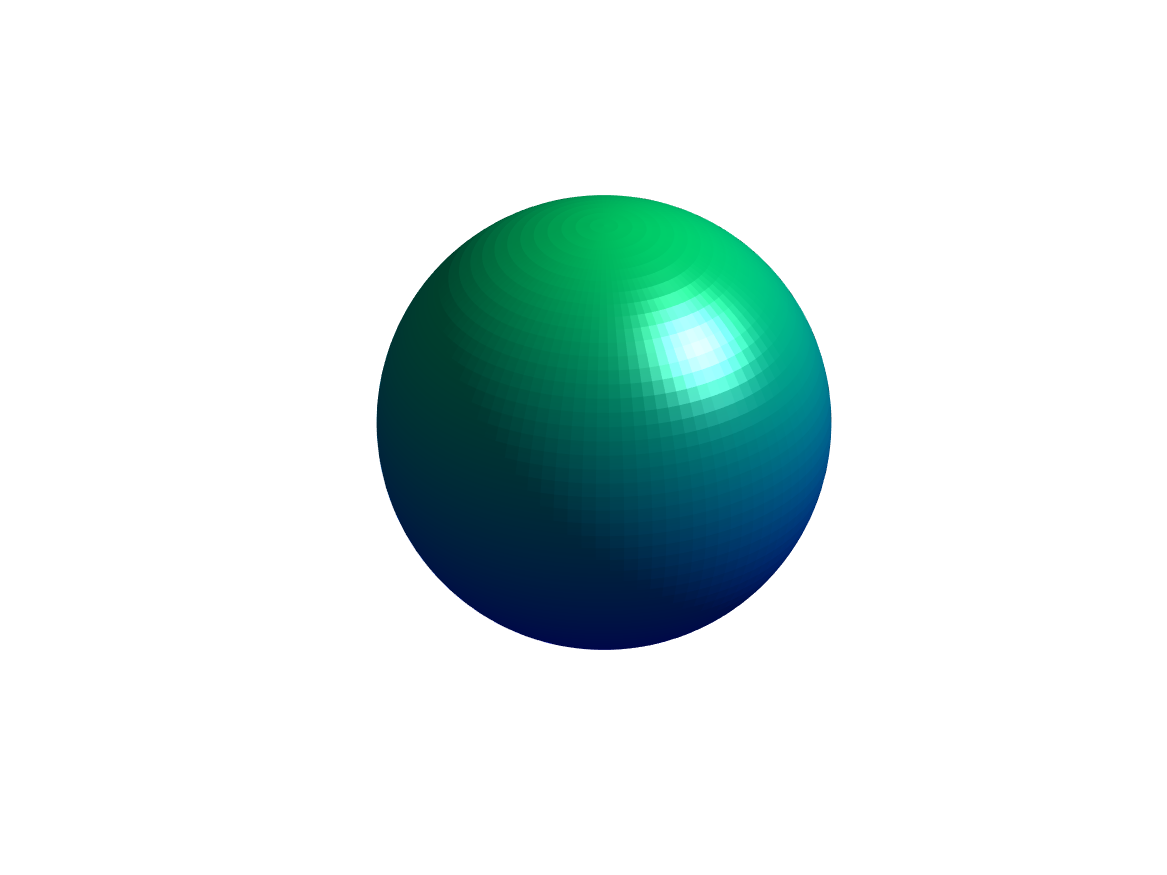}}
\raisebox{0.025in}{\subfloat{\includegraphics[width = .5\linewidth, trim=0 80 0 80, clip]{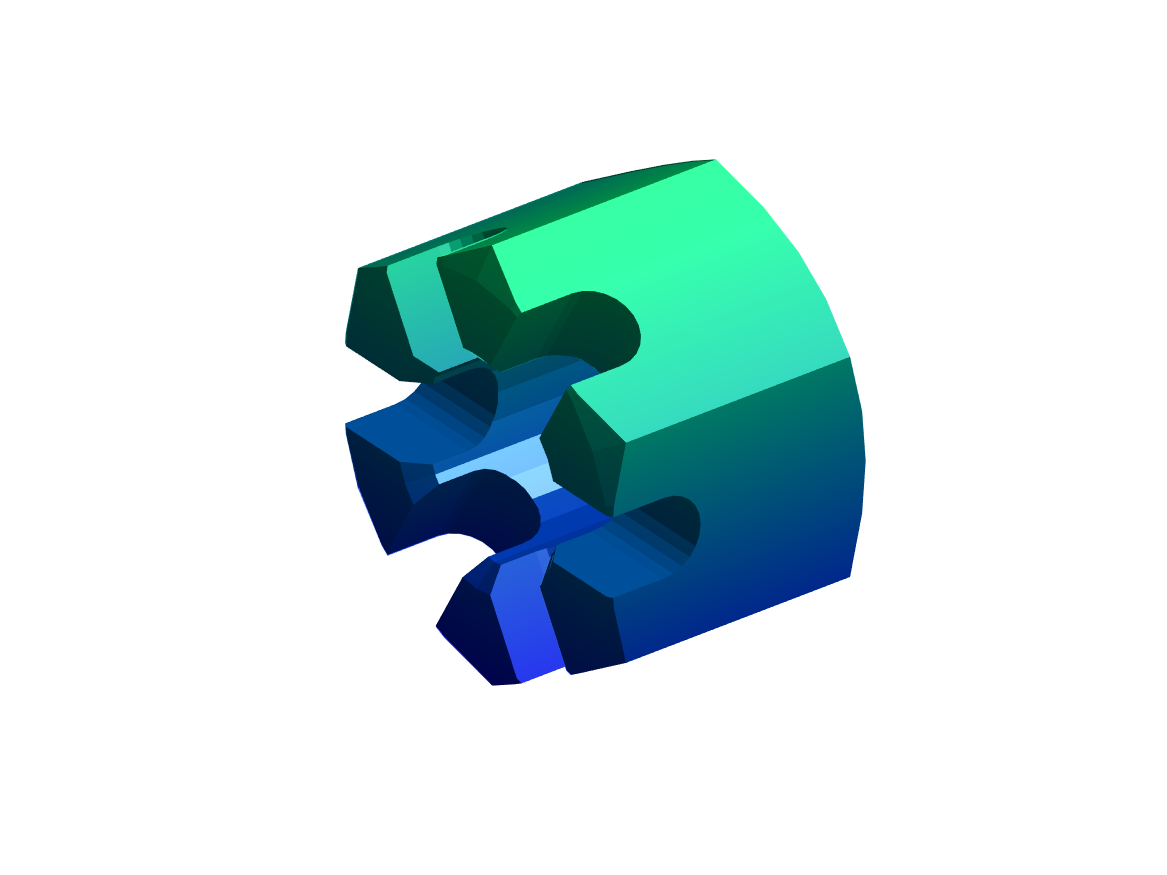}}}
\caption{The sphere and the nut domains}
\vskip -.1in 
 \label{Fig1} 
\end{figure}

\begin{ex}\label{Exsphere}    Solve Problem~\ref{Prob1} with the Laplace operator
on the  sphere shown in \frl{Fig1}, 
with $f$ and $g$ such that
the true solution is $\sin(5x)\sin(5y)\sin(5z)$.
Use both the IPBF and IPBC methods with tensor-product splines of degree (d,d,d) 
defined on $m\times m\times m$ grids covering the bounding box of $\Omega$ 
for a sequence of increasing values of m.
\end{ex}

\dis 
 We use the scripts \eword{TP3ipbfLr} and \eword{TP3ipbcLr}, respectively,
with  the file \eword{sphere50.dat}, and with
\eword{pde = 5} and  \eword{nb = 2000}. We run these scripts with both
\eword{d = 4} and \eword{d = 5}. 
The errors here are computed
on a set of 103,283 points covering $\Omega$.

\sm\noin A) IPBF: $d = 4$

\centerline{\vbox{
\ms
\offinterlineskip
\halign{\strut \hfil # &  \hfil #  & \hfil # \hfil  &  \hfil # \hfil  &  ~~# & ~~#  & ~~#\cr
m & setup & solve & nc& CN  & emax  & rms \cr
\noalign{\hrule}
\noalign{\vskip 1 pt}
5 & 0.19 & 0.13 & 512 & 1.76e+04 &  3.48e-03 &  5.27e-04 \cr 
6 & 0.22 & 0.19 & 729 & 4.12e+04 &  6.82e-04 &  1.55e-04 \cr 
7 & 0.34 & 0.28 & 1000 & 7.78e+04 &  2.98e-04 &  5.62e-05 \cr 
8 & 0.59 & 0.41 & 1331 & 1.74e+05 &  1.53e-04 &  2.67e-05 \cr 
9 & 0.95 & 0.61 & 1728 & 3.92e+05 &  6.24e-05 &  1.28e-05 \cr 
10 & 1.52 & 0.32 & 2197 & 7.63e+05 &  3.46e-05 &  6.50e-06 \cr 
}}
\hskip 2pc
\vbox{\offinterlineskip
\halign{\strut
\hfil  \hfil #  &  ~~#&   ~~#   \cr
\omit\span \hfil rates \hfil \cr
\noalign{\vskip 4 pt}
\noalign{\hrule}
\noalign{\vskip 1 pt}
7.30 & 5.48 \cr 
 4.54 & 5.57 \cr 
 4.33 & 4.82 \cr 
 6.71 & 5.53 \cr 
 5.01 & 5.72 \cr 
}}}

\sm\noin B) IPBF: $d = 5$

\centerline{\vbox{
\offinterlineskip
\halign{\strut \hfil # &  \hfil #  & \hfil # \hfil  &  \hfil # \hfil  &  ~~# & ~~#  & ~~#\cr
m & setup & solve & nc& CN  & emax  & rms \cr
\noalign{\hrule}
\noalign{\vskip 1 pt}
5 & 0.29 & 0.37 & 729 & 2.35e+05 &  9.61e-04 &  1.20e-04 \cr 
6 & 0.55 & 0.54 & 1000 & 2.06e+05 &  2.33e-04 &  2.46e-05 \cr 
7 & 0.94 & 0.83 & 1331 & 5.87e+05 &  5.87e-05 &  7.25e-06 \cr 
8 & 1.60 & 1.28 & 1728 & 1.47e+06 &  1.95e-05 &  2.81e-06 \cr 
9 & 2.56 & 0.38 & 2197 & 2.75e+06 &  9.54e-06 &  1.44e-06 \cr 
10 & 4.23 & 0.58 & 2744 & 3.67e+06 &  5.52e-06 &  6.94e-07 \cr 
}}
\hskip 2pc
\vbox{\offinterlineskip
\halign{\strut
\hfil  \hfil #  &  ~~#&   ~~#   \cr
\omit\span \hfil rates \hfil \cr
\noalign{\vskip 4 pt}
\noalign{\hrule}
\noalign{\vskip 1 pt}
6.36 & 7.11 \cr 
 7.55 & 6.69 \cr 
 7.15 & 6.16 \cr 
 5.35 & 4.98 \cr 
 4.65 & 6.21 \cr 
}}}

\sm\noin C) IPBC: $d = 4$

\centerline{\vbox{
\ms
\offinterlineskip
\halign{\strut \hfil # &  \hfil #  & \hfil # \hfil  &  \hfil # \hfil  &  ~~# & ~~#  & ~~#\cr
m & setup & solve & nc& CN  & emax  & rms \cr
\noalign{\hrule}
\noalign{\vskip 1 pt}
5 & 0.10 & 0.10 & 512 & 3.94e+03 &  5.91e-03 &  1.10e-03 \cr 
6 & 0.15 & 0.15 & 729 & 4.54e+03 &  1.85e-03 &  3.18e-04 \cr 
7 & 0.27 & 0.21 & 1000 & 4.71e+03 &  6.56e-04 &  1.37e-04 \cr 
8 & 0.50 & 0.32 & 1331 & 6.82e+03 &  4.26e-04 &  7.03e-05 \cr 
9 & 0.95 & 0.43 & 1728 & 9.19e+03 &  2.23e-04 &  4.01e-05 \cr 
10 & 1.62 & 0.69 & 2197 & 1.61e+04 &  1.31e-04 &  2.46e-05 \cr 
}}
\hskip 2pc
\vbox{\offinterlineskip
\halign{\strut
\hfil  \hfil #  &  ~~#&   ~~#   \cr
\omit\span \hfil rates \hfil \cr
\noalign{\vskip 4 pt}
\noalign{\hrule}
\noalign{\vskip 1 pt}
5.20 & 5.55 \cr 
 5.68 & 4.63 \cr 
 2.81 & 4.32 \cr 
 4.83 & 4.21 \cr 
 4.55 & 4.14 \cr 
}}}

\sm\noin D) IPBC: $d = 5$

\centerline{\vbox{
\offinterlineskip
\halign{\strut \hfil # &  \hfil #  & \hfil # \hfil  &  \hfil # \hfil  &  ~~# & ~~#  & ~~#\cr
m & setup & solve & nc& CN  & emax  & rms \cr
\noalign{\hrule}
\noalign{\vskip 1 pt}
5 & 0.17 & 0.24 & 729 & 2.58e+04 &  2.42e-03 &  3.93e-04 \cr 
6 & 0.33 & 0.33 & 1000 & 2.95e+04 &  6.65e-04 &  1.09e-04 \cr 
7 & 0.64 & 0.53 & 1331 & 4.05e+04 &  2.36e-04 &  4.78e-05 \cr 
8 & 1.28 & 0.82 & 1728 & 3.36e+04 &  1.43e-04 &  2.42e-05 \cr 
9 & 2.37 & 1.36 & 2197 & 3.82e+04 &  7.62e-05 &  1.37e-05 \cr 
10 & 4.37 & 2.11 & 2744 & 4.28e+04 &  4.45e-05 &  8.32e-06 \cr 
}}
\hskip 2pc
\vbox{\offinterlineskip
\halign{\strut
\hfil  \hfil #  &  ~~#&   ~~#   \cr
\omit\span \hfil rates \hfil \cr
\noalign{\vskip 4 pt}
\noalign{\hrule}
\noalign{\vskip 1 pt}
5.78 & 5.74 \cr 
 5.69 & 4.54 \cr 
 3.23 & 4.40 \cr 
 4.72 & 4.30 \cr 
 4.56 & 4.21 \cr 
}}}
\ms

We now solve the same BVP, but on  a more complicated domain.

\begin{ex}\label{Exnut}    Solve Problem~\ref{Prob1} with the Laplace operator 
on the domain shown in \frr{Fig1}
with $f$ and $g$ such that
the true solution is $\sin(5x)\sin(5y)\sin(5z)$.
Use both IPBF and IPBC with tensor-product splines of degree (d,d,d)
defined on $m\times m\times m$ grids covering the bounding box of $\Omega$ 
for a sequence of increasing values of m.
\end{ex}

\dis 
 We use the scripts  \eword{TP3ipbfLr} and \eword{TP3ipbcLr}, respectively, 
with  the file \eword{nut.stl},
and set \eword{pde = 5},  \eword{nb = 4000}. We run the scripts for both 
 \eword{d = 4} and \eword{d = 5}. 
Here the errors are computed on  a set of 87,025 points covering $\Omega$.

\sm\noin A) IPBF: $d = 4$

\centerline{\vbox{
\offinterlineskip
\halign{\strut \hfil # &  \hfil #  & \hfil # \hfil  &  \hfil # \hfil  &  ~~# & ~~#  & ~~#\cr
m & setup & solve & nc& CN  & emax  & rms \cr
\noalign{\hrule}
\noalign{\vskip 1 pt}
5 & 0.17 & 0.18 & 512 & 2.74e+04 &  1.27e-03 &  3.03e-04 \cr 
6 & 0.28 & 0.26 & 729 & 5.73e+04 &  3.46e-04 &  8.89e-05 \cr 
7 & 0.43 & 0.36 & 1000 & 7.92e+04 &  1.38e-04 &  3.48e-05 \cr 
8 & 0.69 & 0.50 & 1331 & 1.23e+05 &  6.72e-05 &  1.41e-05 \cr 
9 & 1.08 & 0.71 & 1728 & 2.05e+05 &  3.35e-05 &  6.70e-06 \cr 
10 & 1.70 & 1.04 & 2197 & 2.58e+05 &  1.52e-05 &  3.47e-06 \cr 
}}
\hskip 2pc
\vbox{\offinterlineskip
\halign{\strut
\hfil  \hfil #  &  ~~#&   ~~#   \cr
\omit\span \hfil rates \hfil \cr
\noalign{\vskip 4 pt}
\noalign{\hrule}
5.83 & 5.50 \cr 
 5.06 & 5.15 \cr 
 4.64 & 5.87 \cr 
 5.21 & 5.55 \cr 
 6.72 & 5.59 \cr 
}}}

\sm\noin B) IPBF: $d = 5$

\centerline{\vbox{
\offinterlineskip
\halign{\strut \hfil # &  \hfil #  & \hfil # \hfil  &  \hfil # \hfil  &  ~~# & ~~#  & ~~#\cr
m & setup & solve & nc& CN  & emax  & rms \cr
\noalign{\hrule}
\noalign{\vskip 1 pt}
5 & 0.39 & 0.55 & 729 & 1.88e+05 &  3.41e-04 &  7.01e-05 \cr 
6 & 0.70 & 0.72 & 1000 & 2.69e+05 &  8.34e-05 &  1.46e-05 \cr 
7 & 1.11 & 0.98 & 1331 & 7.12e+05 &  2.48e-05 &  5.00e-06 \cr 
8 & 1.95 & 1.45 & 1728 & 1.43e+06 &  1.21e-05 &  1.89e-06 \cr 
9 & 3.06 & 2.03 & 2197 & 2.68e+06 &  4.79e-06 &  7.57e-07 \cr 
10 & 4.78 & 2.78 & 2744 & 4.58e+06 &  2.09e-06 &  3.63e-07 \cr 
}}
\hskip 2pc
\vbox{\offinterlineskip
\halign{\strut
\hfil  \hfil #  &  ~~#&   ~~#   \cr
\omit\span \hfil rates \hfil \cr
\noalign{\vskip 4 pt}
\noalign{\hrule}
\noalign{\vskip 1 pt}
6.31 & 7.04 \cr 
 6.65 & 5.86 \cr 
 4.66 & 6.30 \cr 
 6.93 & 6.86 \cr 
 7.03 & 6.23 \cr 
}}}

\sm\noin C) IPBC: $d = 4$

\centerline{\vbox{
\offinterlineskip
\halign{\strut \hfil # & \hfil # &  \hfil #  & ~~#  &  ~~# &  ~~# & ~~#  & ~~#\cr
m & setup & solve & nc& CN  & emax  & rms \cr
\noalign{\hrule}
\noalign{\vskip 1 pt}
5 & 0.12 & 0.16 & 512 & 2.05e+03 &  1.71e-03 &  3.11e-04 \cr 
6 & 0.22 & 0.22 & 729 & 2.25e+03 &  4.19e-04 &  8.76e-05 \cr 
7 & 0.36 & 0.30 & 1000 & 2.28e+03 &  2.42e-04 &  3.49e-05 \cr 
8 & 0.70 & 0.39 & 1331 & 3.71e+03 &  1.07e-04 &  1.64e-05 \cr 
9 & 1.12 & 0.52 & 1728 & 6.41e+03 &  6.85e-05 &  1.02e-05 \cr 
10 & 1.83 & 0.75 & 2197 & 1.27e+04 &  4.16e-05 &  6.05e-06 \cr 
}}
\hskip 2pc
\vbox{\offinterlineskip
\halign{\strut
\hfil  \hfil #  &  ~~#&   ~~#   \cr
\omit\span \hfil rates \hfil \cr
\noalign{\vskip 4 pt}
\noalign{\hrule}
\noalign{\vskip 1 pt}
6.29 & 5.68 \cr 
 3.02 & 5.05 \cr 
 5.27 & 4.92 \cr 
 3.36 & 3.56 \cr 
 4.24 & 4.40 \cr 
}}}

\sm\noin D) IPBC: $d = 5$

\centerline{\vbox{
\offinterlineskip
\halign{\strut \hfil # & \hfil # &  \hfil #  & ~~#  &  ~~# &  ~~# & ~~#  & ~~#\cr
m & setup & solve & nc& CN  & emax  & rms \cr
\noalign{\hrule}
\noalign{\vskip 1 pt}
5 & 0.25 & 0.38 & 729 & 1.87e+04 &  4.60e-04 &  7.72e-05 \cr 
6 & 0.48 & 0.49 & 1000 & 2.46e+04 &  1.00e-04 &  1.71e-05 \cr 
7 & 0.82 & 0.69 & 1331 & 4.08e+04 &  4.83e-05 &  6.48e-06 \cr 
8 & 1.52 & 1.05 & 1728 & 6.84e+04 &  2.46e-05 &  3.29e-06 \cr 
9 & 2.74 & 1.58 & 2197 & 1.55e+05 &  1.22e-05 &  1.72e-06 \cr 
10 & 5.11 & 2.43 & 2744 & 3.04e+05 &  7.83e-06 &  1.02e-06 \cr 
}}
\hskip 2pc
\vbox{\offinterlineskip
\halign{\strut
\hfil  \hfil #  &  ~~#&   ~~#   \cr
\omit\span \hfil rates \hfil \cr
\noalign{\vskip 4 pt}
\noalign{\hrule}
\noalign{\vskip 1 pt}
6.84 & 6.75 \cr 
 3.99 & 5.33 \cr 
 4.38 & 4.40 \cr 
 5.25 & 4.84 \cr 
 3.78 & 4.44 \cr 
}}} 
\ms\noin
Before proceeding,  we make some comments on these two examples. First we note
that  as we increase {\tt m} the mesh size decreases and
the accuracy improves,  while the
number of coefficients and computational times increase somewhat.
The condition numbers also increase somewhat for both methods, and they are generally
smaller for the IPBC method. In both cases they are larger for $d = 5$ than
for $d = 4$. For the IPBF method the rates of convergence 
are close to optimal, i.e., order $d+1$ when using
splines of degree $d$, see Remarks~\ref{Rates} -- \ref{Appxtet}.
The tables suggest that
IPBC has suboptimal convergence, at least one order less than optimal.
For a given accuracy, IPBC is often a little faster than IPBF.

\subsection{Some additional examples on other domains with other PDEs}
In this section we give several additional examples involving
different domains, different operators $L$, and more
complicated or less smooth true solutions. In all of these examples we use
\tp{} splines defined on the bounding box of the domain, and due to space restrictions,
we run each example for only one of the two methods.
Our first two examples involve domains with holes, and work with
PDEs with variable coefficients.

\begin{figure} 
\centering
\subfloat{\includegraphics[width = .5\linewidth,trim=0 45 0 45, clip]{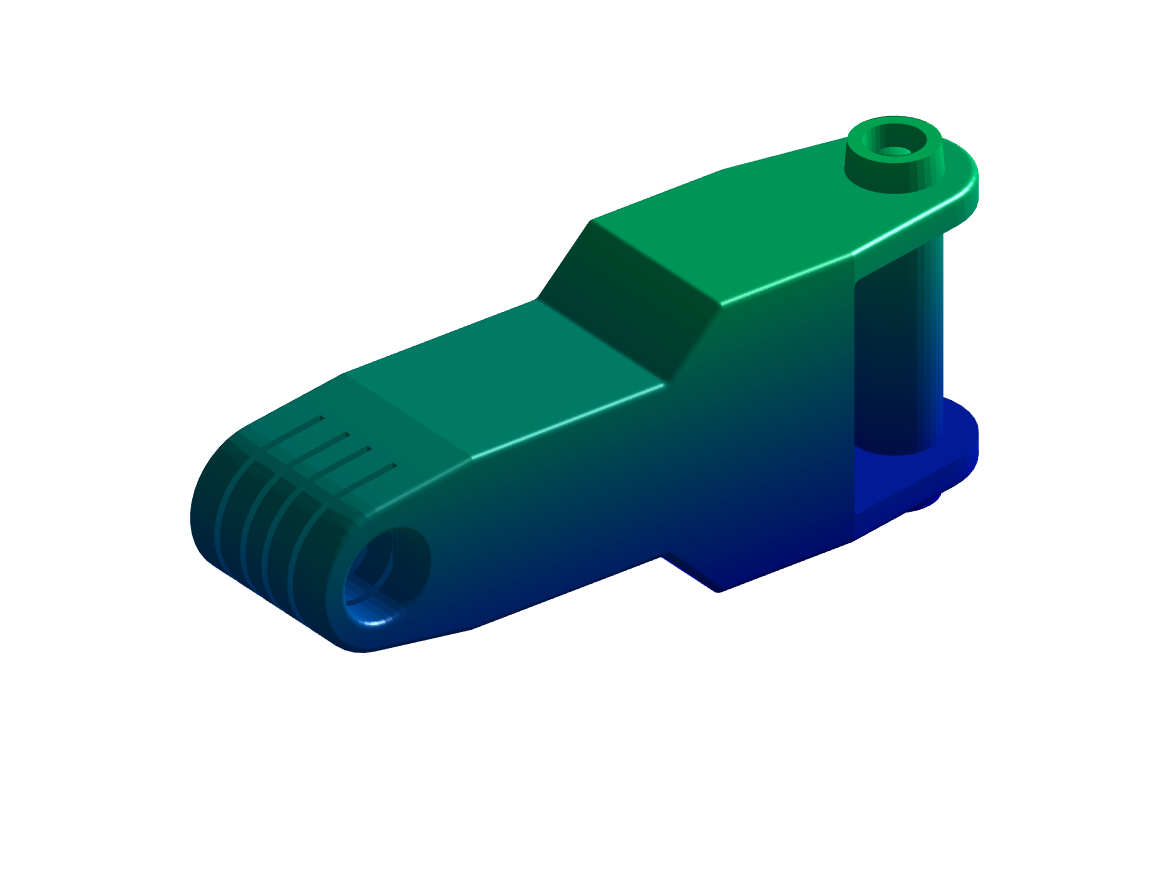}}
\raisebox{0.2in}{\subfloat{\includegraphics[width = .5\linewidth,trim=0 45 0 45, clip]{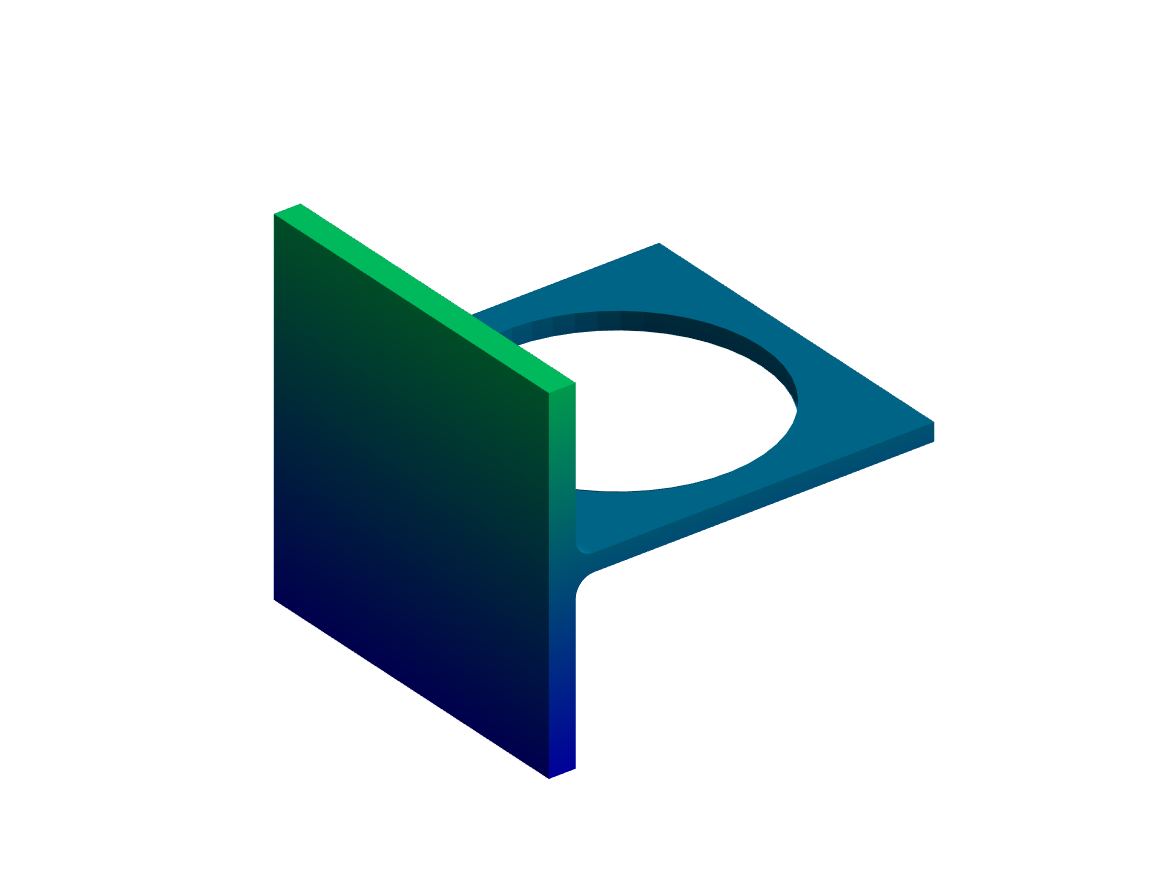}}}
\vskip -.1in
\caption{The forearmlink and bracket}  
\label{Fig2} 
\end{figure}

\begin{ex}\label{Exfore}    Solve Problem~\ref{Prob1} 
on the domain shown in \frl{Fig2} 
where  we take $a = (10+\cos(x) ,\exp(y), 10+\sin(x),0,0,0)$
for the coefficient vector defining $L$ and
$u := \sin(8x)+\sin(12y)+\sin(14z)$ for the true solution.
Use the IPBF method with \tp{} splines  of degree (d,d,d) on
 $m\times m \times m$   grids covering the  bounding box of $\Omega$.
\end{ex}

\dis    We use the script \eword{TP3ipbfr} with  the file \eword{ForearmLink.stl}, 
\eword{pde = 12},  \eword{nb = 5000},  
and \eword{d = 5} and \eword{d = 6}. 
The errors in this example are measured on a set of   92,477 points covering $\Omega$.
\sm

\noin A) IPBF: $d = 5$

\centerline{\vbox{
\offinterlineskip
\halign{\strut \hfil # & \hfil # &  \hfil #  & ~~#  &  ~~# &  ~~# & ~~#  & ~~#\cr
m & setup & solve & nc& CN  & emax  & rms \cr
\noalign{\hrule}
\noalign{\vskip 1 pt}
5 & 0.44 & 0.58 & 729 & 6.11e+05 &  2.79e-02 &  4.76e-03 \cr 
6 & 0.74 & 0.75 & 1000 & 7.75e+05 &  4.58e-03 &  1.03e-03 \cr 
7 & 1.24 & 1.08 & 1331 & 1.23e+06 &  1.37e-03 &  2.60e-04 \cr 
8 & 2.04 & 1.58 & 1728 & 2.20e+06 &  5.10e-04 &  9.88e-05 \cr 
9 & 3.34 & 2.10 & 2197 & 3.27e+06 &  1.98e-04 &  3.57e-05 \cr 
10 & 5.24 & 2.89 & 2744 & 5.51e+06 &  1.01e-04 &  1.63e-05 \cr 
}}
\hskip 2pc
\vbox{\offinterlineskip
\halign{\strut
\hfil  \hfil #  &  ~~#&   ~~#   \cr
\omit\span \hfil rates \hfil \cr
\noalign{\vskip 4 pt}
\noalign{\hrule}
\noalign{\vskip 1 pt}
8.10 & 6.88 \cr 
 6.62 & 7.54 \cr 
 6.40 & 6.27 \cr 
 7.10 & 7.62 \cr 
 5.69 & 6.68 \cr 
}}}

\sm\noin B) IPBF: $d = 6$

\centerline{\vbox{
\offinterlineskip
\halign{\strut \hfil # & \hfil # &  \hfil #  & ~~#  &  ~~# &  ~~# & ~~#  & ~~#\cr
m & setup & solve & nc& CN  & emax  & rms \cr
\noalign{\hrule}
\noalign{\vskip 1 pt}
5 & 0.88 & 1.40 & 1000 & 5.94e+06 &  4.22e-03 &  1.11e-03 \cr 
6 & 1.64 & 1.91 & 1331 & 7.31e+06 &  1.16e-03 &  3.12e-04 \cr 
7 & 2.97 & 2.80 & 1728 & 1.46e+07 &  1.93e-04 &  6.31e-05 \cr 
8 & 4.99 & 3.87 & 2197 & 2.56e+07 &  6.89e-05 &  1.86e-05 \cr 
9 & 7.99 & 5.49 & 2744 & 3.52e+07 &  3.19e-05 &  6.03e-06 \cr 
10 & 12.44 & 7.23 & 3375 & 4.78e+07 &  1.56e-05 &  2.46e-06 \cr 
}}
\hskip 2pc
\vbox{\offinterlineskip
\halign{\strut
\hfil  \hfil #  &  ~~#&   ~~#   \cr
\omit\span \hfil rates \hfil \cr
\noalign{\vskip 4 pt}
\noalign{\hrule}
\noalign{\vskip 1 pt}
5.81 & 5.68 \cr 
 9.81 & 8.76 \cr 
 6.69 & 7.92 \cr 
 5.76 & 8.44 \cr 
 6.08 & 7.60 \cr 
}}}
\ms

\begin{ex}\label{Exbracket}    Solve Problem~\ref{Prob1} on the domain
shown in \frr{Fig2}
with true solution
$\sin(5x) + \sin(5y) + \sin(5z)$.  Suppose the coefficients of $L$
are  given by\break
$a = (10+\cos(x) ,e^y, 10+\sin(x),\cos(x),\cos(y),\cos(z))$.
Use the IPBC method with \tp{} splines 
of degree (d,d,d) on $m\times m\times m$ grids covering
the bounding box, and choose $m_c = 4$.
\end{ex}

 \dis We use the script \eword{TP3ipbcr} with  the file \eword{BracketWithHole.stl}, 
\eword{pde = 11},  \eword{nb = 2000},  \eword{mc = 4}, and \eword{d = 4} and \eword{d = 5}. 
With these coefficients the operator $L$ is elliptic for all
$(x,y,z) \in \Omega$, and we measure the errors on a set of
50,000 points on the boundary and  2,869 points inside the domain.

\sm\noin A) IPBC: $d = 4$

\centerline{\vbox{
\offinterlineskip
\halign{\strut \hfil # & \hfil # &  \hfil #  & ~~#  &  ~~# &  ~~# & ~~#  & ~~#\cr
m & setup & solve & nc& CN  & emax  & rms \cr
\noalign{\hrule}
\noalign{\vskip 1 pt} 
5 & 0.14 & 0.14 & 512 & 1.95e+04 &  1.52e-01 &  1.73e-02 \cr 
6 & 0.29 & 0.25 & 729 & 2.26e+04 &  1.75e-02 &  2.15e-03 \cr 
7 & 0.64 & 0.41 & 1000 & 1.99e+04 &  5.19e-03 &  5.93e-04 \cr 
8 & 1.29 & 0.60 & 1331 & 2.62e+04 &  1.79e-03 &  2.15e-04 \cr 
9 & 2.66 & 0.92 & 1728 & 3.00e+04 &  6.09e-04 &  6.39e-05 \cr 
10 & 6.19 & 1.50 & 2197 & 3.01e+04 &  1.99e-04 &  1.86e-05 \cr 
}}
\hskip 2pc
\vbox{\offinterlineskip
\halign{\strut
\hfil  \hfil #  &  ~~#&   ~~#   \cr
\omit\span \hfil rates \hfil \cr
\noalign{\vskip 4 pt}
\noalign{\hrule}
\noalign{\vskip 1 pt}
9.68 & 9.34 \cr 
 6.67 & 7.06 \cr 
 6.90 & 6.60 \cr 
 8.08 & 9.07 \cr 
 9.50 & 10.46 \cr 
}}}

\sm\noin B) IPBC: $d = 5$

\centerline{\vbox{
\offinterlineskip
\halign{\strut \hfil # & \hfil # &  \hfil #  & ~~#  &  ~~# &  ~~# & ~~#  & ~~#\cr
m & setup & solve & nc& CN  & emax  & rms \cr
\noalign{\hrule}
\noalign{\vskip 1 pt}
5 & 0.23 & 0.31 & 729 & 8.11e+04 &  1.78e-02 &  2.70e-03 \cr 
6 & 0.62 & 0.57 & 1000 & 7.70e+04 &  5.98e-03 &  7.74e-04 \cr 
7 & 1.54 & 1.05 & 1331 & 1.00e+05 &  1.81e-03 &  1.68e-04 \cr 
8 & 3.79 & 1.94 & 1728 & 1.21e+05 &  3.12e-04 &  2.40e-05 \cr 
9 & 10.26 & 3.16 & 2197 & 1.57e+05 &  6.06e-05 &  4.95e-06 \cr 
10 & 24.37 & 5.01 & 2744 & 2.11e+05 &  2.47e-05 &  1.87e-06 \cr 
}}
\hskip 2pc
\vbox{\offinterlineskip
\halign{\strut
\hfil  \hfil #  &  ~~#&   ~~#   \cr
\omit\span \hfil rates \hfil \cr
\noalign{\vskip 4 pt}
\noalign{\hrule}
\noalign{\vskip 1 pt}
4.88 & 5.60 \cr 
 6.56 & 8.39 \cr 
 11.41 & 12.62 \cr 
 12.27 & 11.83 \cr 
 7.61 & 8.25 \cr 
}}}
\ms

\begin{figure} 
\centering
\subfloat{\includegraphics[width = .55\linewidth,trim=0 100 0 100, clip]{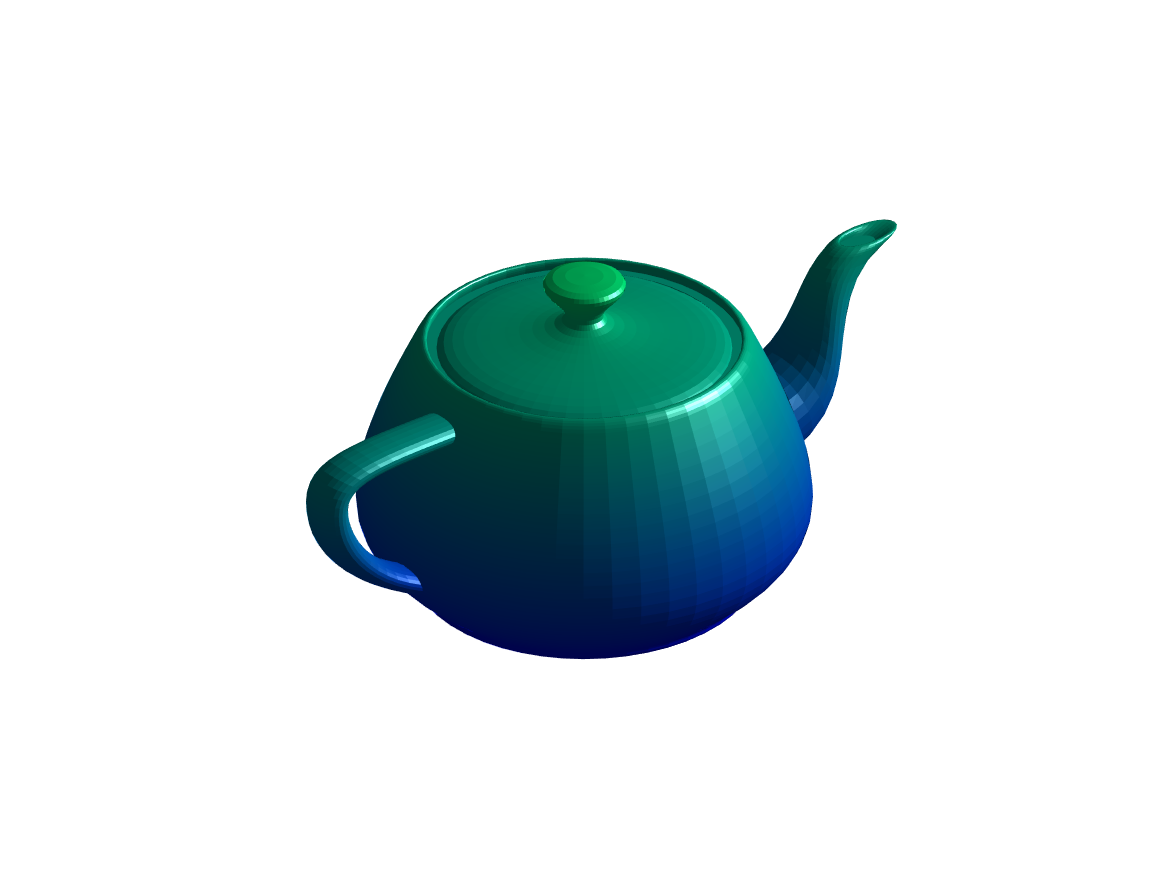}}
\raisebox{0.1in}{\subfloat{\includegraphics[width = .45\linewidth,trim=0 100 0 100, clip]{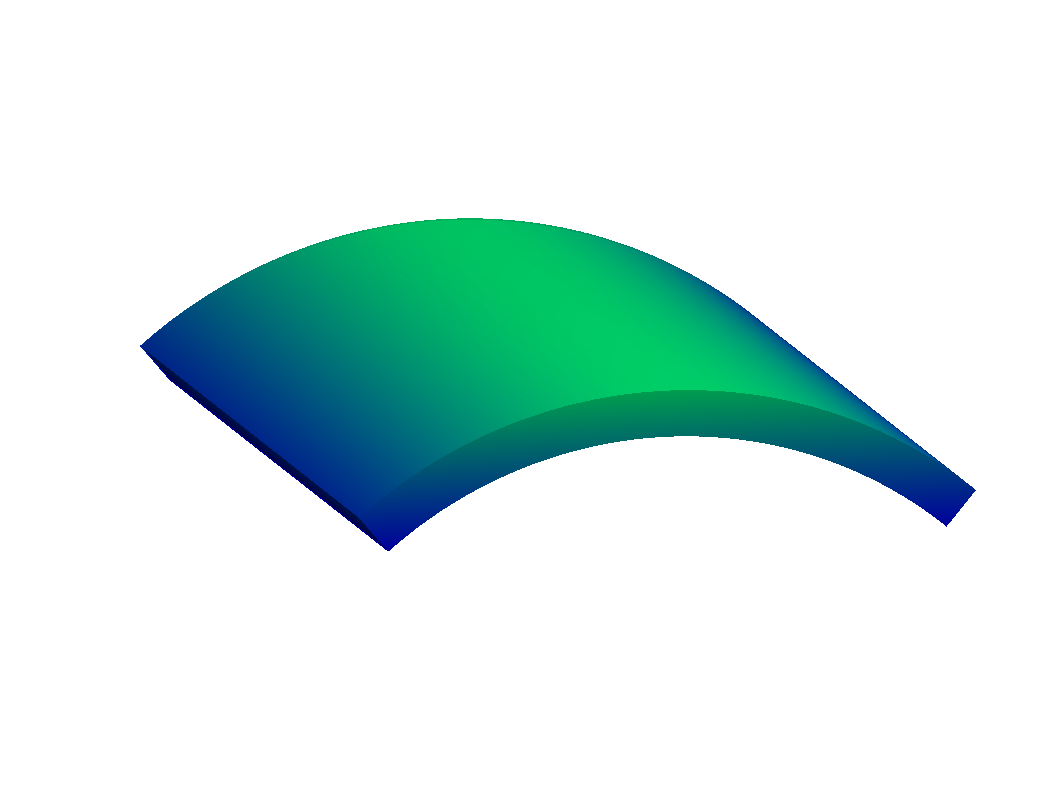}}}
\caption{The teapot and  curved roof domains}
\vskip -.1in
\label{Fig3} 
\end{figure}

\begin{ex}\label{Exteapot}    Solve Problem~\ref{Prob1} 
on the domain shown in \frl{Fig3} 
where the coefficient vector defining $L$ is
$a = (10+\cos(x) ,\exp(y), 10+\sin(x),0,0,0)$.
Suppose the true solution is $u := \sin(8x) + \sin(12y) + \sin(14z)$.
Use the IPBC method using \tp{} splines 
of degree (d,d,d) on
 $m\times m \times m$   grids covering the  bounding box of $\Omega$.
\end{ex}

\dis    We use the script \eword{TP3ipbcr} with  the file \eword{teapot.stl}, 
\eword{pde = 12},  \eword{nb = 2000},  \eword{mc = 3},
and run both \eword{d = 4} and \eword{d = 5}.
With this coefficient vector, the operator $L$ is elliptic on $\Omega$.
The errors in this example are measured on 50,000 boundary points and 83,699 points in the
interior of $\Omega$.
\sm

\noin A) IPBC: $d = 4$

\centerline{\vbox{
\offinterlineskip
\halign{\strut \hfil # & \hfil # &  \hfil #  & ~~#  &  ~~# &  ~~# & ~~#  & ~~#\cr
m & setup & solve & nc& CN  & emax  & rms \cr
\noalign{\hrule}
\noalign{\vskip 1 pt}
5 & 0.08 & 0.10 & 512 & 9.19e+04 &  1.64e-01 &  2.58e-02 \cr 
6 & 0.15 & 0.14 & 729 & 1.10e+05 &  1.77e-02 &  3.85e-03 \cr 
7 & 0.29 & 0.23 & 1000 & 1.03e+05 &  6.78e-03 &  1.39e-03 \cr 
8 & 0.53 & 0.31 & 1331 & 1.40e+05 &  2.88e-03 &  5.93e-04 \cr 
9 & 0.99 & 0.46 & 1728 & 1.84e+05 &  1.39e-03 &  2.81e-04 \cr 
10 & 1.69 & 0.69 & 2197 & 2.74e+05 &  9.71e-04 &  2.05e-04 \cr 
}}
\hskip 2pc
\vbox{\offinterlineskip
\halign{\strut
\hfil  \hfil #  &  ~~#&   ~~#   \cr
\omit\span \hfil rates \hfil \cr
\noalign{\vskip 4 pt}
\noalign{\hrule}
\noalign{\vskip 1 pt}
9.97 & 8.52 \cr 
 5.26 & 5.58 \cr 
 5.56 & 5.54 \cr 
 5.47 & 5.59 \cr 
 3.03 & 2.70 \cr 
}}}

\noin B) IPBC: $d = 5$

\centerline{\vbox{
\offinterlineskip
\halign{\strut \hfil # & \hfil # &  \hfil #  & ~~#  &  ~~# &  ~~# & ~~#  & ~~#\cr
m & setup & solve & nc& CN  & emax  & rms \cr
\noalign{\hrule}
\noalign{\vskip 1 pt}
5 & 0.16 & 0.23 & 729 & 5.81e+05 &  5.87e-02 &  9.84e-03 \cr 
6 & 0.33 & 0.32 & 1000 & 6.35e+05 &  9.73e-03 &  1.67e-03 \cr 
7 & 0.68 & 0.51 & 1331 & 5.36e+05 &  2.27e-03 &  4.59e-04 \cr 
8 & 1.27 & 0.80 & 1728 & 7.02e+05 &  9.46e-04 &  2.01e-04 \cr 
9 & 2.27 & 1.34 & 2197 & 8.69e+05 &  4.63e-04 &  9.22e-05 \cr 
10 & 4.36 & 2.06 & 2744 & 1.19e+06 &  2.75e-04 &  5.62e-05 \cr 
}}
\hskip 2pc
\vbox{\offinterlineskip
\halign{\strut
\hfil  \hfil #  &  ~~#&   ~~#   \cr
\omit\span \hfil rates \hfil \cr
\noalign{\vskip 4 pt}
\noalign{\hrule}
\noalign{\vskip 1 pt}
8.06 & 7.94 \cr 
 7.98 & 7.10 \cr 
 5.68 & 5.36 \cr 
 5.35 & 5.82 \cr 
 4.41 & 4.21 \cr 
}}}
\ms

The previous three examples worked with domains defined by stl files.
The next three work with point clouds.
The  domain in the following example was designed with
NURBS and then converted to a point cloud. 

\begin{ex}\label{Exroof}    
 Solve Problem~\ref{Prob1} on the domain
shown in \frr{Fig3} for the Laplace operator with true solution
$ \sin(8x)\sin(12y)\sin(14z)$.
Use the IPBC method based on \tp{} splines
of degree (d,d,d) on $m\times m\times m$ grids covering
the bounding box, and choose $m_c = 4$.
\end{ex}

 \dis We use the script \eword{TP3ipbcLr} with  the file \eword{roof.dat},
\eword{pde = 10},  \eword{nb = 5000}, and \eword{mc = 4}.
We run the script for both \eword{d = 4} and \eword{d = 5}. 
The errors here are measured on a set of  124,782 points in $\Omega$.

\sm\noin A) IPBC: $d = 4$  

\centerline{\vbox{  
\offinterlineskip
\halign{\strut \hfil # & \hfil # &  \hfil #  & ~~#  &  ~~# &  ~~# & ~~#  & ~~#\cr
m & setup & solve & nc& CN  & emax  & rms \cr
\noalign{\hrule}
\noalign{\vskip 1 pt}
6 & 0.38 & 0.38 & 729 & 6.87e+03 &  1.63e-01 &  2.88e-02 \cr 
7 & 0.82 & 0.58 & 1000 & 1.29e+04 &  7.36e-02 &  8.66e-03 \cr 
8 & 1.52 & 0.78 & 1331 & 3.59e+04 &  2.49e-02 &  3.20e-03 \cr 
9 & 2.93 & 1.10 & 1728 & 6.84e+04 &  1.41e-02 &  1.50e-03 \cr 
10 & 6.78 & 1.62 & 2197 & 1.28e+05 &  6.70e-03 &  7.26e-04 \cr 
}}
\hskip 2pc
\vbox{\offinterlineskip
\halign{\strut
\hfil  \hfil #  &  ~~#&   ~~#   \cr
\omit\span \hfil rates \hfil \cr
\noalign{\vskip 4 pt}
\noalign{\hrule}
\noalign{\vskip 1 pt}
 4.36 & 6.59 \cr 
 7.04 & 6.47 \cr 
 4.27 & 5.67 \cr 
 6.30 & 6.16 \cr 
}}}

\sm\noin B) IPBC: $d = 5$ 

\centerline{\vbox{
\offinterlineskip
\halign{\strut \hfil # & \hfil # &  \hfil #  & ~~#  &  ~~# &  ~~# & ~~#  & ~~#\cr
m & setup & solve & nc& CN  & emax  & rms \cr
\noalign{\hrule}
\noalign{\vskip 1 pt}
6 & 0.92 & 0.89 & 1000 & 2.72e+04 &  1.16e-01 &  8.17e-03 \cr 
7 & 2.08 & 1.44 & 1331 & 4.54e+04 &  3.21e-02 &  4.21e-03 \cr 
8 & 4.25 & 2.38 & 1728 & 7.62e+04 &  1.95e-02 &  1.54e-03 \cr 
9 & 11.46 & 3.85 & 2197 & 1.80e+05 &  1.08e-02 &  7.10e-04 \cr 
10 & 26.89 & 5.86 & 2744 & 3.17e+05 &  7.47e-03 &  4.12e-04 \cr 
}}
\hskip 2pc
\vbox{\offinterlineskip
\halign{\strut
\hfil  \hfil #  &  ~~#&   ~~#   \cr
\omit\span \hfil rates \hfil \cr
\noalign{\vskip 4 pt}
\noalign{\hrule}
\noalign{\vskip 1 pt}
7.04 & 3.63 \cr 
 3.23 & 6.52 \cr 
 4.44 & 5.81 \cr 
 3.12 & 4.62 \cr 
}}}
\ms

\begin{figure}
\centering
\subfloat{\includegraphics[width = .5\linewidth,trim=0 80 0 80, clip]{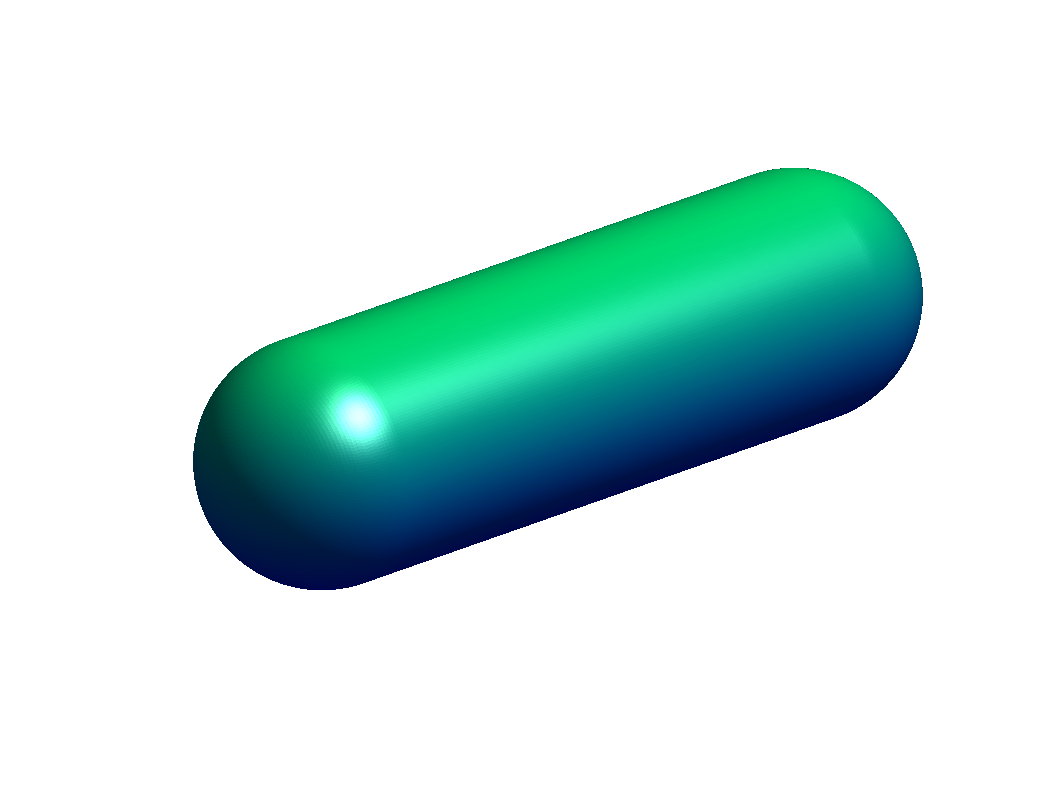}}
\hspace{0.02\linewidth}
\raisebox{0.2in}{
\subfloat{\includegraphics[width = .4\linewidth,trim=0 80 0 80, clip]{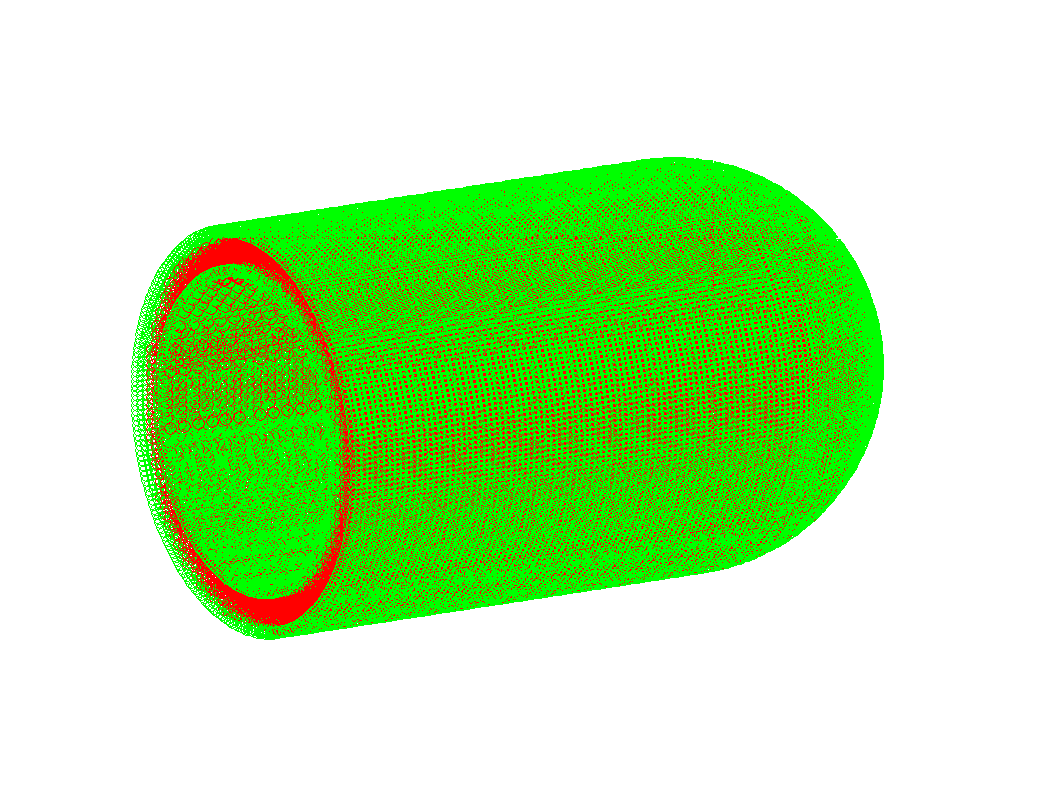}}}
\caption{A domain with a cavity and a slice through it} 
\vskip -.1in
\label{Fig4} 
\end{figure}

We now give an  example to show that our methods also work well for
domains with cavities. 
The boundary defining the inner surface
is just a scaled version of the outer surface. A slice of the domain
showing the cavity can be seen in \frr{Fig4}.

\begin{ex}\label{Exvessel}    
 Solve Problem~\ref{Prob1} on the domain
shown in \frl{Fig4} with the Laplace operator.  Choose $f$ and $g$ so that
the true solution is 
$\sin(5x)\sin(5y)\sin(5z)$.  
Use the IPBF method based on \tp{} splines
of degree (d,d,d) on $m\times m\times m$ grids covering
the bounding box. 
\end{ex}

\dis We use the script \eword{TP3ipbhfLr} with  the file \eword{tube.dat},
\eword{pde = 5} and   \eword{nb = 5000}.
We run the script for both \eword{d = 4} and \eword{d = 5}. 
The errors here are measured on a set of  185,900 points in $\Omega$.
\sm

\noin A) IPBF: $d = 4$

\centerline{\vbox{
\offinterlineskip
\halign{\strut \hfil # & \hfil # &  \hfil #  & ~~#  &  ~~# &  ~~# & ~~#  & ~~#\cr
m & setup & solve & nc& CN  & emax  & rms \cr
\noalign{\hrule}
\noalign{\vskip 1 pt} 
5 & 0.21 & 0.24 & 512 & 1.43e+04 &  9.95e-04 &  2.26e-04 \cr 
6 & 0.35 & 0.32 & 729 & 2.20e+04 &  2.62e-04 &  6.98e-05 \cr 
7 & 0.54 & 0.43 & 1000 & 4.34e+04 &  9.90e-05 &  2.52e-05 \cr 
8 & 0.80 & 0.58 & 1331 & 5.99e+04 &  4.40e-05 &  1.08e-05 \cr 
9 & 1.20 & 0.79 & 1728 & 9.34e+04 &  2.24e-05 &  5.21e-06 \cr 
10 & 1.84 & 1.12 & 2197 & 1.45e+05 &  1.14e-05 &  2.77e-06 \cr 
}}
\hskip 2pc
\vbox{\offinterlineskip
\halign{\strut
\hfil  \hfil #  &  ~~#&   ~~#   \cr
\omit\span \hfil rates \hfil \cr
\noalign{\vskip 4 pt}
\noalign{\hrule}
\noalign{\vskip 1 pt}
5.98 & 5.26 \cr 
 5.34 & 5.58 \cr 
 5.26 & 5.51 \cr 
 5.06 & 5.45 \cr 
 5.73 & 5.35 \cr 
}}}
\sm\noin B) IPBF: $d = 5$

\centerline{\vbox{
\offinterlineskip
\halign{\strut \hfil # & \hfil # &  \hfil #  & ~~#  &  ~~# &  ~~# & ~~#  & ~~#\cr
m & setup & solve & nc& CN  & emax  & rms \cr
\noalign{\hrule}
5 & 0.46 & 0.62 & 729 & 1.45e+05 &  2.37e-04 &  5.34e-05 \cr 
6 & 0.84 & 0.85 & 1000 & 1.22e+05 &  8.10e-05 &  1.36e-05 \cr 
7 & 1.38 & 1.20 & 1331 & 2.14e+05 &  1.63e-05 &  3.81e-06 \cr 
8 & 2.25 & 1.72 & 1728 & 5.44e+05 &  6.43e-06 &  1.39e-06 \cr 
9 & 3.37 & 2.20 & 2197 & 8.60e+05 &  2.79e-06 &  5.89e-07 \cr 
10 & 5.37 & 3.08 & 2744 & 1.43e+06 &  1.58e-06 &  2.78e-07 \cr 
}}
\hskip 2pc
\vbox{\offinterlineskip
\halign{\strut
\hfil  \hfil #  &  ~~#&   ~~#   \cr
\omit\span \hfil rates \hfil \cr
\noalign{\vskip 4 pt}
\noalign{\hrule}
\noalign{\vskip 1 pt}
4.81 & 6.12 \cr 
 8.80 & 6.98 \cr 
 6.04 & 6.56 \cr 
 6.24 & 6.41 \cr 
 4.81 & 6.38 \cr 
}}}
\ms

\begin{figure} 
\centering
\hspace{0.1\linewidth}
\subfloat{\includegraphics[width = .25\linewidth,trim=0 0 0 0, clip]{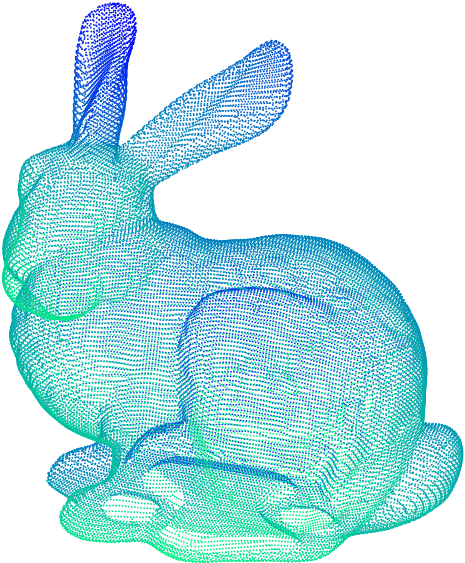}}
\subfloat{\includegraphics[width = .5\linewidth,trim=0 70 0 70, clip]{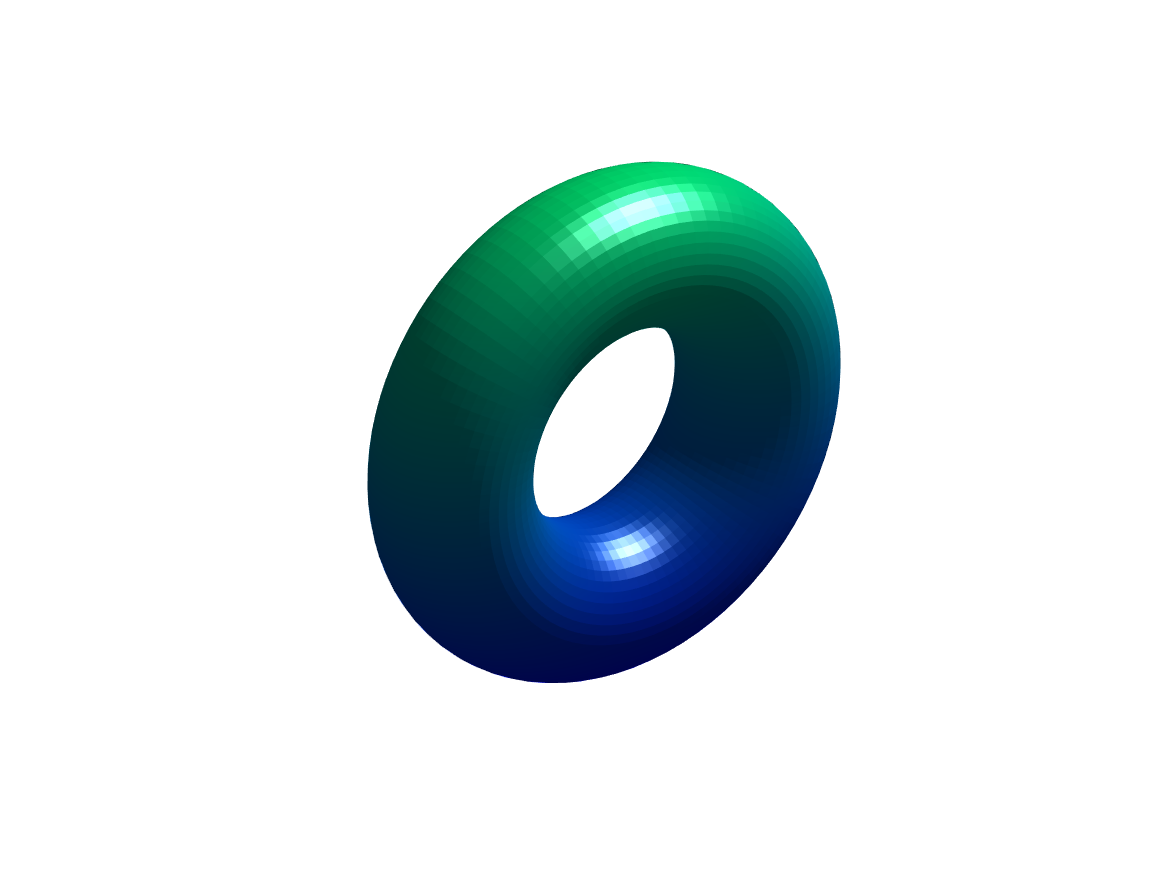}}
\caption{The Stanford bunny and a torus}
\vskip -.1in
\label{Fig5} 
\end{figure}

The results in this example are very similar to those obtained in Example~\ref{Exsphere}
where the same problem was solved on a domain without a cavity, namely
the sphere.
The following  example
works with a domain defined by a ply file
obtained by scanning a physical object,
namely the famous Stanford bunny.

\begin{ex}\label{Exbunny}    
 Solve Problem~\ref{Prob1} on the domain
shown in \frl{Fig5} with the Laplace operator.  Choose $f$ and $g$ so that
the true solution is 
$\sin(8x)\sin(12y)\sin(14z)$.  
Use the IPBC method based on \tp{} splines
of degree (d,d,d) on $m\times m\times m$ grids covering
the bounding box. 
\end{ex}

 \dis We use the script \eword{TP3ipbcLr}  with  the file \eword{bunny.ply},
\eword{pde = 10},  \eword{nb = 2000}, and \eword{mc = 3}. 
We run the script for both \eword{d = 4} and \eword{d = 5}. 
The errors here are measured on a set of  56,948 points inside $\Omega$.

\sm\noin A) IPBC: $d = 4$

\centerline{\vbox{
\offinterlineskip
\halign{\strut \hfil # & \hfil # &  \hfil #  & ~~#  &  ~~# &  ~~# & ~~#  & ~~#\cr
m & setup & solve & nc& CN  & emax  & rms \cr
\noalign{\hrule}
\noalign{\vskip 1 pt}
5 & 0.07 & 0.10 & 512 & 1.82e+04 &  4.46e-01 &  7.57e-02 \cr 
6 & 0.14 & 0.15 & 729 & 2.49e+04 &  1.31e-01 &  2.45e-02 \cr 
7 & 0.27 & 0.21 & 1000 & 3.17e+04 &  5.64e-02 &  7.93e-03 \cr 
8 & 0.51 & 0.31 & 1331 & 5.00e+04 &  2.03e-02 &  3.31e-03 \cr 
9 & 0.96 & 0.47 & 1728 & 1.16e+05 &  9.66e-03 &  1.43e-03 \cr 
10 & 1.66 & 0.68 & 2197 & 2.59e+05 &  6.21e-03 &  8.76e-04 \cr 
}}
\hskip 2pc
\vbox{\offinterlineskip
\halign{\strut
\hfil  \hfil #  &  ~~#&   ~~#   \cr
\omit\span \hfil rates \hfil \cr
\noalign{\vskip 4 pt}
\noalign{\hrule}
\noalign{\vskip 1 pt}
5.50 & 5.06 \cr 
 4.61 & 6.18 \cr 
 6.63 & 5.67 \cr 
 5.55 & 6.28 \cr 
 3.75 & 4.16 \cr 
}}}
\sm\noin B) IPBC: $d = 5$

\centerline{\vbox{
\offinterlineskip
\halign{\strut \hfil # & \hfil # &  \hfil #  & ~~#  &  ~~# &  ~~# & ~~#  & ~~#\cr
m & setup & solve & nc& CN  & emax  & rms \cr
\noalign{\hrule}
\noalign{\vskip 1 pt}
5 & 0.15 & 0.23 & 729 & 1.03e+06 &  5.90e-01 &  4.92e-02 \cr 
6 & 0.32 & 0.33 & 1000 & 6.32e+05 &  1.06e-01 &  9.95e-03 \cr 
7 & 0.65 & 0.50 & 1331 & 6.33e+05 &  3.25e-02 &  4.45e-03 \cr 
8 & 1.33 & 0.84 & 1728 & 5.94e+05 &  8.24e-03 &  1.43e-03 \cr 
9 & 2.39 & 1.36 & 2197 & 6.03e+05 &  3.61e-03 &  5.56e-04 \cr 
10 & 4.54 & 2.15 & 2744 & 8.67e+05 &  2.13e-03 &  3.17e-04 \cr 
}}
\hskip 2pc
\vbox{\offinterlineskip
\halign{\strut
\hfil  \hfil #  &  ~~#&   ~~#   \cr
\omit\span \hfil rates \hfil \cr
\noalign{\vskip 4 pt}
\noalign{\hrule}
\noalign{\vskip 1 pt}
7.69 & 7.16 \cr 
 6.50 & 4.41 \cr 
 8.89 & 7.38 \cr 
 6.17 & 7.07 \cr 
 4.47 & 4.77 \cr 
}}}
\ms

All of the above examples are for BVPs where 
the true solutions are infinitely differentiable.
We now explore what happens when the true solution has limited smoothness.
In this case we expect that our methods will provide a reduced order of convergence
since  it is well known that when approximating functions that are only 
$C^r$ smooth by splines of degree 
$d\ge r$, we will get
only order $r+1$ convergence  instead of the optimal order $d+1$.

\begin{ex}\label{Extorus}    Solve Problem~\ref{Prob1}
on the domain shown in \frr{Fig5} for the Laplace operator
where $f$ and $g$ are chosen so that the
true solution is $u := |x+y-z|^3$.
Use the IPBF method with tensor-product splines of degree (d,d,d)
defined on $m\times m \times m$ 
grids covering the  bounding box of $\Omega$.
\end{ex}

 \dis We use the script \eword{TP3ipbfLr} with  the file \eword{Torus.stl},
\eword{pde = 7} and  \eword{nb = 2000}.
We run the script for  both \eword{d = 2} and \eword{d = 5}.
For this example the
function $f$ is only $C^0$ and the true solution is only $C^2$.
The errors here are measured on 99,522 points.

\sm\noin A) IPBF: $d = 2$

\centerline{\vbox{
\offinterlineskip
\halign{\strut \hfil # & \hfil # &  \hfil #  & ~~#  &  ~~# &  ~~# & ~~#  & ~~#\cr
m & setup & solve & nc& CN  & emax  & rms \cr
\noalign{\hrule}
\noalign{\vskip 1 pt}
5 & 0.02 & 0.01 & 216 & 3.18e+03 &  4.58e-03 &  1.23e-03 \cr 
6 & 0.03 & 0.02 & 343 & 4.52e+03 &  2.47e-03 &  6.67e-04 \cr 
7 & 0.05 & 0.03 & 512 & 5.50e+03 &  1.55e-03 &  4.07e-04 \cr 
8 & 0.06 & 0.04 & 729 & 7.44e+03 &  1.07e-03 &  2.67e-04 \cr 
9 & 0.10 & 0.06 & 1000 & 9.73e+03 &  7.81e-04 &  1.87e-04 \cr 
10 & 0.16 & 0.09 & 1331 & 1.24e+04 &  5.38e-04 &  1.33e-04 \cr 
}}
\hskip 2pc
\vbox{\offinterlineskip
\halign{\strut
\hfil  \hfil #  &  ~~#&   ~~#   \cr
\omit\span \hfil rates \hfil \cr
\noalign{\vskip 4 pt}
\noalign{\hrule}
\noalign{\vskip 1 pt}
2.76 & 2.72 \cr 
 2.57 & 2.71 \cr 
 2.37 & 2.75 \cr 
 2.38 & 2.64 \cr 
 3.17 & 2.91 \cr 
}}}

\sm\noin B) IPBF: $d = 5$

\centerline{\vbox{
\offinterlineskip
\halign{\strut \hfil # & \hfil # &  \hfil #  & ~~#  &  ~~# &  ~~# & ~~#  & ~~#\cr
m & setup & solve & nc& CN  & emax  & rms \cr
\noalign{\hrule}
\noalign{\vskip 1 pt} 
5 & 0.28 & 0.37 & 729 & 5.68e+05 &  4.36e-04 &  9.44e-05 \cr 
6 & 0.54 & 0.57 & 1000 & 1.15e+06 &  3.19e-04 &  6.56e-05 \cr 
7 & 0.95 & 0.84 & 1331 & 2.68e+06 &  2.31e-04 &  4.30e-05 \cr 
8 & 1.61 & 1.27 & 1728 & 5.63e+06 &  1.73e-04 &  3.05e-05 \cr 
9 & 2.79 & 0.38 & 2197 & 1.02e+07 &  1.10e-04 &  2.11e-05 \cr 
10 & 4.25 & 0.60 & 2744 & 2.13e+07 &  8.50e-05 &  1.58e-05 \cr 
}}
\hskip 2pc
\vbox{\offinterlineskip
\halign{\strut
\hfil  \hfil #  &  ~~#&   ~~#   \cr
\omit\span \hfil rates \hfil \cr
\noalign{\vskip 4 pt}
\noalign{\hrule}
\noalign{\vskip 1 pt}
1.40 & 1.63 \cr 
 1.78 & 2.32 \cr 
 1.87 & 2.23 \cr 
 3.40 & 2.76 \cr 
 2.18 & 2.46 \cr 
}}}

\section{Using splines on type-5 tetrahedral partitions}
In this section we examine the use of our methods 
with trivariate polynomial splines on tetrahedral partitions. 
Given a curved domain $\Omega$ in $\RR^3$, we first 
immerse it in the bounding box 
$D:= [ax,bx] \otimes [ay,by] \otimes [az,bz]$ 
surrounding $\Omega$. Then
we create a tetrahedral partition $\tri$ of $D$ as follows.
First we build
a box partition of $D$ using grid lines defined by three
vectors  $ax = x_1 < \cdots < x_\nx = bx$,
$ay = y_1 < \cdots  <y_\ny = by$,
and $az = z_1 < \cdots  <z_\nz = bz$. We then split each sub-box into
five tetrahedra as discussed in Example~6.3 of \cite{Scompb}.
If this is done carefully we get a \tpart{} of $D$ with
$n_t := 5(\nx-1)(\ny-1)(\nz-1)$ tetrahedra. 
We call it a {\sl type-5 tetrahedral partition}.

To approximate the solution of a BVP we will use the space
$$  \Szd := \{s \in C^0(\Omega) \st s|_{T_i} \in \Pd, \quad
                   i = 1,\ldots,\nt\}.  \label{zd} $$
Here $\Pd$ is the space of trivariate polynomials of degree $d$.
The space $\Szd$  is a finite dimensional linear space. Its dimension $N$
is given in Theorem~6.23 of \cite{Scompb}.  Moreover, it has 
a very convenient set of nonnegative, locally supported basis functions $\phi_1,\dots,
\phi_N$ that form a partition of unity.

It turns out that IPBM methods do not perform well with
spaces of splines that
are only continuous.  We need to work with splines that are $C^1$ continuous,
or we need to force the splines to be near $C^1$ continuous.
This can be accomplished by adding the
additional equation 
$$ \lambda_S S c = 0\label{E}$$
to the linear systems in Section~2,
where $\lambda_S$ is a user-selected positive number controlling how close
to $C^1$ we want the spline to be.
Here $S$ is a matrix such that a spline $s \in \Szd$ with coefficient vector $c$
belongs to $C^1$ if and only if $Sc = 0$.  The matrix $S$ is called a 
{\sl smoothness matrix}, 
and can be constructed by the methods described in Section~6.10 of
\cite{Scompb}. 
The algebraic enforcement of spline smoothness has been used in
some related spline methods; see \cite{AwanouLW06}, \cite{LaiL22}, and Remark \ref{AWL}.


The construction of the smoothness matrix $S$ can be avoided if we work with $C^1$
subspaces of $\Szd$. Such spaces are typically called macro-element spaces.
For a comprehensive treatment of them in the trivariate case, see Chapter~18 of
\cite{LaiS}. Computing with them requires the computation of
so-called minimal determining sets and transformation matrices.
Since in the bivariate
case the use of macro-element spaces did not turn out to be worth the trouble, 
see Section~11.2 of \cite{Scompb}, we do not go into them in detail here.

\subsection{Implementing IPBF with splines on tetrahedral partitions}
To implement IPBF with splines on tetrahedral partitions we need to
compute the integrals appearing in (\ref{PBFa}).
These integrals are defined over all of $D$, and
can be computed one tetrahedron at a time using well-known Gaussian quadrature
rules for tetrahedra, see Section~6.6.6 of \cite{Scompb} and Remark~\ref{quadpts}.
Here is a function for computing the coefficient vector \eword{c} of a 
polynomial spline of degree \eword{d} that approximates the solution of the
boundary-value problem~\ref{Prob1}. It uses the IPBF method with
$C^0$ splines defined 
on a type-5 tetrahedral partition of
the bounding box surrounding the curved domain $\Omega$. 

\ms 
\sb{4.3}\vbox{\ob
 [c,H,t1,t2] = ipbf0d3(d,x,y,z,TET,E,F,evertex,fvertex,...
    a1,a2,a3,a4,a5,a6,f,g,xb,yb,zb,S,lambda,lams)
}\eb 
\ms\noin Here \eword{d} is the degree of the polynomial splines,
and \eword{x,y,z} are Cartesian coordinates of the vertices of the partition.
The inputs \eword{TET,..,fvertex} are the lists describing the partition as produced
by the function \eword{tetlistsb} introduced on page 121 of 
\cite{Scompb}. The inputs
{\tt  a1,...,a6,f,g},  describe the BVP and the 
vectors \eword{xb,yb,zb} contain the
Cartesian coordinates of the points in the set $B$.
Finally, \eword{lambda} and \eword{lams} are nonnegative weights associated with
the least-squares problem being solved. The main output is the desired vector \eword{c}
of B-coefficients of the splines, while \eword{H} gives the matrix of 
the linear system being solved
and \eword{t1,t2} are the times to assemble and solve it.

We have also written a simplified version of this function that can be used when the
PDE is the Laplace operator.
\ms  

\sb{4.3}\vbox{\ob 
 [c,H,t1,t2] = ipbf0d3L(d,x,y,z,TET,E,F,evertex,fvertex,...
    f,g,xb,yb,zb,S,lambda,lams)
}\eb 
\ms\noin 

When calling these functions we will use the default
values $\lambda = \lambda_S = .01$.
To test these functions we have written the scripts \eword{T3ipbf} and
\eword{T3ipbfL}. To explore rates of convergence of these methods we
also have the scripts  \eword{T3ipbfr} and
\eword{T3ipbfLr}.

\subsection{Implementing IPBC with splines on tetrahedral partitions}
Let $\Szd$ be a space of $C^0$ splines defined on a tetrahedral
partition of the bounding box of $\Omega$.
To use the IPBC method with this space,
we have to choose a set $\Gamma$ of collocation points lying in the immersing
domain $D$.
Since we are working on tetrahedral partitions, a natural approach to forming $\Gamma$ is
to take it to be the union of the sets
${\cal D}_{m_c,T}$ of domain points of degree $m_c$
(see Remark~\ref{DP})  
associated with each tetrahedron in $\tri$.
Here $m_c$ has to be chosen large enough so that the matrix $H$ has full rank.
Since the equations involve collocating second derivatives of
splines, it is reasonable to conjecture that taking $m_c = d-2$ would be
large enough to make this happen.  We have verified this with several numerical
examples, and so when using IPBC with splines defined on tetrahedral
partitions, we will require $m_c \ge d-2$. The examples in the following
section use values of $d = 4$ and $d = 5$, so we use $m_c = 3$ for these examples.
In the supplement we have presented numerical experiments to
explore what happens for
larger values of  $m_c$. The experiments there show that
the computational times increase somewhat as we increase
$m_c$, but  the  condition numbers and errors do not change too much.

Here is a function that uses the IPBC method to  compute 
the coefficients of a spline in $\Szd$
that approximately solves Problem~\ref{Prob1}.
\ms

\sb{4.2}\vbox{\ob
[c,H,t1,t2] = ipbc0d3(d,x,y,z,TET,E,F,evertex,fvertex,...
      a1,a2,a3,a4,a5,a6,f,g,xb,yb,zb,S,lambda,lams,mc)
}\eb
\ms\noin
The inputs are the same as for \eword{ipbf0d3} except for
the integer \eword{mc} which controls the selection of the collocation points.

We have also written a simplified version of this function that can be used when the
PDE is the Laplace operator.

\ms
\sb{4.3}\vbox{\ob
[c,H,t1,t2] = ipbc0d3L(d,x,y,z,TET,E,F,evertex,fvertex,...
      f,g,xb,yb,zb,S,lambda,lams,mc)
}\eb
\ms\noin

Based on numerical experiments, see the supplement, we recommend calling
these two functions with the default values $\lambda =\lambda_s  = .01$.
To test these functions we have written the scripts \eword{T3ipbc} and
\eword{T3ipbcL}. To explore rates of convergence of these methods we
can employ the scripts  \eword{T3ipbcr} and
\eword{T3ipbcLr}.

\section{Performance of the methods with splines on type-5 tetrahedral partitions}
In this section we give two examples to illustrate
the performance of our methods using trivariate polynomial splines
defined on the immersing domain.
The tables in this section present the same information as in all of the
previous examples.

\begin{ex}\label{ExsphereT}    Solve Problem~\ref{Prob1} with the Laplace operator
on the  sphere shown in \frl{Fig1}
with $f$ and $g$ such that
the true solution is $\sin(5x)\sin(5y)\sin(5z)$.
Use both the IPBF and IPBC methods with trivariate splines $\Szd$
defined on type-5 tetrahedral partitions corresponding to
$m\times m\times m$ grids covering the bounding box of $\Omega$.
\end{ex}

\dis  We use the scripts  \eword{T3ipbfLr} and \eword{T3ipbcLr}, respectively, 
with  the file \eword{sphere50.dat}, 
\eword{pde = 5},  \eword{nb = 2000}. We run these scripts with both
\eword{d = 4} and \eword{d = 5}. 
The errors here are computed on a set of 103,283 points covering $\Omega$.

\sm\noin A) IPBF: $d = 4$

\centerline{\vbox{
\ms
\offinterlineskip
\halign{\strut \hfil # &  \hfil #  & \hfil # \hfil  &  \hfil # \hfil  &  ~~# & ~~#  & ~~#\cr
m & setup & solve & nc& CN  & emax  & rms \cr
\noalign{\hrule}
\noalign{\vskip 1 pt}
 5 & 0.44 & 0.13 &  4273 & 2.01e+03 & 4.08e-03 & 1.03e-03 \cr 
 6 & 0.67 & 0.30 &  8011 & 3.42e+03 & 1.42e-03 & 3.87e-04 \cr 
 7 & 1.67 & 0.86 &  13465 & 5.00e+03 & 7.58e-04 & 1.72e-04 \cr 
 8 & 3.64 & 1.34 &  20959 & 7.22e+03 & 3.96e-04 & 8.91e-05 \cr 
 9 & 9.70 & 2.48 &  30817 & 1.03e+04 & 1.93e-04 & 4.71e-05 \cr 
 10 & 29.37 & 4.57 &  43363 & 1.46e+04 & 1.16e-04 & 2.72e-05 \cr 
}}
\hskip 2pc
\vbox{\offinterlineskip
\halign{\strut
\hfil  \hfil #  &  ~~#&   ~~#   \cr
\omit\span \hfil rates \hfil \cr
\noalign{\vskip 4 pt}
\noalign{\hrule}
\noalign{\vskip 1 pt}
4.73 & 4.37  \cr 
 3.44 & 4.46  \cr 
 4.22 & 4.26  \cr 
 5.38 & 4.78  \cr 
 4.28 & 4.65  \cr 
}}}

\sm\noin B) IPBF: $d = 5$

\centerline{\vbox{
\offinterlineskip
\halign{\strut \hfil # &  \hfil #  & \hfil # \hfil  &  \hfil # \hfil  &  ~~# & ~~#  & ~~#\cr
m & setup & solve & nc& CN  & emax  & rms \cr
\noalign{\hrule}
\noalign{\vskip 1 pt}
 5 & 0.56 & 0.38 &  7981 & 7.31e+03 & 3.46e-04 & 1.04e-04 \cr 
 6 & 1.63 & 0.97 &  15076 & 1.38e+04 & 1.35e-04 & 2.86e-05 \cr 
 7 & 5.70 & 2.16 &  25471 & 2.31e+04 & 4.63e-05 & 1.00e-05 \cr 
 8 & 25.72 & 3.92 &  39796 & 3.88e+04 & 2.07e-05 & 4.07e-06 \cr 
 9 & 86.25 & 7.49 &  58681 & 6.13e+04 & 8.09e-06 & 1.94e-06 \cr 
 10 & 200.21 & 13.17 &  82756 & 9.60e+04 & 4.75e-06 & 9.92e-07 \cr 
}}
\hskip 2pc
\vbox{\offinterlineskip
\halign{\strut
\hfil  \hfil #  &  ~~#&   ~~#   \cr
\omit\span \hfil rates \hfil \cr
\noalign{\vskip 4 pt}
\noalign{\hrule}
\noalign{\vskip 1 pt}
4.22 & 5.78  \cr 
 5.87 & 5.75  \cr 
 5.24 & 5.87  \cr 
 7.02 & 5.53  \cr 
 4.52 & 5.71  \cr 
}}}

\sm\noin C) IPBC: $d = 4$

\centerline{\vbox{
\ms
\offinterlineskip
\halign{\strut \hfil # & \hfil # &  \hfil #  & ~~#  &  ~~# &  ~~# & ~~#  & ~~#\cr
m & setup & solve & nc& CN  & emax  & rms \cr
\noalign{\hrule}
\noalign{\vskip 1 pt}
 5 & 0.19 & 0.10 &  4273 & 9.51e+04 & 1.59e-02 & 3.74e-03 \cr 
 6 & 0.31 & 0.26 &  8011 & 1.16e+05 & 7.22e-03 & 1.50e-03 \cr 
 7 & 0.71 & 0.58 &  13465 & 1.46e+05 & 4.06e-03 & 7.42e-04 \cr 
 8 & 1.32 & 1.22 &  20959 & 1.92e+05 & 2.13e-03 & 4.15e-04 \cr 
 9 & 2.84 & 2.10 &  30817 & 2.45e+05 & 1.26e-03 & 2.46e-04 \cr 
 10 & 5.36 & 4.12 &  43363 & 3.09e+05 & 8.27e-04 & 1.56e-04 \cr 
}}
\hskip 2pc
\vbox{\offinterlineskip
\halign{\strut
\hfil  \hfil #  &  ~~#&   ~~#   \cr
\omit\span \hfil rates \hfil \cr
\noalign{\vskip 4 pt}
\noalign{\hrule}
\noalign{\vskip 1 pt}
3.55 & 4.10  \cr 
 3.15 & 3.85  \cr 
 4.18 & 3.77  \cr 
 3.93 & 3.91  \cr 
 3.58 & 3.88  \cr 
}}}

\sm\noin D) IPBC: $d = 5$

\centerline{\vbox{
\offinterlineskip
\halign{\strut \hfil # & \hfil # &  \hfil #  & ~~#  &  ~~# &  ~~# & ~~#  & ~~#\cr
m & setup & solve & nc& CN  & emax  & rms \cr
\noalign{\hrule}
\noalign{\vskip 1 pt}
 5 & 0.19 & 0.25 &  7981 & 2.75e+05 & 6.95e-03 & 1.70e-03 \cr 
 6 & 0.39 & 0.62 &  15076 & 3.99e+05 & 4.54e-03 & 5.60e-04 \cr 
 7 & 0.95 & 1.44 &  25471 & 5.48e+05 & 1.41e-03 & 2.04e-04 \cr 
 8 & 1.81 & 3.04 &  39796 & 7.74e+05 & 6.74e-04 & 9.16e-05 \cr 
 9 & 3.85 & 6.31 &  58681 & 1.07e+06 & 2.52e-04 & 4.78e-05 \cr 
 10 & 7.54 & 11.96 &  82756 & 1.46e+06 & 1.50e-04 & 2.84e-05 \cr 
}}
\hskip 2pc
\vbox{\offinterlineskip
\halign{\strut
\hfil  \hfil #  &  ~~#&   ~~#   \cr
\omit\span \hfil rates \hfil \cr
\noalign{\vskip 4 pt}
\noalign{\hrule}
\noalign{\vskip 1 pt}
1.91 & 4.97  \cr 
 6.41 & 5.53  \cr 
 4.80 & 5.19  \cr 
 7.36 & 4.88  \cr 
 4.39 & 4.40  \cr 
}}}
\ms

\begin{ex}\label{ExnutT}    Solve Problem~\ref{Prob1} with the Laplace operator 
on the domain shown in \frr{Fig1}
with $f$ and $g$ such that
the true solution is $\sin(5x)\sin(5y)\sin(5z)$.
Use both IPBF and IPBC with splines in the space $\Szd$ on a type-5 tetrahedral partition 
associated with 
$m\times m\times m$ grids covering the bounding box of $\Omega$. 
\end{ex}

\dis   We use the scripts  \eword{T3ipbfLr} and \eword{T3ipbcLr}, respectively, 
with  the file \eword{nut.stl},
\eword{pde = 5}, and \eword{nb = 4000}. We run the scripts for \eword{d = 4} and \eword{d = 5}. 
Here the errors are computed on  a set of 87,025 points covering $\Omega$.

\sm\noin A) IPBF: $d = 4$

\centerline{\vbox{
\offinterlineskip
\halign{\strut \hfil # & \hfil # &  \hfil #  & ~~#  &  ~~# &  ~~# & ~~#  & ~~#\cr
m & setup & solve & nc& CN  & emax  & rms \cr
\noalign{\hrule}
\noalign{\vskip 1 pt}
 5 & 0.28 & 0.14 &  4273 & 1.08e+03 & 1.95e-03 & 3.91e-04 \cr 
 6 & 0.73 & 0.30 &  8011 & 1.65e+03 & 6.97e-04 & 1.50e-04 \cr 
 7 & 1.80 & 0.61 &  13465 & 2.61e+03 & 3.60e-04 & 6.74e-05 \cr 
 8 & 3.83 & 1.32 &  20959 & 3.70e+03 & 2.00e-04 & 3.30e-05 \cr 
 9 & 9.73 & 2.48 &  30817 & 5.14e+03 & 9.49e-05 & 1.80e-05 \cr 
 10 & 26.63 & 4.62 &  43363 & 6.99e+03 & 5.73e-05 & 1.07e-05 \cr 
}}
\hskip 2pc
\vbox{\offinterlineskip
\halign{\strut
\hfil  \hfil #  &  ~~#&   ~~#   \cr
\omit\span \hfil rates \hfil \cr
\noalign{\vskip 4 pt}
\noalign{\hrule}
\noalign{\vskip 1 pt}
4.62 & 4.29  \cr 
 3.63 & 4.39  \cr 
 3.81 & 4.63  \cr 
 5.59 & 4.55  \cr 
 4.28 & 4.44  \cr 
}}}

\sm\noin B) IPBF: $d = 5$

\centerline{\vbox{
\offinterlineskip
\halign{\strut \hfil # & \hfil # &  \hfil #  & ~~#  &  ~~# &  ~~# & ~~#  & ~~#\cr
m & setup & solve & nc& CN  & emax  & rms \cr
\noalign{\hrule}
\noalign{\vskip 1 pt}
 5 & 0.56 & 0.35 &  7981 & 3.39e+03 & 1.98e-04 & 3.64e-05 \cr 
 6 & 1.71 & 0.93 &  15076 & 6.25e+03 & 5.74e-05 & 1.05e-05 \cr 
 7 & 5.35 & 1.86 &  25471 & 1.11e+04 & 2.19e-05 & 3.98e-06 \cr 
 8 & 18.25 & 3.90 &  39796 & 1.79e+04 & 8.39e-06 & 1.64e-06 \cr 
 9 & 52.40 & 7.34 &  58681 & 3.00e+04 & 5.17e-06 & 7.83e-07 \cr 
 10 & 153.59 & 13.25 &  82756 & 4.53e+04 & 2.44e-06 & 4.15e-07 \cr 
}}
\hskip 2pc
\vbox{\offinterlineskip
\halign{\strut
\hfil  \hfil #  &  ~~#&   ~~#   \cr
\omit\span \hfil rates \hfil \cr
\noalign{\vskip 4 pt}
\noalign{\hrule}
\noalign{\vskip 1 pt}
5.56 & 5.58  \cr 
 5.28 & 5.32  \cr 
 6.23 & 5.77  \cr 
 3.62 & 5.53  \cr 
 6.39 & 5.39  \cr 
}}}

\sm\noin C) IPBC: $d = 4$

\centerline{\vbox{
\ms
\offinterlineskip
\halign{\strut \hfil # & \hfil # &  \hfil #  & ~~#  &  ~~# &  ~~# & ~~#  & ~~#\cr
m & setup & solve & nc& CN  & emax  & rms \cr
\noalign{\hrule}
\noalign{\vskip 1 pt}
 5 & 0.18 & 0.11 &  4273 & 4.42e+04 & 1.30e-02 & 2.37e-03 \cr 
 6 & 0.41 & 0.25 &  8011 & 5.87e+04 & 7.33e-03 & 1.11e-03 \cr 
 7 & 0.86 & 0.59 &  13465 & 7.93e+04 & 4.15e-03 & 5.49e-04 \cr 
 8 & 1.49 & 1.21 &  20959 & 1.01e+05 & 2.33e-03 & 2.84e-04 \cr 
 9 & 3.00 & 2.13 &  30817 & 1.25e+05 & 1.13e-03 & 1.55e-04 \cr 
 10 & 5.78 & 4.29 &  43363 & 1.51e+05 & 7.29e-04 & 9.40e-05 \cr 
}}
\hskip 2pc
\vbox{\offinterlineskip
\halign{\strut
\hfil  \hfil #  &  ~~#&   ~~#   \cr
\omit\span \hfil rates \hfil \cr
\noalign{\vskip 4 pt}
\noalign{\hrule}
\noalign{\vskip 1 pt}
2.58 & 3.39  \cr 
 3.12 & 3.88  \cr 
 3.75 & 4.27  \cr 
 5.39 & 4.55  \cr 
 3.75 & 4.24  \cr 
}}}

\sm\noin D) IPBC: $d = 5$

\centerline{\vbox{
\offinterlineskip
\halign{\strut \hfil # & \hfil # &  \hfil #  & ~~#  &  ~~# &  ~~# & ~~#  & ~~#\cr
m & setup & solve & nc& CN  & emax  & rms \cr
\noalign{\hrule}
\noalign{\vskip 1 pt}
 5 & 0.22 & 0.25 &  7981 & 1.43e+05 & 2.30e-03 & 4.28e-04 \cr 
 6 & 0.52 & 0.70 &  15076 & 2.07e+05 & 9.95e-04 & 1.46e-04 \cr 
 7 & 1.07 & 1.47 &  25471 & 2.62e+05 & 2.96e-04 & 5.65e-05 \cr 
 8 & 2.15 & 3.15 &  39796 & 3.79e+05 & 1.63e-04 & 2.49e-05 \cr 
 9 & 4.03 & 6.52 &  58681 & 5.22e+05 & 8.93e-05 & 1.20e-05 \cr 
 10 & 7.70 & 12.20 &  82756 & 6.81e+05 & 4.51e-05 & 6.21e-06 \cr 
}}
\hskip 2pc
\vbox{\offinterlineskip
\halign{\strut
\hfil  \hfil #  &  ~~#&   ~~#   \cr
\omit\span \hfil rates \hfil \cr
\noalign{\vskip 4 pt}
\noalign{\hrule}
\noalign{\vskip 1 pt}
3.75 & 4.81  \cr 
 6.65 & 5.21  \cr 
 3.87 & 5.32  \cr 
 4.51 & 5.49  \cr 
 5.80 & 5.56  \cr 
}}}
\ms\noin
These two examples should be compared with Examples~\ref{Exsphere} and
\ref{Exnut} where the same problems were solved using {\tp} splines.
Here we are getting comparable accuracies, but with many more coefficients and 
significantly more computational time.  In both examples we see
that the IPBC method has larger condition numbers than the IPBF method, 
which is the opposite of what happens with \tp{} splines.
In addition, here the rates are less stable, but 
as with \tp{} splines, it appears that IPBF gives near
optimal rates while IPBC is at least one lower.

\begin{figure}
\centering
\subfloat{\includegraphics[width = .5\linewidth,trim=0 70 0 70, clip]{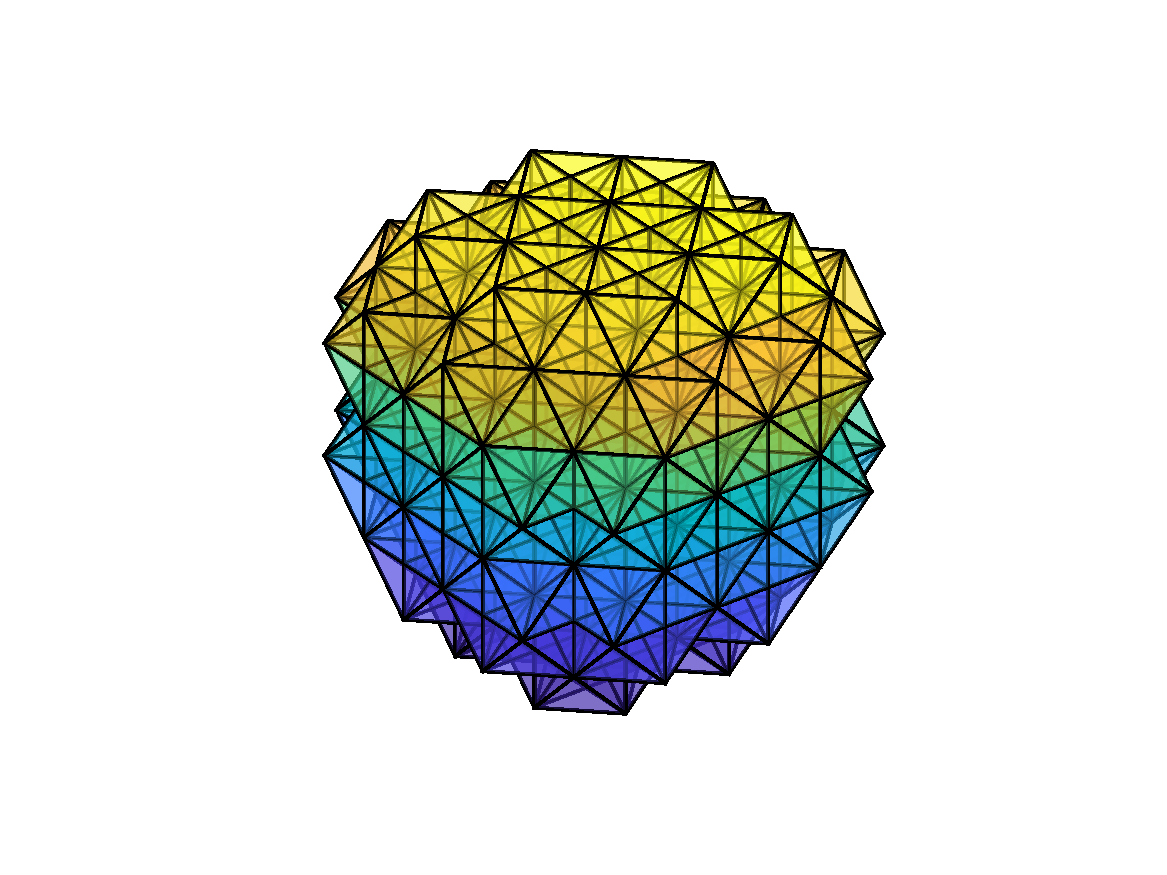}}
\hspace{0.09\linewidth}
\raisebox{0em}{\subfloat{\includegraphics[width = .28\linewidth,trim=0 0 0 0, clip]{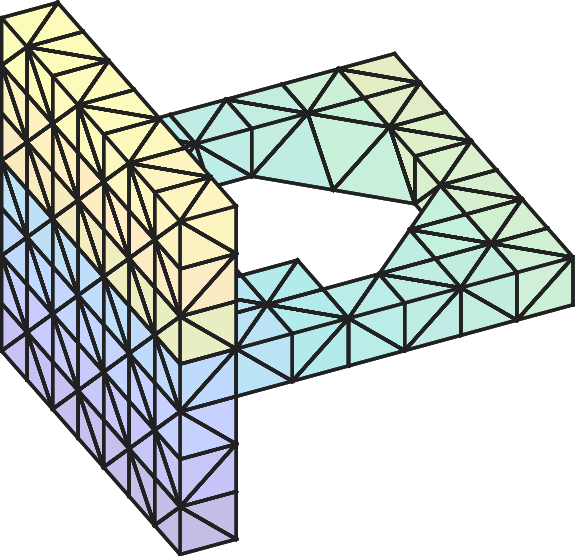}}}
\caption{Trimmed tetrahedral immersing domains for the sphere and bracket}
\vskip -.1in
\label{Fig6} 
\end{figure}

\section{Trimming tetrahedral partitions}
In working with splines on tetrahedral partitions, we can make our methods
more efficient by reducing the size of the immersing domain $D$. When using type-5
partitions covering the entire bounding box, it frequently happens
that some of the tetrahedra do not intersect $\Omega$. 
Our package includes the  following
function for removing such tetrahedra  from the immersing domain.
It makes use of the large set of points {\tt xr,yr,zr} covering the interior
of the domain created by {\tt inpts3}.

\ms 
\sb{3.8}\vbox{\ob
[x,y,z,TET,tetout] = trimtet(xo,yo,zo,TETo,xr,yr,zr)
}\eb
\ms\noin 
The inputs {\tt xo,yo,zo,TETo} 
describe the original tetrahedral partition while the outputs
{\tt x,y,z,TET} describe the corresponding trimmed tetrahedral partition.
The vector {\tt tetout} gives a list of the tetrahedra that were removed.

\begin{ex}\label{Exspheretrim}    Solve Problem~\ref{Prob1} on the spherical domain
shown in \frl{Fig1} with the Laplace operator. Choose the 
true solution to be
$\sin(8x) + \sin(12y) + \sin(14z)$, and 
use the IPBF method with splines of degree six on a type-5 tetrahedral
partition of the bounding box. Redo the problem after trimming this partition.

\end{ex}

\dis  First we use the script \eword{T3ipbfLr} with \eword{sphere50.dat},
\eword{pde = 6}, \eword{nb = 2000}, and  \eword{d = 6}. 
Then we rerun the problem with the script \eword{T3ftrimLr}.
The columns labelled {\tt nto} and {\tt nt} in B)
give the original number of tetrahedra, and the number after trimming,
respectively.
\ms

\noin A) untrimmed

\centerline{\vbox{
\offinterlineskip
\halign{\strut \hfil # & \hfil # &  \hfil #  & ~~#  &  ~~# &  ~~# & ~~#  & ~~#\cr
m & setup & solve & nc& CN  & emax  & rms \cr
\noalign{\hrule}
\noalign{\vskip 1 pt}
 5 & 1.08 & 0.86 &  13385 & 3.20e+04 & 9.08e-04 & 2.07e-04 \cr 
 6 & 7.14 & 2.45 &  25416 & 7.70e+04 & 2.60e-04 & 5.16e-05 \cr 
 7 & 22.40 & 4.96 &  43093 & 1.52e+05 & 6.95e-05 & 1.46e-05 \cr 
 8 & 89.17 & 11.02 &  67502 & 2.78e+05 & 2.28e-05 & 5.14e-06 \cr 
}}
\hskip 2pc
\vbox{\offinterlineskip
\halign{\strut
\hfil  \hfil #  &  ~~#&   ~~#   \cr
\omit\span \hfil rates \hfil \cr
\noalign{\vskip 4 pt}
\noalign{\hrule}
\noalign{\vskip 1 pt}
5.60 & 6.23  \cr 
 7.24 & 6.93  \cr 
 7.23 & 6.75  \cr 
}}}
\ms

\noin B) trimmed

\centerline{\vbox{
\offinterlineskip
\halign{\strut \hfil # &  \hfil #  & \hfil # \hfil  &  \hfil # \hfil  
& \hfil # \hfil & \hfil # \hfil  &  \hfil # \hfil  &  \hfil # \hfil  &  \hfil # \hfil  \cr
m & setup & solve & nto&  nt & nc & CN  & emax  & rms \cr
\noalign{\hrule}
\noalign{\vskip 1 pt}
 5 & 0.93 & 0.92 & 320 & 288 & 11845 &  8.19e+03 & 9.14e-04 & 2.06e-04 \cr 
 6 & 3.16 & 1.73 & 625 & 501 & 20020 &  6.92e+03 & 2.34e-04 & 4.90e-05 \cr 
 7 & 16.27 & 4.19 & 1080 & 872 & 34613 &  1.03e+04 & 7.11e-05 & 1.41e-05 \cr 
 8 & 43.19 & 6.34 & 1715 & 1239 & 48406 &  8.57e+03 & 2.16e-05 & 4.77e-06 \cr 
}}
\hskip 2pc
\vbox{\offinterlineskip
\halign{\strut
\hfil  \hfil #  &  ~~#&   ~~#   \cr
\omit\span \hfil rates \hfil \cr
\noalign{\vskip 4 pt}
\noalign{\hrule}
\noalign{\vskip 1 pt}
6.10 & 6.44  \cr 
 6.54 & 6.84  \cr 
 7.74 & 7.03  \cr 
}}}
\ms\noin \frl{Fig6} shows the trimmed tetrahedral partition
corresponding to $m = 8$.
The tables show that after trimming
the computational times and
numbers of coefficients are much smaller, while
the accuracy is almost the same.  

Here is a second example where trimming is even more effective.

\begin{ex}\label{Exbrackettrim}    Solve Problem~\ref{Prob1} on the domain
shown in \frr{Fig2}
with true solution
$\sin(5x) + \sin(5y) + \sin(5z)$.  Suppose the coefficients of $L$
are  given by\break
$a = (10+\cos(x) ,e^y, 10+\sin(x),\cos(x),\cos(y),\cos(z))$.
Use the IPBF method with splines of degree five on a type-5 tetrahedral
partition of the bounding box. Redo the problem after trimming this partition.

\end{ex}

\dis  First we use the script \eword{T3ipbfr} with \eword{BracketWithHole.stl},
\eword{pde = 11}, \eword{nb = 2000}, and  \eword{d = 5}. 
Then we rerun the problem with the script \eword{T3ftrimr}.
The columns labelled {\tt nto} and {\tt nt} in the resulting table
give the original number of tetrahedra, and the number after trimming,
respectively.
\ms

\noin A) untrimmed

\centerline{\vbox{
\offinterlineskip
\halign{\strut \hfil # & \hfil # &  \hfil #  & ~~#  &  ~~# &  ~~# & ~~#  & ~~#\cr
m & setup & solve & nc& CN  & emax  & rms \cr
\noalign{\hrule}
\noalign{\vskip 1 pt}
 5 & 0.63 & 0.33 &  7981 & 1.10e+05 & 8.41e-05 & 1.37e-05 \cr 
 6 & 1.76 & 0.91 &  15076 & 1.70e+05 & 2.46e-05 & 3.58e-06 \cr 
 7 & 4.72 & 1.74 &  25471 & 2.72e+05 & 1.10e-05 & 1.31e-06 \cr 
 8 & 21.61 & 4.02 &  39796 & 4.79e+05 & 3.99e-06 & 5.13e-07 \cr 
}}
\hskip 2pc
\vbox{\offinterlineskip
\halign{\strut
\hfil  \hfil #  &  ~~#&   ~~#   \cr
\omit\span \hfil rates \hfil \cr
\noalign{\vskip 4 pt}
\noalign{\hrule}
\noalign{\vskip 1 pt}
5.50 & 6.00  \cr 
 4.42 & 5.50  \cr 
 6.58 & 6.11  \cr 
}}}
\ms

\noin B) trimmed

\centerline{\vbox{
\offinterlineskip
\halign{\strut \hfil # &  \hfil #  & \hfil # \hfil  &  \hfil # \hfil  
& \hfil # \hfil & \hfil # \hfil  &  \hfil # \hfil  &  \hfil # \hfil  &  \hfil # \hfil  \cr
m & setup & solve & nto&  nt & nc & CN  & emax  & rms \cr
\noalign{\hrule}
\noalign{\vskip 1 pt}
 5 & 0.22 & 0.12 &  320 & 147 & 4230 &  2.46e+04 & 7.59e-05 & 1.47e-05 \cr 
 6 & 0.29 & 0.14 &  625 & 189 & 5490 &  1.34e+04 & 2.69e-05 & 4.14e-06 \cr 
 7 & 0.63 & 0.27 &  1080 & 337 & 9305 &  3.48e+04 & 1.02e-05 & 1.37e-06 \cr 
 8 & 0.75 & 0.30 &  1715 & 370 & 10430 &  2.85e+04 & 3.71e-06 & 6.30e-07 \cr 
}}
\hskip 2pc
\vbox{\offinterlineskip
\halign{\strut
\hfil  \hfil #  &  ~~#&   ~~#   \cr
\omit\span \hfil rates \hfil \cr
\noalign{\vskip 4 pt}
\noalign{\hrule}
\noalign{\vskip 1 pt}
4.65 & 5.67  \cr 
 5.33 & 6.05  \cr 
 6.55 & 5.06  \cr 
}}}
\ms\noin \frr{Fig6} shows the trimmed tetrahedral partition
corresponding to $m = 8$. As in the previous example
the computational times and numbers of coefficients  after trimming
are much smaller, while
the accuracy is almost the same.  
We see that using IPBF with trimmed trivariate splines
gives a given accuracy  in even less time than either of the 
methods based on \tp{} splines discussed in
Example~\ref{Exbracket}.

\section{Conclusions}
Compared with other methods in the literature, the IPBM method
has several advantages: namely,

\begin{itemize}
\item there is no need to construct a mapping 
onto the desired domain $\Omega$ from some
computational domain such as a collection of patches,
\item there is no need to mesh complicated curved domains,
\item it is not required to construct basis functions that vanish on the boundary,
\item it is not required to evaluate integrals over curved sets,
\item the approximating spaces do not involve rational functions.
\end{itemize}

\ms

When using splines (piecewise polynomials) to solve boundary-value problems,
we would like to have theorems that give the exact 
relationship between the mesh sizes of the spaces and the approximation errors. 
No such theorems exist for IPBM methods as of this writing, which is why we
have included detailed information on rates of convergence in all of our examples.
Ideally, we would like the methods to produce optimal order convergence, i.e.,
the rates of convergence in solving a BVP are the same as the rates of convergence
when simply approximating a known function,
see Remarks~\ref{appxTP} and \ref{Appxtet}.
Our tables show that when using the IPBF method,
the rates of convergence of the RMS
errors are near optimal. For IPBC, they appear to be suboptimal.


The examples presented here show that working with \tp{} splines is
generally much more efficient than working with splines on tetrahedral
partitions.  Tensor-product splines also have the advantage that they
are inherently very smooth -- the smoothness in each direction is of order
one less than the degree in that direction. There are, however, certain
situations where working with splines on tetrahedral partitions may be advantageous,
for example when we use trimmed partitions as was done in the examples
of Section~8.

Comparing  IPBF with IPBC, we note that IPBF seems to have better convergence rates
than IPBC, but to produce a given accuracy, IPBC may often do it in less time than IPBF.
This makes IPBC attractive for low accuracies, but
when high accuracy is required,  IPBF  might be a better choice. 
Another way of getting high accuracy is to 
work with  splines of higher degrees, see Section~SM9 of the supplement. 
Our codes can be used for degrees up to nine.

\section{Remarks} 

\begin{remark} \label{immerse}
\rm  There are many papers in the literature which solve BVPs with piecewise
polynomials defined on a mesh covering an immersing domain.  We do not attempt
to list them here, but  see \cite{Burman15} or \cite{Burman25}. The
latter paper summarizes many of these methods, including 
the immersed boundary method, the fictitious domain method, the boundary penalty
method, the volume penalty method, the immersed finite element method,
the finite cell method, and
the cut finite element method. The paper includes  several hundred references.
While our methods  share the idea of an immersing domain, 
our approach to handling boundary conditions is very different. Here
we simply collocate the boundary conditions instead of having to compute
integrals over curved sets.
\end{remark}

\begin{remark} \label{AWL}
\rm  There are a number of papers in the literature that solve BVPs on a 3D
polygonal domain $\Omega$
using polynomial spline spaces defined on tetrahedral partitions of $\Omega$.
\cite{AwanouLW06} is such a paper
where the spline coefficients are determined by minimizing a certain quadratic 
functional  subject to side constraints enforcing smoothness of the spline and 
making it approximately fit the boundary conditions at points on the boundary. 
\cite{LaiL22} describes a method based on minimizing a quadratic functional based on
boundary conditions and smoothness, where now the side constraint involves
collocation of the PDE at points inside the domain.
Neither paper uses immersing domains. 
Their methods require solving constrained
minimization problems, while here we simply solve a least-squares minimization problem.
\end{remark}

\begin{remark} \label{IGA}
\rm A very different approach to solving BVPs on curved domains is to
first find 
a mapping from a simpler computational domain that covers the desired domain.
The heavily studied IGA methods, see \cite{Cotrell}, fall into this
category. As pointed out in the literature, the construction of such maps
(especially for domains with holes) can be difficult. To get an idea of
current topics along with many references, see  \cite{Evans24}
and the rest of that special issue on Isogeometric Analysis which contains
13 papers. Topics 
include immersed IGA, Boundary-Element IGA, and isogeometric collocation methods.
For some recent papers on isogeometric collocation, see \cite{Kapl26} and
\cite{Shill13} and references therein.
\end{remark}

\begin{remark} \label{RBF}
\rm  So-called radial basis functions (RBFs) provide an interesting
alternative to splines for approximation purposes, and there is a growing
literature on their use for solving boundary value problems.  For an
introduction to collocation with RBFs for BVPs 
in the 2D case, see  Chapter 38 of \cite{Fass07}. It contains
several examples and quite a few references. For a recent paper
that works with 3D domains, see \cite{Chen20} where basis functions defined
on an immersing domain are used.  RBF methods always involve a parameter
controlling the shape of the basis functions and the authors present a heuristic
for choosing it. This parameter has a major
impact on the accuracy of the methods, and choosing a good value for it  is critical.
\end{remark}

\begin{remark} \label{extension}
\rm As observed in Sect.~1, to use the IPBM methods discussed here we need
the coefficient functions of $L$ and the 
function $f$ in (\ref{Lu}) to
be defined on the immersing domain $D$, not just on $\Omega$.  
We believe that for the vast majority of practical problems, they will 
be typically defined on all of $\RR^3$. 
\end{remark}

\begin{remark}\label{2D}
\rm An extensive set of  examples showing the efficacy
of immersed penalized boundary  methods 
in the bivariate case can be found in Chapters~10--14 of
\cite{Scompb}, where the method was
applied to solve BVPs on curved domains in $\RR^2$  using
bivariate tensor-product splines as well as polynomial
splines on various partitions of the immersing domain, including ordinary
triangulations, curved triangulations, triangulations with hanging vertices,
 and T-meshes.
\end{remark}

\begin{remark}\label{scale}
\rm To make sure the domains in the examples in this paper have comparable
dimensions, before solving a boundary-value problem for a given domain, 
we first scale it so that the
corresponding bounding box has largest side length equal to  one.
This is important for comparison purposes since a true solution involving high frequency
sine functions would look very different on a box with maximum side length 50
as compared with a box with maximum side length one.
\end{remark}

\begin{remark}\label{truebdy}
\rm 
To use the method described here the user must provide a description of
the domain of interest. Our codes do not make any change to this model
unlike parametric methods where for example the boundary of the domain is
often modified. For our methods the boundary conditions are enforced
only at points on the user's original boundary.
\end{remark}

\begin{remark} \label{times}
\rm Each example in this paper required setting up 
a linear system of equations and finding a least-squares solution of it.
No attempt was made to preprocess these systems --
we simply used the Matlab backslash operator to solve them.
In most examples we have reported the run times for both the setup and
for the solution of the system.  We include these just to give a general idea
of how fast the methods are for different domains and different spline spaces.
\end{remark} 

\begin{remark} \label{backslash}
\rm  As we saw in Section~2, the IPBF method always leads to an overdetermined
linear system of equations which we then solve by least squares using either
the backslash operator directly on $H$, or by converting it to a square system
with the Gram matrix $G$.  For the examples in this paper the two methods give
almost identical results. Using $G$ is somewhat faster, but solve times are already
quite small. In situations where
we need to solve the systems to high accuracy, such as
when checking the patch test where we expect
essentially zero errors, it turns out that working with the Gram matrix $G$
is not appropriate due to the well-known fact that Gram matrices are 
notoriously ill conditioned. In fact, the condition numbers of $G$ are typically
approximately the square of the condition numbers of the associated
$H$ matrix.  We should also note that if the matrix $H$ is not of full rank,
the Gram matrix $G$ will be singular. This can happen 
for example with IPBC if we do not use enough  collocation points, see
Example~SM3.2 in the supplement.  
\end{remark}

\begin{remark} \label{CN}
\rm  When solving an overdetermined system of equations defined by a 
sparse matrix $H$ with the Matlab backslash operator, it would be useful to
have an estimate for the 2-norm condition number $CN_2$ of $H$. 
It  can be computed exactly as {\tt cond(full(H))}. 
However, for problems where $H$ is very large, 
this computation can take a considerable
amount of time, or even run out of memory, as happens for
the matrices arising in Example~\ref{ExsphereT} for larger values of $m$.
To avoid these problems we can
make use of the fact that $CN_2(H) = \sqrt{{CN_2(G)}}$ where
$G$ is  the Gram matrix $G = H^TH$, coupled with the fact that
an estimate for the  1-norm condition number $CN_1$ of $G$ can be quickly computed
with {\tt condest(G)}, and it is typically quite close to $CN_2(G)$. 
Thus, in all of the examples in the paper we report 
$CN  = \sqrt{CN_1(H^T H)}$ as our estimate of $CN_2(H)$.
\end{remark}

\begin{remark} \label{ellip}
\rm A second order partial differential operator $L$ as in (\ref{L1})
is defined to be {\sl elliptic} on a domain $\Omega$ provided 
that the eigenvalues of the matrix
\begin{equation}
E = \begin{pmatrix}
    a_1 & a_4 & a_5 \\  
    a_4 &  a_2 & a_6 \\ 
    a_5 & a_6 & a_3 \\
  \end{pmatrix}
\end{equation}
are positive for all points in the domain $\Omega$.
The Laplace operator where $a_1 = a_2 = a_3 = 1$ and $a_4 = a_5 = a_6 = 0$
is elliptic since in this case the eigenvalues are just $(1,1,1)$.
All of the PDEs used in this paper are elliptic in this sense.
Surprisingly, we have found that IPBM methods can give good results even for
non-elliptic PDEs, see Example~SM7.1 in the supplement.
\end{remark}

\begin{remark}\label{spacing}
\rm For the examples presented here we have required that the points in the set $B$ 
used in the penalty term for the
Dirichlet boundary conditions be reasonably well-spaced. To show that 
this is essential, we have
redone Example~\ref{Exsphere} using point sets $B$ where all the points lie
on the bottom half of the sphere.  This gives very large errors.
\end{remark}

\begin{remark}\label{quadpts}
\rm For a discussion of how to do numerical quadrature over a 
tetrahedron, see Section~6.6.6 of \cite{Scompb}. The function 
{\tt quadpts3} described there can be used to read in quadrature points and weights
to get formulae that are exact for trivariate polynomials of degree $m_q$ for
all choices of even $m_q$ from 4 to 20. For $m_q = 10$ and $12$ the numbers of points
produced are 81 and 168, respectively.
\end{remark}


\begin{remark}\label{aniso}
\rm When working with \tp{} splines, throughout this paper we have 
chosen the numbers of knots in all three directions to be the same, and have
worked only with equally-spaced knots.  However, our codes can be used
with different numbers of knots in the three directions, and with
arbitrarily spaced knots. This could be 
useful for PDEs where the solution is anisotropic, i.e., behaves rather
differently in the three variables.
\end{remark}

\begin{remark}\label{DP}
\rm  Given a tetrahedron $T$ with the four vertices $v_1,\ldots,v_4$, 
and a positive integer
$d$, the set   ${\cal D}_{d,T}$  
of {\sl domain points}  associated with $T$
is just the set of points of the form
$\xi_{ijkl} : = (iv_1 + jv_2 + kv_3 + lv_4)/d$ for $i+j+k+l = d$.  These
are a total of  $({d+3 \atop 3})$
equally-spaced points distributed throughout $T$,
and are very useful in working with trivariate polynomials defined on $T$,
see \cite{Scomp} and \cite{Scompb}.
\end{remark}

\begin{remark} \label{Rates}
\rm All of our examples include estimates of the rates of convergence
for both the max and RMS errors. Based on what is known about the finite-element
methods, we expect these errors to behave like $\cO(h^r)$ where $h$ is a
measure of the mesh size of the spline space and $r$ is a positive number.
To estimate these rates of convergence we can run an example 
with spline spaces with two different mesh sizes $h_2< h_1$, and then 
compute $\log_2(e_2/e_1)/\log_2(h_2/h_1)$, where $e_2$ and $e_1$ are the errors
of interest.  In the FEM literature errors are 
usually computed in some $H^1$ norm, but here we have opted
to present max and RMS errors. Moreover, it is common to draw graphs to
illustrate rates of convergence, but we believe the tables presented here
provide more detailed information.
\end{remark}

\begin{remark}\label{appxTP}
\rm 
It is known  that when approximating a smooth function $F$ by
\tp{} splines of degree $(d_x,d_y,d_z)$
with mesh sizes $h_x,h_y,h_z$ the error in the uniform norm
has the form
$\cO (h_x^{dx+1})$ +  $\cO (h_y^{dy+1})$ +  $\cO (h_z^{dz+1})$. 
This is called the optimal rate of approximation as it cannot be improved
no matter how smooth $F$ might be.
\end{remark}

\begin{remark} \label{Appxtet}
\rm It is known that when approximating a function 
$F \in C^r(\Omega)$  with $0 \le r \le d$
by  a trivariate polynomial spline on a tetrahedral partition with mesh size $h$
the error in the uniform norm behaves like ${\cal O}(h^{r+1})$.  
The constant in this error bound depends on the shape parameter of the partition,
and the highest rate achievable 
with splines of degree $d$ is $d+1$ which is called
the {\sl optimal order}.
\end{remark}

\begin{remark}\label{patch}
\rm For a computational method for solving a BVP using splines
to have optimal order convergence,
it has to give exact results when the true solution is a polynomial
in the spline space. Engineers sometimes refer to this as the patch test.
Our IPBM methods are set up to give exact results not only for polynomials
but also in the case where the true solution is a spline in the space
being used.  For some examples see sections SM5 and SM6 of the supplement.
\end{remark}

\begin{remark}\label{mixed}
\rm In the 2D case we gave a detailed treatment of the use of the immersed
penalized boundary method for second order problems with mixed boundary
conditions, see Sects.~10.3, 11.4, and 12.3 in \cite{Scompb}.
We saw that the results with mixed boundary conditions were
very similar to those with only Dirichlet boundary conditions.  
Our IPBM methods can be easily adapted to
work with mixed boundary conditions also in the  3D case.
\end{remark}

\begin{remark}\label{curved tets}
\rm It is possible to use both the IPBF and IPBC methods on curved
domains without using an immersing domain by creating 
partitions of $\Omega$ that are allowed to have curved elements.
To see how this works in the planar case, see \cite{Scurve} and
Chapter~12 of \cite{Scompb}.
\end{remark}

\begin{remark} \label{sources}
\rm Our experiments with the sphere are based on a
cloud of 50,000 points lying exactly on the surface of the sphere.
The files  ForearmLink.stl, Torus.stl, and BracketWithHole.stl
are in the Matlab PDE toolbox. The files
 teapot.stl  and nut.stl can be downloaded from 
https://www.thingiverse.com/thing:852078 and\break
https://3dexport.com/3dmodel-hexagon-slotted-nut-din-935-m8-394076.htm,
respectively. 
The files roof.dat and tube.dat were provided
by Dr. Emily Johnson of Notre Dame University.
Our computations with the Stanford bunny were based on
the file bun\_zipper.ply, see
https://graphics.stanford.edu/pub/3Dscanrep/bunny.tar.gz.
\end{remark}

\begin{remark}\label{spak2}
\rm All of the examples in this paper and its supplement
were computed with Matlab on a standard desktop. All of the software
can be downloaded from
https://github.com/AussieGreene/splinepak2.
For convenience, this software package also includes  Matlab code for several functions
described in \cite{Scompb}
for dealing with trivariate tensor-product splines and with trivariate polynomial
splines on tetrahedral partitions.
Note that when working with stl files, our function \eword{setdomain3} calls on
the built-in Matlab function \eword{mesh2pc} which
is not part of a basic Matlab installation. To get access
to it the user must also install the Matlab lidar toolbox.
\end{remark}

\bibliographystyle{siamplain} 
\bibliography{references}
\end{document}


\maketitle

\section{Introduction}
This document is a supplement to the associated SISC paper {\sl Solving PDEs on Curved
Domains in 3D using the Immersed Penalized Boundary Method with Trivariate Spline Spaces}.
In it we will explore various aspects of its use to better understand how it works.
These include choice of the set B and parameter $\lambda$, choice of collocation
points for IPBC, choice of the parameter $\lambda_s$ when using $C^0$ splines,
the patch test and exactness when the true solution is a spline,
nonelliptic PDEs, the calculation of interior points for error calculations,
and a comparison of refinement strategies. 

\section{Choice of the set \texorpdfstring{$B$}{B} and the parameter \texorpdfstring{$\lambda$}{lambda}}
Both the IPBF method and IPBC method require the user to 
choose a set of $n_b$ well-spaced points $B$
on the boundary of $\Omega$ along with a parameter $\lambda$. 
In this section
we give several examples to illustrate the interplay between the size of
$\nb$  and the choice of the parameter $\lambda$.
The first three examples work with \tp{} splines.

\begin{ex}\label{SM21}    
Solve Problem~\ref{Prob1} with the Laplace operator $L$ on the 
sphere shown in \frl{Fig1}, 
where $f$ and $g$ are chosen so that
the true solution is $\sin(5x)+\sin(5y)+\sin(5z)$.
Use the IPBF method with tensor-product splines of degree (5,5,5)
defined on a $7\times 7\times 7$ grid covering the bounding box of $\Omega$.
Explore the errors  for a sequence of $\lambda$ values
using 2,000, 5,000, and 10,000 boundary points.
\end{ex}

\dis We use the script \eword{TP3flam} with  \eword{sphere50.dat},
\eword{pde = 4}, \eword{nb = 2000},  \eword{dx = 5}, \eword{mx = 7},
and set the starting value of $\lambda$ to .001.
The following table shows results for a sequence of values of $\lambda$.
For each $\lambda$ it shows
the time to assemble the system and solve it, the number nc of coefficients,
an estimate CN for the 2-condition number of the matrix $H$ of the linear
system being solved, and the Max and RMS errors on
a set of 103,283
points covering $\Omega$ and its boundary.

A) nb = 2000

\centerline{\vbox{
\offinterlineskip
\halign{\strut \hfil # &  \hfil #  & \hfil # \hfil  &  ~~# &  ~~# & ~~#  & ~~#\cr
$\lambda$ & time & nc& CN  & emax  & rms \cr
\noalign{\hrule}
\noalign{\vskip 1 pt}
.001 & 1.76 & 1331 & 2.25e+06 &  5.98e-04 &  9.96e-05 \cr 
.01 & 1.72 & 1331 & 5.88e+05 &  6.84e-05 &  1.60e-05 \cr 
.1 & 1.71 & 1331 & 2.33e+05 &  1.26e-04 &  1.73e-05 \cr 
1 & 1.75 & 1331 & 4.59e+05 &  2.05e-04 &  2.50e-05 \cr 
}}}

B) nb = 5000

\centerline{\vbox{
\offinterlineskip
\halign{\strut \hfil # &  \hfil #  & \hfil # \hfil  &  ~~# &  ~~# & ~~#  & ~~#\cr
$\lambda$ & time & nc& CN  & emax  & rms \cr
\noalign{\hrule}
\noalign{\vskip 1 pt}
.001 & 2.41 & 1331 & 1.74e+06 &  3.19e-04 &  5.96e-05 \cr 
.01 & 2.36 & 1331 & 4.23e+05 &  6.73e-05 &  1.46e-05 \cr 
.1 & 2.37 & 1331 & 2.12e+05 &  1.40e-04 &  1.90e-05 \cr 
1 & 2.43 & 1331 & 5.90e+05 &  2.17e-04 &  2.62e-05 \cr 
}}}

C) nb = 10000

\centerline{\vbox{
\offinterlineskip
\halign{\strut \hfil # &  \hfil #  & \hfil # \hfil  &  ~~# &  ~~# & ~~#  & ~~#\cr
$\lambda$ & time & nc& CN  & emax  & rms \cr
\noalign{\hrule}
\noalign{\vskip 1 pt}
.001 & 3.85 & 1331 & 1.45e+06 &  1.87e-04 &  4.01e-05 \cr 
.01 & 4.15 & 1331 & 3.32e+05 &  7.45e-05 &  1.44e-05 \cr 
.1 & 3.69 & 1331 & 2.14e+05 &  1.53e-04 &  2.05e-05 \cr 
1 & 3.62 & 1331 & 6.94e+05 &  2.35e-04 &  2.84e-05 \cr 
}}}
\ms
\noin Using more points increases the computational times with no
essential gain in accuracy.  In all three cases 
the RMS errors and the condition numbers go down and then up, and their
minima happen for $\lambda = .01$. \eop

Here is a second example  with  a more complicated domain and a more complicated
true solution. It shows the same behavior.

\begin{ex}\label{SM22}    Repeat Example~\ref{SM21} with the sphere replaced
by the domain shown in \frl{Fig2} and with the true solution
given by $u = \sin(8x) + \sin(12y) + \sin(14z)$.
\end{ex}

\dis   We use the script \eword{TP3flam} with  \eword{ForearmLink.stl},
\eword{pde = 6}, \eword{nb = 2000}, \eword{dx = 5}, and \eword{mx = 7}.
Here the error is measured on 92,477 points.

A) nb = 2000

\centerline{\vbox{
\offinterlineskip
\halign{\strut \hfil # &  \hfil #  & \hfil # \hfil  &  ~~# &  ~~# & ~~#  & ~~#\cr
$\lambda$ & time & nc& CN  & emax  & rms \cr
\noalign{\hrule}
\noalign{\vskip 1 pt}
.001 & 1.77 & 1331 & 5.62e+06 &  2.53e-03 &  3.06e-04 \cr 
.01 & 1.76 & 1331 & 1.75e+06 &  8.90e-04 &  2.29e-04 \cr 
.1 & 1.73 & 1331 & 5.59e+05 &  1.57e-03 &  3.14e-04 \cr 
1 & 1.71 & 1331 & 7.51e+05 &  2.75e-03 &  4.20e-04 \cr 
}}}

B) nb = 5000

\centerline{\vbox{
\offinterlineskip
\halign{\strut \hfil # &  \hfil #  & \hfil # \hfil  &  ~~# &  ~~# & ~~#  & ~~#\cr
$\lambda$ & time & nc& CN  & emax  & rms \cr
\noalign{\hrule}
\noalign{\vskip 1 pt}
.001 & 2.39 & 1331 & 4.32e+06 &  1.57e-03 &  2.60e-04 \cr 
.01 & 2.36 & 1331 & 1.31e+06 &  8.08e-04 &  2.37e-04 \cr 
.1 & 2.39 & 1331 & 4.25e+05 &  1.67e-03 &  3.37e-04 \cr 
1 & 2.35 & 1331 & 1.01e+06 &  2.77e-03 &  4.40e-04 \cr 
}}}

C) nb = 10000

\centerline{\vbox{
\offinterlineskip
\halign{\strut \hfil # &  \hfil #  & \hfil # \hfil  &  ~~# &  ~~# & ~~#  & ~~#\cr
$\lambda$ & time & nc& CN  & emax  & rms \cr
\noalign{\hrule}
\noalign{\vskip 1 pt}
.001 & 3.55 & 1331 & 3.54e+06 &  1.18e-03 &  2.41e-04 \cr 
.01 & 3.45 & 1331 & 1.08e+06 &  8.76e-04 &  2.47e-04 \cr 
.1 & 3.52 & 1331 & 4.00e+05 &  1.65e-03 &  3.52e-04 \cr 
1 & 3.47 & 1331 & 1.30e+06 &  3.52e-03 &  4.86e-04 \cr 
}}}
\ms

Our next example makes use of the IPBC method with \tp{} splines.

\begin{ex}\label{SM23}    
Solve Problem~\ref{Prob1} with the Laplace operator $L$ on the
domain shown in \frl{Fig5}, 
where $f$ and $g$ are such that
the true solution is $\sin(8x)+\sin(12y)+\sin(14z)$.
Use the IPBC method with tensor-product splines of degree (5,5,5)
defined on a $7\times 7\times 7$ grid covering the bounding box of $\Omega$.
Explore the errors  for a sequence of $\lambda$ values
using 2,000, 5,000, and 10,000 boundary points.
\end{ex}

\dis   We use the script \eword{TP3clam} with  \eword{bunny.ply},
\eword{pde = 6}, \eword{nb = 2000}, \eword{mc = 3}, \eword{dx = 5}, and \eword{mx = 7}.
Here the errors are measured on 56,948 points inside the domain.

A) nb = 2000

\centerline{\vbox{
\offinterlineskip
\halign{\strut \hfil # &  \hfil #  & \hfil # \hfil  &  ~~# &  ~~# & ~~#  & ~~#\cr
$\lambda$ & time & nc& CN  & emax  & rms \cr
\noalign{\hrule}
\noalign{\vskip 1 pt}
.001 & 1.19 & 1331 & 6.82e+05 &  8.10e-02 &  1.32e-02 \cr 
.01 & 1.17 & 1331 & 6.33e+05 &  4.08e-02 &  8.06e-03 \cr 
.1 & 1.17 & 1331 & 2.09e+06 &  4.59e-02 &  7.83e-03 \cr 
1 & 1.15 & 1331 & 1.83e+07 &  1.26e-01 &  1.57e-02 \cr 
}}}

B) nb = 5000

\centerline{\vbox{
\offinterlineskip
\halign{\strut \hfil # &  \hfil #  & \hfil # \hfil  &  ~~# &  ~~# & ~~#  & ~~#\cr
$\lambda$ & time & nc& CN  & emax  & rms \cr
\noalign{\hrule}
\noalign{\vskip 1 pt}
.001 & 1.81 & 1331 & 6.61e+05 &  5.96e-02 &  1.06e-02 \cr 
.01 & 1.79 & 1331 & 6.26e+05 &  3.94e-02 &  7.71e-03 \cr 
.1 & 1.81 & 1331 & 3.27e+06 &  5.81e-02 &  8.69e-03 \cr 
1 & 1.73 & 1331 & 2.74e+07 &  1.59e-01 &  1.86e-02 \cr 
}}}

C) nb = 10000

\centerline{\vbox{
\offinterlineskip
\halign{\strut \hfil # &  \hfil #  & \hfil # \hfil  &  ~~# &  ~~# & ~~#  & ~~#\cr
$\lambda$ & time & nc& CN  & emax  & rms \cr
\noalign{\hrule}
\noalign{\vskip 1 pt}
.001 & 3.03 & 1331 & 6.55e+05 &  5.18e-02 &  9.60e-03 \cr 
.01 & 2.95 & 1331 & 6.23e+05 &  3.85e-02 &  7.45e-03 \cr 
.1 & 2.99 & 1331 & 4.58e+06 &  6.89e-02 &  9.49e-03 \cr 
1 & 2.98 & 1331 & 3.69e+07 &  1.69e-01 &  1.91e-02 \cr 
}}}
\ms
\noin Once again the errors and condition numbers go down and then up. \eop
\ms

These examples suggest that in 
using IPBF or IPBC
with \tp{} splines good default values are {\tt nb = 2000} and
$\lambda = .01$ which is what we use for most of the examples in both the paper
and the supplement.

We now give two examples using splines on tetrahedral partitions of the
bounding box. We again see that increasing nb doesn't change the errors much,
and that the accuracy and condition numbers both go down and then up as we
increase $\lambda$.

\begin{ex}\label{SM24}   Repeat Example~\ref{SM21} but 
now using IPBF with the space $\cS^0_5(\tri)$ defined on a type-5
tetrahedral partition of the bounding box corresponding to a 
$7 \times 7 \times 7$ grid.
\end{ex}

\dis   We use the script \eword{T3flam} with  \eword{sphere50.dat},
\eword{pde = 4}, \eword{nb = 2000}, \eword{dx = 5}, \eword{mx = 7},
and set $\lambda_s = \lambda$.
The following table shows the same information as in the previous examples.

A) nb = 2000

\centerline{\vbox{
\offinterlineskip
\halign{\strut \hfil # &  \hfil #  & \hfil # \hfil  &  ~~# &  ~~# & ~~#  & ~~#\cr
$\lambda$ & time & nc& CN  & emax  & rms \cr
\noalign{\hrule}
\noalign{\vskip 1 pt}
 .001 & 8.96 & 25471 & 8.80e+04 & 1.14e-05 & 1.35e-06 \cr 
 .01 & 6.81 & 25471 & 2.34e+04 & 3.29e-06 & 6.98e-07 \cr 
 .1 & 6.69 & 25471 & 5.29e+04 & 2.76e-06 & 5.55e-07 \cr 
 1 & 6.21 & 25471 & 4.78e+05 & 5.70e-06 & 7.62e-07 \cr 
}}}

B) nb = 5000

\centerline{\vbox{
\offinterlineskip
\halign{\strut \hfil # &  \hfil #  & \hfil # \hfil  &  ~~# &  ~~# & ~~#  & ~~#\cr
$\lambda$ & time & nc& CN  & emax  & rms \cr
\noalign{\hrule}
\noalign{\vskip 1 pt}
 .001 & 7.55 & 25471 & 7.53e+04 & 6.91e-06 & 1.11e-06 \cr 
 .01 & 7.10 & 25471 & 2.01e+04 & 2.66e-06 & 6.16e-07 \cr 
 .1 & 7.39 & 25471 & 5.22e+04 & 2.75e-06 & 5.47e-07 \cr 
 1 & 7.04 & 25471 & 4.71e+05 & 6.85e-06 & 1.06e-06 \cr 
}}}

C) nb = 10000

\centerline{\vbox{
\offinterlineskip
\halign{\strut \hfil # &  \hfil #  & \hfil # \hfil  &  ~~# &  ~~# & ~~#  & ~~#\cr
$\lambda$ & time & nc& CN  & emax  & rms \cr
\noalign{\hrule}
\noalign{\vskip 1 pt}
 .001 & 7.71 & 25471 & 6.89e+04 & 5.40e-06 & 9.95e-07 \cr 
 .01 & 7.50 & 25471 & 1.86e+04 & 2.24e-06 & 5.75e-07 \cr 
 .1 & 8.31 & 25471 & 5.28e+04 & 2.76e-06 & 5.49e-07 \cr 
 1 & 7.81 & 25471 & 4.75e+05 & 1.06e-05 & 1.43e-06 \cr 
}}}

\begin{ex}\label{SM25}  Repeat Example~\ref{SM23} but now using IPBC 
with splines in the space $\cS^0_5(\tri)$ defined on a type-5 tetrahedral
partition of the bounding box for the domain corresponding to a
$7 \times 7 \times 7$ grid covering the bounding box.
\end{ex}

\dis   We use the script \eword{T3clam} with  \eword{bunny.ply},
\eword{pde = 6}, \eword{nb = 2000}, \eword{dx = 5}, and \eword{mx = 7}.
Here the errors are measured on 56,948 points inside the domain.

A) nb = 2000

\centerline{\vbox{
\offinterlineskip
\halign{\strut \hfil # &  \hfil #  & \hfil # \hfil  &  ~~# &  ~~# & ~~#  & ~~#\cr
$\lambda$ & time & nc& CN  & emax  & rms \cr
\noalign{\hrule}
\noalign{\vskip 1 pt}
 .001 &  2.31 & 25471 & 4.73e+07 & 3.19e-03 & 6.39e-04 \cr 
 .01 &  2.31 & 25471 & 5.72e+06 & 3.19e-03 & 6.39e-04 \cr 
 .1 &  2.51 & 25471 & 8.70e+05 & 3.17e-03 & 6.34e-04 \cr 
 1 &  2.33 & 25471 & 2.31e+05 & 3.01e-03 & 5.82e-04 \cr 
}}}

B) nb = 5000

\centerline{\vbox{
\offinterlineskip
\halign{\strut \hfil # &  \hfil #  & \hfil # \hfil  &  ~~# &  ~~# & ~~#  & ~~#\cr
$\lambda$ & time & nc& CN  & emax  & rms \cr
\noalign{\hrule}
\noalign{\vskip 1 pt}
 .001 &  2.54 & 25471 & 4.09e+07 & 3.04e-03 & 5.86e-04 \cr 
 .01 &  2.58 & 25471 & 4.97e+06 & 3.04e-03 & 5.86e-04 \cr 
 .1 &  2.50 & 25471 & 7.54e+05 & 3.01e-03 & 5.83e-04 \cr 
 1 &  2.53 & 25471 & 2.09e+05 & 2.84e-03 & 5.31e-04 \cr 
}}}

C) nb = 10000

\centerline{\vbox{
\offinterlineskip
\halign{\strut \hfil # &  \hfil #  & \hfil # \hfil  &  ~~# &  ~~# & ~~#  & ~~#\cr
$\lambda$ & time & nc& CN  & emax  & rms \cr
\noalign{\hrule}
\noalign{\vskip 1 pt}
 .001 &  2.77 & 25471 & 3.76e+07 & 2.87e-03 & 5.35e-04 \cr 
 .01 &  2.81 & 25471 & 4.57e+06 & 2.87e-03 & 5.37e-04 \cr 
 .1 &  2.81 & 25471 & 6.93e+05 & 2.85e-03 & 5.36e-04 \cr 
 1 &  2.85 & 25471 & 1.98e+05 & 2.73e-03 & 4.92e-04 \cr 
}}}
\ms

This set of examples explains why we use $n_b$ between 
2000 and 5000 
in all of the examples in the paper and the supplement, and why
we use $\lambda = .01$ as a default value.

\section{Choice of collocation points for the IPBC method}
To make use of the  IPBC method discussed in Section~2.2,
the user has
to select a set of collocation points. In this section we give examples
to illustrate how the accuracies depend on the number of points selected.
First we deal with the use of  trivariate \tp{} splines.
As discussed in Section~4.2, in this case
we take the collocation points  to be
the union of sets of $m_c^3$ points on equally-spaced grids covering each
sub-box of the partition.

\begin{ex}\label{SM31}    Solve Problem~\ref{Prob1} with the Laplace operator
on the  sphere shown in \frl{Fig1},
where the true solution is $\sin(5x)\sin(5y)\sin(5z)$.
Use the IPBC method with tensor-product splines of degree (5,5,5)
defined on $m\times m\times m$ grids covering the bounding box of $\Omega$,
and explore the errors for $m_c = 2,\dots,5$.
\end{ex}

\dis   We use the script \eword{TP3ipbcmc} with  the file \eword{sphere50.dat}, 
and set \eword{pde = 5},  \eword{nb = 2000},   and \eword{dx = 5}. 
The errors here are calculated on 103,283 points.

\noin A) $m = 6$
\ms

\centerline{\vbox{
\offinterlineskip
\halign{\strut \hfil # \hfil &  \hfil #  & \hfil # \hfil  &  \hfil # \hfil  &  ~~# & ~~#  & ~~#\cr
$m_c$ & setup & solve & nc& CN  & emax  & rms \cr
\noalign{\hrule}
\noalign{\vskip 1 pt}
2 & 0.19 & 0.20 & 1000 & 2.00e+09 &  7.96e-03 &  1.33e-03 \cr 
3 & 0.33 & 0.32 & 1000 & 2.90e+04 &  6.65e-04 &  1.09e-04 \cr 
4 & 0.69 & 0.63 & 1000 & 7.90e+03 &  5.23e-04 &  8.26e-05 \cr 
5 & 1.50 & 1.18 & 1000 & 8.63e+03 &  4.43e-04 &  7.66e-05 \cr 
}}}

\noin B) $m = 7$

\centerline{\vbox{
\offinterlineskip
\halign{\strut \hfil # \hfil &  \hfil #  & \hfil # \hfil  &  \hfil # \hfil  &  ~~# & ~~#  & ~~#\cr
$m_c$ & setup & solve & nc& CN  & emax  & rms \cr
\noalign{\hrule}
\noalign{\vskip 1 pt}
2 & 0.26 & 0.25 & 1331 & 1.06e+10 &  2.12e-03 &  2.86e-04 \cr 
3 & 0.63 & 0.50 & 1331 & 4.03e+04 &  2.36e-04 &  4.78e-05 \cr 
4 & 1.58 & 1.10 & 1331 & 1.23e+04 &  1.72e-04 &  3.41e-05 \cr 
5 & 4.13 & 2.26 & 1331 & 1.30e+04 &  1.42e-04 &  2.81e-05 \cr 
}}}
\ms\noindent
The value $m_c = 2$ gives the worst results.
As we increase $m_c$,  the computational times increase and  the condition numbers decrease somewhat,
but the errors do not change much.  \eop

As discussed in Section~6.2, when carrying out the IPBC method with trivariate
polynomial splines on tetrahedral partitions, we select the set $\Gamma$ of
collocation points to be the union of the sets of domain points of degree
$m_c$ defined on each of the tetrahedra in the partition.  
Here is an example to illustrate what happens when we use different values of  $m_c$.

\begin{ex}\label{SM32}    Solve Problem~\ref{Prob1} with the Laplace operator
on the  sphere shown in \frl{Fig1},
where 
the true solution is $\sin(5x)\sin(5y)\sin(5z)$.
Use the IPBC method with trivariate  splines of degree $d = 5$
defined on type-5 partitions associated with
$m\times m\times m$ grids covering the bounding box of $\Omega$.
Explore the errors for $m_c = 2,\dots,5$.

\end{ex}

 \dis We use the script \eword{T3ipbcmc} with  the file \eword{sphere50.dat}, 
\eword{pde = 5},  \eword{nb = 2000},  
and \eword{d = 5}. 

\noin {m = 5}
\ms

\centerline{\vbox{
\offinterlineskip
\halign{\strut \hfil # & \hfil # &  \hfil #  & ~~#  &  ~~# &  ~~# & ~~#  & ~~#\cr
$m_c$ & setup & solve & nc& CN  & emax  & rms \cr
\noalign{\hrule}
\noalign{\vskip 1 pt}
 2 & 0.12 & 0.25 &  7981 & Inf & 2.37e-02 & 1.80e-03 \cr 
 3 & 0.15 & 0.22 &  7981 & 2.75e+05 & 7.31e-03 & 1.72e-03 \cr 
 4 & 0.22 & 0.25 &  7981 & 3.22e+05 & 7.12e-03 & 1.64e-03 \cr 
 5 & 0.33 & 0.31 &  7981 & 3.76e+05 & 6.96e-03 & 1.60e-03 \cr 
}}}

\noin {m = 7}

\centerline{\vbox{
\offinterlineskip
\halign{\strut \hfil # & \hfil # &  \hfil #  & ~~#  &  ~~# &  ~~# & ~~#  & ~~#\cr
$m_c$ & setup & solve & nc& CN  & emax  & rms \cr
\noalign{\hrule}
\noalign{\vskip 1 pt}
 2 & 0.51 & 1.57 &  25471 & Inf & 1.21e-01 & 5.04e-03 \cr 
 3 & 0.90 & 1.45 &  25471 & 5.46e+05 & 1.41e-03 & 2.05e-04 \cr 
 4 & 1.45 & 1.40 &  25471 & 6.25e+05 & 1.96e-03 & 2.47e-04 \cr 
 5 & 2.48 & 1.69 &  25471 & 7.14e+05 & 2.36e-03 & 2.93e-04 \cr 
}}}
\ms\noindent
The choice $m_c = 2$  leads to rank deficient matrices in both cases.
As we increase $m_c$ to three and beyond, we see that
the setup times increase, but  the solve times,
condition numbers, and accuracy do not change much. \eop

Experiments on other PDEs and domains show a similar
behavior, and so in the interest of keeping computational time to a minimum, we 
suggest using
$m_c  = 3$ or $m_c = 4$ when working with either \tp{} spline spaces
or splines on tetrahedral partitions.

\section{Choice of the parameter \texorpdfstring{$\lambda_s$}{lambdas} when using \texorpdfstring{$C^0$}{C0} splines}
In Section~6 we asserted that when using IPBM methods with
splines that are only continuous, to get good results we need  to include a term
in the least-squares problem to force the spline to be approximately
$C^1$. This term includes a parameter $\lambda_s$. In this section we
explore how the choice of this parameter affects the behavior of the methods.

\begin{ex}\label{SM41}    
Solve Problem~\ref{Prob1} with the Laplace operator $L$ on the
sphere shown in \frl{Fig1},
where 
the true solution is $\sin(5x)\sin(5y)\sin(5z)$.
Use the IPBF method with splines of degree 5 
defined on a type-5 tetrahedral partition corresponding to
a $5\times 5\times 5$ grid
on the bounding box of $\Omega$.
Explore the errors  for a sequence of $\lambda_s$ values.
\end{ex}
    
\dis   We use the script \eword{T3flams} with  \eword{sphere50.dat},
\eword{pde = 5}, \eword{nb = 2000}, \eword{d = 5}, and \eword{mx = 5}.
The following table includes a column showing the  $\|Sc\|_\infty$ which
is a measure of smoothness.

\centerline{\vbox{
\offinterlineskip
\halign{\strut \hfil # &  \hfil # \hfil   & \hfil # \hfil  &  \hfil # \hfil
& \hfil  #\hfil  &  \hfil #\hfil  & \hfil #\hfil  & \hfil # \hfil\cr
$\lambda_s$ & time & nc &  CN  & $\|Sc\|$  & emax  & rms \cr
\noalign{\hrule}
\noalign{\vskip 1 pt}
0 & 0.72 & 7981 & 1.98e+11 & 8.23e+10 & 3.83e+09 & 1.90e+08 \cr 
 .001 & 0.87 & 7981 & 6.98e+04 & 1.39e-01 & 1.25e-03 & 1.98e-04 \cr 
 .01 & 0.87 & 7981 & 2.57e+04 & 2.99e-02 & 5.63e-04 & 9.83e-05 \cr 
 .1 & 0.82 & 7981 & 1.97e+04 & 2.67e-03 & 8.38e-04 & 1.43e-04 \cr 
 1 & 0.87 & 7981 & 8.20e+04 & 3.43e-05 & 8.96e-04 & 1.51e-04 \cr 
}}}
\sm\noin The first row of this table shows that with $\lambda_s = 0$ we
get very bad results. For this example $\lambda_s = .01$
gives the minimal RMS error but not the minimal condition number. \eop

\begin{ex}\label{SM42} 
Solve Problem~\ref{Prob1} with the Laplace operator $L$ on the
domain shown in \frl{Fig2}, where 
the true solution is $\sin(8x)+\sin(12y)+\sin(14z)$.
Use the IPBF method with the space $\cS^0_5(\tri)$
defined on a type-5 tetrahedral partition  associated with
a $5\times 5\times 5$ grid covering the bounding box of $\Omega$.
Explore the errors  for a sequence of $\lambda_s$ values.
\end{ex}
    
\dis   We use the script \eword{T3flams} with  \eword{ForearmLink.stl},
\eword{pde = 6}, \eword{nb = 2000}, \eword{dx = 5}, and \eword{mx = 5}.
\ms

\centerline{\vbox{
\offinterlineskip
\halign{\strut \hfil # &  \hfil # \hfil   & \hfil # \hfil  &  \hfil # \hfil
& \hfil  #\hfil  &  \hfil #\hfil  & \hfil #\hfil  & \hfil # \hfil\cr
$\lambda_s$ & time & nc &  CN  & $\|Sc\|$  & emax  & rms \cr
\noalign{\hrule}
\noalign{\vskip 1 pt} 
 0 & 0.72 & 7981 & 1.25e+12 & 3.34e+10 & 4.98e+04 & 1.89e+03 \cr 
 .001 & 1.03 & 7981 & 1.84e+05 & 7.93e-02 & 2.80e-03 & 3.70e-04 \cr 
 .01 & 0.83 & 7981 & 5.35e+04 & 1.99e-02 & 6.33e-04 & 1.13e-04 \cr 
 .1 & 0.80 & 7981 & 3.71e+04 & 2.20e-03 & 7.54e-04 & 1.02e-04 \cr 
 1 & 0.78 & 7981 & 1.27e+05 & 3.19e-05 & 8.16e-04 & 1.08e-04 \cr
}}}

\begin{ex}\label{SM43} Repeat Example~\ref{SM41} using
the IPBC method instead of the IPBF method.
\end{ex}

\dis   We use the script \eword{T3clams} with  \eword{sphere50.dat},
\eword{pde = 5}, \eword{nb = 2000}, \eword{dx = 5}, and \eword{mx = 5}.

\centerline{\vbox{
\offinterlineskip
\halign{\strut \hfil # &  \hfil #  & \hfil # \hfil  &  ~~# &  ~~# & ~~#  & ~~#\cr
$\lambda_s$ & time & nc &  CN  & $\|Sc\|$  & emax  & rms \cr
\noalign{\hrule}
\noalign{\vskip 1 pt}
 0 & 0.29 & 7981 & 9.00e+10 & 1.62e+05 & 6.91e+03 & 2.89e+02 \cr 
 .001 & 0.38 & 7981 & 1.88e+06 & 4.92e-01 & 7.37e-03 & 1.04e-03 \cr 
 .01 & 0.40 & 7981 & 2.74e+05 & 3.67e-01 & 7.10e-03 & 1.73e-03 \cr 
 .1 & 0.38 & 7981 & 1.28e+05 & 1.12e-01 & 1.23e-02 & 2.53e-03 \cr 
 1 & 0.39 & 7981 & 1.03e+05 & 2.86e-02 & 4.15e-03 & 1.02e-03 \cr 
}}}
\ms

In view of these examples we have opted to use $\lambda_s = .01$
as a default value in most of the examples involving $C^0$ splines,
both in the paper and in the supplement.

\section{Patch tests}
In this section we show that with appropriate choice of parameters,
both the IPBF and IPBC methods satisfy 
the patch test when using spline spaces.
In the first two examples we use the IPBF method with
\tp{} splines.

\begin{ex}\label{SM51}    Solve Problem~\ref{Prob1} with the Laplace operator on the 
sphere shown in \frl{Fig1}, 
where $f$ and $g$ are chosen so that
the true solution is $u = xy^2z^3$.
Use the IPBF method with tensor-product splines of degree (1,2,3)
defined on $m\times m \times m$ 
grids covering the bounding box of $\Omega$.
\end{ex}

\dis    We use the script \eword{TP3ipbfL} with  the file \eword{sphere50.dat}, 
\eword{pde = 9},  \eword{nb = 2000},   and \eword{(dx,dy,dz) = (1,2,3)}. 

\ms
\centerline{\vbox{
\offinterlineskip
\halign{\strut \hfil # &  \hfil #  & \hfil # \hfil  &  \hfil # \hfil  &  ~~# & ~~#  & ~~#\cr
mx & setup & solve & nc& CN  & emax  & rms \cr
\noalign{\hrule}
\noalign{\vskip 1 pt}
5 & 0.02 & 0.01 & 210 & 8.71e+02 &  4.16e-16 &  5.98e-17 \cr 
6 & 0.02 & 0.01 & 336 & 1.35e+03 &  4.72e-16 &  5.96e-17 \cr 
}}}
\ms

\begin{ex}\label{SM52}    Solve Problem~\ref{Prob1} with the Laplace operator
on the domain shown in \frl{Fig2},
where $f$ and $g$ are chosen so that
the true solution is $u = x^5 + 2y^5 + 3z^5$.
Use the IPBF method with tensor-product splines of degree (5,5,5) 
defined on $m\times m \times m$ 
grids covering the bounding box of $\Omega$.
\end{ex}

\dis    We use the script \eword{TP3ipbfL} with  the file \eword{ForearmLink.stl},
\eword{pde = 3}, \eword{nb = 2000}, and \eword{(dx,dy,dz) = (5,5,5)}. 

\ms
\centerline{\vbox{
\offinterlineskip
\halign{\strut \hfil # &  \hfil #  & \hfil # \hfil  &  \hfil # \hfil  &  ~~# & ~~#  & ~~#\cr
m & setup & solve & nc& CN  & emax  & rms \cr
\noalign{\hrule}
\noalign{\vskip 1 pt} 
5 & 0.25 & 0.35 & 729 & 4.13e+05 &  2.89e-15 &  3.73e-16 \cr 
6 & 0.50 & 0.50 & 1000 & 5.48e+05 &  3.44e-15 &  3.69e-16 \cr 
}}}
\sm

In the next two examples we use the IPBC method with \tp{} splines.

\begin{ex}\label{SM53}    Repeat Example~\ref{SM51} using
the IPBC method.
\end{ex}

\dis    We use the script \eword{TP3ipbcL} with  the file \eword{sphere50.dat},  
\eword{pde = 9},  \eword{nb = 2000}, $m_c = 3$,  and \eword{(dx,dy,dz) = (1,2,3)}.

\ms
\centerline{\vbox{
\offinterlineskip
\halign{\strut \hfil # &  \hfil #  & \hfil # \hfil  &  \hfil # \hfil  &  ~~# & ~~#  & ~~#\cr
mx & setup & solve & nc& CN  & emax  & rms \cr
\noalign{\hrule}
\noalign{\vskip 1 pt}
5 & 0.01 & 0.01 & 210 & 1.09e+02 &  2.22e-16 &  2.00e-17 \cr 
6 & 0.02 & 0.01 & 336 & 1.07e+02 &  2.78e-16 &  2.21e-17 \cr 
}}}

\begin{ex}\label{SM54}    Repeat Example~\ref{SM52} using the IPBC method.
\end{ex}

\dis    We use the script \eword{TP3ipbcL} with  the file \eword{ForearmLink.stl},
\eword{pde = 3}, \eword{nb = 2000}, $m_c = 3$, and \eword{(dx,dy,dz) = (5,5,5)}. 

\ms
\centerline{\vbox{
\offinterlineskip
\halign{\strut \hfil # &  \hfil #  & \hfil # \hfil  &  \hfil # \hfil  &  ~~# & ~~#  & ~~#\cr
m & setup & solve & nc& CN  & emax  & rms \cr
\noalign{\hrule}
\noalign{\vskip 1 pt} 
5 & 0.14 & 0.22 & 729 & 1.77e+04 &  2.22e-15 &  2.90e-16 \cr 
6 & 0.30 & 0.30 & 1000 & 3.85e+04 &  4.18e-15 &  5.68e-16 \cr 
}}}

In the following four  examples we approximate  with splines on 
tetrahedral partitions. These spline spaces have many more coefficients
than the ones in the previous four examples, but still show exact results.

\begin{ex}\label{SM55}    Repeat Example~\ref{SM51} using the IPBF method
with the spline spaces ${\cal S}^0_6(\tri)$  
on a type-5 partition  associated with $m\times m\times m$ grids
on the bounding box for $\Omega$.
\end{ex}

\dis  
We use the script \eword{T3ipbfL} with  the file \eword{sphere50.dat}, 
\eword{pde = 9},  \eword{nb = 2000},   and \eword{dx = 6}. 

\ms
\centerline{\vbox{
\offinterlineskip
\halign{\strut \hfil # \hfil & \hfil # &  \hfil #  & ~~#  &  ~~# &  ~~# & ~~#  & ~~#\cr
m & setup & solve & nc& CN  & emax  & rms \cr
\noalign{\hrule}
\noalign{\vskip 1 pt}
 4 & 0.28 & 0.30 &  5914 & 1.50e+04 & 3.89e-16 & 4.15e-17 \cr 
 5 & 1.14 & 0.84 &  13385 & 3.16e+04 & 6.66e-16 & 4.76e-17 \cr 
}}}

\begin{ex}\label{SM56} Repeat Example~\ref{SM52} 
using the IPBF method with splines in $\cS^0_5(\tri$)
defined on type-5 tetrahedral partitions of $m\times m \times m$ 
grids covering the bounding box of $\Omega$.
\end{ex}

\dis    We use the script \eword{T3ipbfL} with  the file \eword{ForearmLink.stl},
and set \eword{pde = 3},  \eword{nb = 2000},   and \eword{dx = 5}. 

\ms
\centerline{\vbox{
\offinterlineskip
\halign{\strut \hfil # & \hfil # &  \hfil #  & ~~#  &  ~~# &  ~~# & ~~#  & ~~#\cr
m & setup & solve & nc& CN  & emax  & rms \cr
\noalign{\hrule}
\noalign{\vskip 1 pt} 
 4 & 0.14 & 0.13 &  3556 & 6.25e+03 & 2.33e-15 & 1.90e-16 \cr 
 5 & 0.48 & 0.35 &  7981 & 1.32e+04 & 2.33e-15 & 2.14e-16 \cr 
}}}

\begin{ex}\label{SM57} Repeat Example~\ref{SM55}
using the IPBC method instead of the IPBF method.
Choose \eword{mc = 4}.
\end{ex}

\dis    We use the script \eword{T3ipbcL} with  the file \eword{sphere50.dat},
and set \eword{pde = 9},  \eword{nb = 2000},   \eword{mc = 4}, \eword{dx = 6},
and \eword{nx = 4}.

\ms
\centerline{\vbox{
\offinterlineskip
\halign{\strut \hfil # & \hfil # &  \hfil #  & ~~#  &  ~~# &  ~~# & ~~#  & ~~#\cr
m & setup & solve & nc& CN  & emax  & rms \cr
\noalign{\hrule}
\noalign{\vskip 1 pt}
 4 & 0.09 & 0.43 &  5914 & 5.39e+05 & 1.16e-14 & 1.71e-15 \cr 
 5 & 0.29 & 0.84 &  13385 & 1.06e+06 & 1.87e-14 & 2.18e-15 \cr 
}}}

\begin{ex}\label{SM58} Repeat Example~\ref{SM56}
using the IPBC method instead of the IPBF method.
\end{ex}

\dis    We use the script \eword{T3ipbcL} with  the file \eword{ForearmLink.stl},
and set \eword{pde = 3},  \eword{nb = 2000},  \eword{mc = 4}, 
\eword{dx = 5}, and \eword{nx = 4}.

\ms
\centerline{\vbox{
\offinterlineskip
\halign{\strut \hfil # & \hfil # &  \hfil #  & ~~#  &  ~~# &  ~~# & ~~#  & ~~#\cr
m & setup & solve & nc& CN  & emax  & rms \cr
\noalign{\hrule}
\noalign{\vskip 1 pt}
 4 & 0.07 & 0.10 &  3556 & 2.85e+05 & 6.15e-14 & 5.53e-15 \cr 
 5 & 0.21 & 0.25 &  7981 & 5.15e+05 & 3.70e-14 & 4.64e-15 \cr 
}}}

\section{Exactness}
In this section we give an example showing that when using
the IPBF and IPBC methods we get exact results when the true solution
lies in the spline space being used.

\begin{ex}\label{SM61}    Solve Problem~\ref{Prob1} with the Laplace operator on the
sphere shown in \frl{Fig1},  
where $f$ and $g$ are chosen so that
the true solution is $u = |x-.5|^3$.
We examine both tricubic tensor-product splines of degree dx = dy = dz = 3
with $m$ knots in all three directions,
and cubic splines on type-5 partitions corresponding to 
an $m\times m \times m$ grid covering the unit box.
\end{ex}

\dis    We use the scripts \eword{TP3exact} and \eword{T3exact}
 with the file \eword{sphere50.dat} and
\eword{nb = 2000}.

\ms
\noindent A) IPBF with TP splines

\ms
\centerline{\vbox{
\offinterlineskip
\halign{\strut \hfil # &  \hfil #  & \hfil # \hfil & ~~#   &  ~~# & ~~#  & ~~#\cr
m & setup & solve & nc & CN  & emax  & rms \cr
\noalign{\hrule}
\noalign{\vskip 1 pt}
5 & 0.03 & 0.03 & 343 & 5.69e+03 &  3.05e-16 &  4.45e-17 \cr 
}}} 

\ms
\noindent B)  IPBC with TP splines

\ms
\centerline{\vbox{
\offinterlineskip
\halign{\strut \hfil # &  \hfil #  & \hfil # \hfil & ~~#   &  ~~# & ~~#  & ~~#\cr
m & setup & solve & nc & CN  & emax  & rms \cr
\noalign{\hrule}
\noalign{\vskip 1 pt}
5 & 0.05 & 0.03 & 343 & 2.10e+02 &  1.39e-16 &  1.81e-17 \cr 
}}} 
\ms
\noindent C) IPBF with cubic splines

\ms
\centerline{\vbox{
\offinterlineskip
\halign{\strut \hfil # &  \hfil #  & \hfil # \hfil & ~~#    &  ~~# & ~~#  & ~~#\cr
m & setup & solve & nc  & CN & emax  & rms \cr
\noalign{\hrule}
\noalign{\vskip 1 pt}
5 & 0.12 & 0.04 & 1941 & 7.34e+02 &  1.25e-16 &  1.62e-17 \cr 
}}} 

\ms
\noindent D)  IPBC with cubic splines

\ms
\centerline{\vbox{
\offinterlineskip
\halign{\strut \hfil # &  \hfil #  & \hfil # \hfil & ~~#    &  ~~# & ~~#  & ~~#\cr
m & setup & solve & nc  & CN & emax  & rms \cr
\noalign{\hrule}
\noalign{\vskip 1 pt}
5 & 0.12 & 0.03 & 1941 & 2.55e+04 &  7.84e-16 &  1.23e-16 \cr 
}}}

\section{Nonelliptic PDEs}
In this section we give an example to illustrate the fact that
IPBM methods can produce good results
even when the operator $L$ is not elliptic.
We hold the domain and
the true solution fixed, and solve the BVP three times with different differential 
operators $L$ that are 1) elliptic, 2) weakly elliptic, and 3) non-elliptic.

\begin{ex}\label{SM71}    
Suppose  $L$ is defined as in (\ref{L1})  with coefficient vector given by
$a = (a_1,\ldots a_6)$. 
Solve Problem~\ref{Prob1} on the domain
shown in \frr{Fig1} where $f$ and $g$ are chosen so that
the true solution is $\sin(5x)\sin(5y)\sin(5z)$.
Use the IPBF method with \tp{} splines of degree (6,6,6) defined on
 $m\times m \times m$   grids covering the  bounding box of $\Omega$.
\end{ex}

 \dis For part A) We use the script \eword{TP3ipbfr} with  the file \eword{nut.stl}, 
\eword{pde = 5},  \eword{nb = 2000}, 
\eword{dx = 6}, and \eword{mx = 5}.  For parts B) and C) we use \eword{pde = 13} and
\eword{pde = 14}, respectively.
Here the errors are measured on 87,025 points.
\ms

\noin A) $a = (1,1,1,0,0,0)$
\sm

\centerline{\vbox{
\offinterlineskip
\halign{\strut \hfil # &  \hfil #  & \hfil # \hfil  &  \hfil # \hfil  &  ~~# & ~~#  & ~~#\cr
m & setup & solve & nc& CN  & emax  & rms \cr
\noalign{\hrule}
\noalign{\vskip 1 pt}
5 & 0.62 & 0.91 & 1000 & 2.10e+06 &  8.15e-05 &  1.30e-05 \cr 
6 & 1.25 & 1.38 & 1331 & 2.28e+06 &  1.43e-05 &  2.59e-06 \cr 
7 & 2.27 & 2.20 & 1728 & 3.01e+06 &  3.41e-06 &  6.85e-07 \cr 
8 & 3.93 & 0.41 & 2197 & 8.40e+06 &  1.34e-06 &  2.26e-07 \cr 
}}
\hskip 2pc
\vbox{\offinterlineskip
\halign{\strut
\hfil  \hfil #  &  ~~#&   ~~#   \cr
\omit\span \hfil rates \hfil \cr
\noalign{\vskip 4 pt}
\noalign{\hrule}
\noalign{\vskip 1 pt}
7.80 & 7.22 \cr 
 7.86 & 7.30 \cr 
 6.08 & 7.19 \cr 
}}}

\noin B) $a= (1,1,1,1,1,1)$
\sm

\centerline{\vbox{
\offinterlineskip
\halign{\strut \hfil # &  \hfil #  & \hfil # \hfil  &  \hfil # \hfil  &  ~~# & ~~#  & ~~#\cr
m & setup & solve & nc& CN  & emax  & rms \cr
\noalign{\hrule}
\noalign{\vskip 1 pt}
5 & 0.59 & 0.88 & 1000 & 1.82e+06 &  7.23e-05 &  1.16e-05 \cr 
6 & 1.26 & 1.37 & 1331 & 2.08e+06 &  1.53e-05 &  2.72e-06 \cr 
7 & 2.25 & 2.15 & 1728 & 3.21e+06 &  4.49e-06 &  7.50e-07 \cr 
8 & 4.03 & 0.41 & 2197 & 8.53e+06 &  2.09e-06 &  3.27e-07 \cr 
}}
\hskip 2pc
\vbox{\offinterlineskip
\halign{\strut
\hfil  \hfil #  &  ~~#&   ~~#   \cr
\omit\span \hfil rates \hfil \cr
\noalign{\vskip 4 pt}
\noalign{\hrule}
\noalign{\vskip 1 pt}
6.96 & 6.49 \cr 
 6.72 & 7.07 \cr 
 4.95 & 5.39 \cr 
}}}
\ms\noin

\noin C) $a= (1,1,1,10,10,10)$
\sm

\centerline{\vbox{
\offinterlineskip
\halign{\strut \hfil # & \hfil # &  \hfil #  & ~~#  &  ~~# &  ~~# & ~~#  & ~~#\cr
m & setup & solve & nc& CN  & emax  & rms \cr
\noalign{\hrule}
\noalign{\vskip 1 pt}
5 & 0.60 & 0.90 & 1000 & 3.40e+06 &  1.07e-04 &  1.57e-05 \cr 
6 & 1.20 & 1.37 & 1331 & 3.14e+06 &  2.52e-05 &  3.90e-06 \cr 
7 & 2.21 & 2.16 & 1728 & 4.72e+06 &  3.96e-06 &  6.36e-07 \cr 
8 & 3.96 & 0.42 & 2197 & 9.93e+06 &  1.64e-06 &  2.65e-07 \cr 
}}
\hskip 2pc
\vbox{\offinterlineskip
\halign{\strut
\hfil  \hfil #  &  ~~#&   ~~#   \cr
\omit\span \hfil rates \hfil \cr
\noalign{\vskip 4 pt}
\noalign{\hrule}
\noalign{\vskip 1 pt}
6.48 & 6.24 \cr 
 10.14 & 9.95 \cr 
 5.74 & 5.68 \cr 
}}}
\ms\noin
The eigenvalues of the ellipticity matrix
are (1,1,1) for case A),
(0,0,3) in case B), and (-9,-9,21) in case C). 
This means that the operator is elliptic in case A), weakly
elliptic in case B), and non-elliptic in case C).
The tables show very similar results for all three cases.  

\section{Generating points for computing errors}
In all of the examples in both the paper and the supplement 
we are reporting errors defined on a large dense
set of points covering the domain $\Omega$. We create these points by starting with
points on a grid covering the bounding box for $\Omega$,  and then removing 
points that lie outside of $\Omega$. Our function
{\tt inpts3} described  in Sect.~3.4 generates points in this way.
However,  for complicated domains with holes, it may produce a few points that
do not actually lie in $\Omega$. This means that the reported errors 
in some examples could be slightly higher than the actual errors. 
To illustrate that this makes very
little difference, consider the following problem where the domain is defined
by an STL file. For such domains we can use 
the MATLAB function {\tt inpolyhedron} (see [9]) instead 
of {\tt inpts3}.  It
 ensures that  no points outside of $\Omega$ are chosen.

\begin{ex}\label{inpts} Solve Problem~\ref{Prob1} with the Laplace operator 
on the domain shown in \frl{Fig2}
with $f$ and $g$ such that
the true solution is $\sin(5x)\sin(5y)\sin(5z)$.
Use IPBF with tensor-product splines of degree (4,4,4)
defined on $m\times m\times m$ grids covering the bounding box of $\Omega$ 
for a sequence of increasing values of m. Compute the errors
on interior points generated by {\tt inpts3}, and compare with
using the interior points generated by {\tt inpolyhedron}.
\end{ex}

\begin{figure} [H]
  \centering
  \includegraphics[width = .6\linewidth]{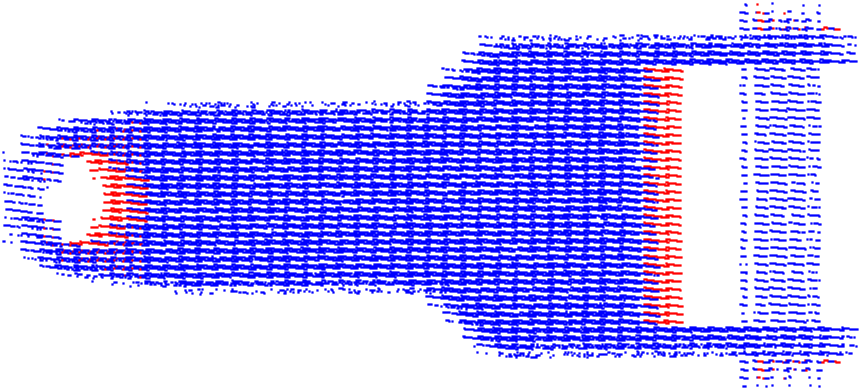}
  \caption{Sets of points inside the domain  of ForearmLink}
  \label{forearm_int}
\end{figure}

 \dis We use the script \eword{TP3ipbfLrint} with  the file \eword{ForearmLink.stl}, 
\eword{pde = 5},  \eword{nb = 5000}, 
\eword{dx = 4}, and \eword{mx = 5}.
In the first table, the errors are evaluated on 42,478 points generated 
with {\tt inpts3} using no additional sensors, along with $50,000$ boundary points. 
The second table finds the  errors on 40,185 interior points 
using {\tt inpolyhedron} in addition to 50,000 boundary points. 
Figure \ref{forearm_int} shows the second set of interior points in blue.
The  2,293 points in the first set that lie outside of the domain 
are shown in red.
\ms

\noin A) inpts3
\sm

\centerline{\vbox{
\offinterlineskip
\halign{\strut \hfil # & \hfil # &  \hfil #  & ~~#  &  ~~# &  ~~# & ~~#  & ~~#\cr
m & setup & solve & nc& CN  & emax  & rms \cr
\noalign{\hrule}
\noalign{\vskip 1 pt}
5 & 0.16 & 0.20 & 512 & 3.70e+04 &  9.65e-04 &  2.10e-04 \cr 
6 & 0.27 & 0.28 & 729 & 1.03e+05 &  3.34e-04 &  7.81e-05 \cr 
7 & 0.41 & 0.37 & 1000 & 1.99e+05 &  1.03e-04 &  2.64e-05 \cr 
8 & 0.66 & 0.50 & 1331 & 3.66e+05 &  5.11e-05 &  1.18e-05 \cr 
9 & 1.01 & 0.69 & 1728 & 5.69e+05 &  2.57e-05 &  5.71e-06 \cr 
10 & 1.53 & 0.93 & 2197 & 8.25e+05 &  1.46e-05 &  3.18e-06 \cr 
}}
\hskip 2pc
\vbox{\offinterlineskip
\halign{\strut
\hfil  \hfil #  &  ~~#&   ~~#   \cr
\omit\span \hfil rates \hfil \cr
\noalign{\vskip 4 pt}
\noalign{\hrule}
\noalign{\vskip 1 pt}
4.76 & 4.43 \cr 
 6.47 & 5.94 \cr 
 4.52 & 5.21 \cr 
 5.16 & 5.47 \cr 
 4.80 & 4.98 \cr 
}}}

\noin B) inpolyhedron
\sm

\centerline{\vbox{
\offinterlineskip
\halign{\strut \hfil # & \hfil # &  \hfil #  & ~~#  &  ~~# &  ~~# & ~~#  & ~~#\cr
m & setup & solve & nc& CN  & emax  & rms \cr
\noalign{\hrule}
\noalign{\vskip 1 pt}
 5 & 0.15 & 0.20 & 512 & 3.70e+04 &  9.65e-04 &  2.12e-04 \cr 
6 & 0.26 & 0.27 & 729 & 1.03e+05 &  3.34e-04 &  7.88e-05 \cr 
7 & 0.41 & 0.36 & 1000 & 1.99e+05 &  1.03e-04 &  2.67e-05 \cr 
8 & 0.65 & 0.60 & 1331 & 3.68e+05 &  5.11e-05 &  1.20e-05 \cr 
9 & 0.98 & 0.66 & 1728 & 5.70e+05 &  2.57e-05 &  5.74e-06 \cr 
10 & 1.50 & 0.96 & 2197 & 8.27e+05 &  1.46e-05 &  3.20e-06 \cr 
}}
\hskip 2pc
\vbox{\offinterlineskip
\halign{\strut
\hfil  \hfil #  &  ~~#&   ~~#   \cr
\omit\span \hfil rates \hfil \cr
\noalign{\vskip 4 pt}
\noalign{\hrule}
\noalign{\vskip 1 pt}
4.76 & 4.44 \cr 
 6.47 & 5.94 \cr 
 4.52 & 5.21 \cr 
 5.16 & 5.50 \cr 
 4.80 & 4.96 \cr 
}}}
\ms \noin
We see that the max errors in A) and B) are identical 
while the RMS errors differ only slightly. Since {\tt inpts3} can be used
with any of the types of inputs considered in our software package, 
we use it to generate interior points for all examples in this paper. 

\section{Refinement strategies}
When working with splines to solve a BVP we can attain higher accuracy
in three ways: 1) by holding the degree of the splines fixed and decreasing
the mesh size ($h$--refinement), 2) by holding the mesh size fixed
and increasing $d$ ($d$--refinement), or 3) by doing both.
Here is an example where these strategies are compared.
It is based on Example~\ref{Exsphere} where $h$ refinement was used

\begin{ex}\label{SM91}    
Solve Problem~\ref{Prob1} with the Laplace operator
on the sphere shown in \frl{Fig1}, 
with $f$ and $g$ such that
the true solution is $\sin(5x)\sin(5y)\sin(5z)$.
Use the IPBF with tensor-product splines of degree (d,d,d) 
defined on $m\times m\times m$ grids covering the bounding box of $\Omega$ 
for a sequence of increasing values of d and m.
\end{ex}

\dis We present three tables. The first two use $d$ refinement and
are  generated by using
\eword{TP3ipbfLd} with the file \eword{sphere50.dat},
\eword{pde = 5} and  \eword{nb = 2000}, first for $m = 4$ and then
with $m = 10$. The third table uses $h$ refinement, and
is generated with
\eword{TP3ipbfLr} with $d=5$ and starting with $m = 10$.
The errors here are computed
on a set of 103,283 points covering $\Omega$.
\ms

\noin A) $m=4$
\sm

\centerline{\vbox{
\offinterlineskip
\halign{\strut \hfil # &  \hfil #  & \hfil # \hfil  &  \hfil # \hfil  &  ~~# & ~~#  & ~~#\cr
d & setup & solve & nc& CN  & emax  & rms \cr
\noalign{\hrule}
\noalign{\vskip 1 pt}
2 & 0.01 & 0.01 & 125 & 4.15e+02 &  2.16e-01 &  5.05e-02 \cr 
3 & 0.03 & 0.03 & 216 & 1.92e+03 &  6.10e-02 &  1.21e-02 \cr 
4 & 0.05 & 0.09 & 343 & 2.19e+04 &  2.53e-02 &  3.05e-03 \cr 
5 & 0.13 & 0.24 & 512 & 2.11e+05 &  2.53e-03 &  3.92e-04 \cr 
6 & 0.32 & 0.70 & 729 & 2.04e+06 &  8.00e-04 &  1.00e-04 \cr 
7 & 0.62 & 1.62 & 1000 & 1.65e+07 &  1.07e-04 &  1.28e-05 \cr 
8 & 1.20 & 3.69 & 1331 & 1.82e+08 &  2.52e-05 &  3.12e-06 \cr 
9 & 1.97 & 9.89 & 1728 & 4.53e+09 &  3.04e-06 &  6.12e-07 \cr 
}}}

\noin B) $m=10$
\sm

\centerline{\vbox{
\offinterlineskip
\halign{\strut \hfil # &  \hfil #  & \hfil # \hfil  &  \hfil # \hfil  &  ~~# & ~~#  & ~~#\cr
d & setup & solve & nc& CN  & emax  & rms \cr
\noalign{\hrule}
\noalign{\vskip 1 pt}
2 & 0.22 & 0.21 & 1331 & 2.00e+03 &  2.75e-02 &  4.16e-03 \cr 
3 & 0.52 & 0.47 & 1728 & 2.40e+04 &  7.52e-04 &  1.43e-04 \cr 
4 & 1.77 & 0.35 & 2197 & 7.65e+05 &  3.01e-05 &  6.50e-06 \cr 
5 & 5.17 & 0.58 & 2744 & 3.44e+06 &  5.54e-06 &  6.93e-07 \cr 
6 & 13.96 & 1.16 & 3375 & 1.86e+07 &  7.77e-07 &  8.19e-08 \cr 
7 & 37.74 & 1.94 & 4096 & 1.40e+08 &  8.12e-08 &  1.19e-08 \cr 
8 & 82.67 & 2.60 & 4913 & 6.70e+08 &  7.62e-09 &  1.16e-09 \cr 
9 & 152.46 & 4.49 & 5832 & 1.80e+10 &  1.13e-09 &  1.27e-10 \cr 
}}}
\ms\noin 

\noin C) $d=5$
\sm

\centerline{\vbox{
\offinterlineskip
\halign{\strut \hfil # &  \hfil #  & \hfil # \hfil  &  \hfil # \hfil  &  ~~# & ~~#  & ~~#\cr
m & setup & solve & nc& CN  & emax  & rms \cr
\noalign{\hrule}
\noalign{\vskip 1 pt}
10 & 5.08 & 0.59 & 2744 & 3.40e+06 &  5.52e-06 &  6.94e-07 \cr 
11 & 8.74 & 0.83 & 3375 & 6.85e+06 &  2.63e-06 &  3.40e-07 \cr 
12 & 17.33 & 1.62 & 4096 & 1.40e+07 &  1.48e-06 &  1.78e-07 \cr 
13 & 28.63 & 2.00 & 4913 & 2.99e+07 &  7.58e-07 &  9.87e-08 \cr 
14 & 32.84 & 3.10 & 5832 & 6.33e+07 &  6.01e-07 &  6.44e-08 \cr 
15 & 70.38 & 4.32 & 6859 & 1.25e+08 &  2.88e-07 &  3.73e-08 \cr 
16 & 103.61 & 8.84 & 8000 & 2.15e+08 &  2.23e-07 &  2.60e-08 \cr 
17 & 144.98 & 9.82 & 9261 & 3.02e+08 &  1.39e-07 &  1.70e-08 \cr 
18 & 201.73 & 16.45 & 10648 & 3.29e+08 &  9.94e-08 &  1.21e-08 \cr 
19 & 337.67 & 18.89 & 12167 & 3.57e+08 &  6.12e-08 &  8.67e-09 \cr 
}}}
\ms\noin 
Table C) is just a continuation of table B in Example~\ref{Exsphere}.
The tables show that by using IPBF with higher degree splines 
we can get very high accuracy
with relatively small numbers of coefficients.
We can also use the tables to compare the different refinement
strategies.
For example, the ($d$,$m$) combination (9,10) gives better accuracy
(an RMS error of 1.27e-10) 
than the combination (5,19)  (an RMS error of 8.67e-09),
but with less than half the number of coefficients and half the computational
time. 

\nocite{soft}
\bibliographystyle{siamplain}
\bibliography{references}